\newif\iffullproof
	\fullprooffalse
\newcommand\proofswitch[2]{\iffullproof#1\else#2\fi}
\newif\ifarxiv
	\arxivtrue
\newcommand\siam[1]{\ifarxiv\relax\else\relax#1\relax\fi\relax}
\newcommand\arxiv[1]{\ifarxiv\relax#1\else\relax\fi\relax}
\newcommand\switch[2]{\siam{#1}\arxiv{#2}}
\ifarxiv
	\newcommand\newarxiv[1]{#1}
	\newcommand\newsiam[1]{}
\else
	\newcommand\newarxiv[1]{}
	\newcommand\newsiam[1]{#1}
\fi
\newcommand\Newswitch[2]{\newsiam{#1}\newarxiv{#2}}

\newcommand\TheTitle{Forward-backward envelope for the sum of two nonconvex functions: further properties and nonmonotone line-search algorithms}
\newcommand\TheShortTitle{Forward-backward envelope for the sum of two nonconvex functions}
\newcommand\TheAbbrTitle{FBE for the sum of two nonconvex functions}

\newcommand\TheAuthors{%
	A. Themelis,
	L. Stella and
	P. Patrinos%
}

\newcommand\TheAddressIMT{%
	IMT School for Advanced Studies Lucca,
	Piazza S. Francesco 19, 55100 Lucca, Italy%
}
\newcommand\TheAddressKU{%
	KU Leuven,
	Department of Electrical Engineering (ESAT-STADIUS) and Optimization in Engineering Center (OPTEC),
	Kasteelpark Arenberg 10, 3001 Leuven, Belgium%
}

\switch{
	\PassOptionsToPackage{usenames,dvipsnames,table}{xcolor}
	\documentclass[review]{siamart0516}
	\usepackage[%
		version=2017/02/21,
		xcolor=noload,
		amsthm=noload,
		abbr,
		appendix=true,
	]{clevethm}
	\usepackage[%
		version=2017/04/25,
	]{myPreamble}
}{
	\documentclass[a4paper, 10pt, USenglish]{amsart}
	\usepackage[%
		version=2017/03/05,
		abbr,
		appendix=true,
	]{clevethm}
	\usepackage[%
		version=2017/04/25,
	]{myPreamble}
	\usepackage[foot]{amsaddr}	
}
\usepackage{pifont}
\newcommand\cmark{\textcolor{ForestGreen}{\ding{51}}}
\newcommand\xmark{\textcolor{red}{\ding{55}}}

\definecolor{color_FBS}{rgb}{0.00000,0.44700,0.74100}%
\definecolor{color_IFBS}{rgb}{0.0,0.0,0.0}%
\definecolor{color_AFBS}{rgb}{0.49400,0.18400,0.55600}%
\definecolor{color_Broyden}{rgb}{0.85000,0.32500,0.09800}%
\definecolor{color_BFGS}{rgb}{0.92900,0.69400,0.12500}%
\newcommand\linestylefull{solid}
\newcommand\linestylelimited{dashed}

\pgfplotsset{compat=1.10}
\let\im\relax
\DeclareMathOperator*\im{Im}

\newcommand\minfbe{{\sf minFBE}}
\newcommand\zerofpr{{\sf ZeroFPR}}
\newcommand\funcf{\func{f}{\R^n}{\R}}
\newcommand\funcg{\func{g}{\R^n}{\Rinf}}
\newcommand\cont[1][]{C^{#1}}
\newcommandx\T[2][1=\gamma,2={}]{T_{#1}^{#2}}
\newcommandx\Res[2][1=\gamma,2={}]{R_{#1}^{#2}}
\newcommandx\lin[4][{1=\gamma,4={}}]{
	\ifstrempty{#3}{%
		\ifstrempty{#2}{%
			\ell_{#1}^{#4}%
		}{%
			\ell_{#1}^{#4}\left(#2\right)
		}%
	}{%
		\ell_{#1}^{#4}\left(#2;\,#3\right)%
	}%
}%
\DeclareMathOperator*\mesh{mesh}
\newcommand\refZ[1][]{\hyperref[alg:zerofpr]{\zerofpr}}

\AtBeginDocument{
	\let\OLDref\ref%
	\let\OLDcref\cref%
	\let\OLDCref\Cref%
	\renewcommand\ref[1]{%
					\def\myVarOne{alg:zerofpr}%
					\ifx\myVar\myVarOne%
						\refZ[]%
					\else%
						\OLDref{#1}%
					\fi%
	}%
	\renewcommand\cref[1]{%
					\def\myVarOne{alg:zerofpr}%
					\ifx\myVar\myVarOne%
						\refZ[]%
					\else%
						\OLDcref{#1}%
					\fi%
	}%
	\renewcommand\Cref[1]{%
					\def\myVarOne{alg:zerofpr}%
					\ifx\myVar\myVarOne%
						\refZ[]%
					\else%
						\OLDCref{#1}%
					\fi%
	}%
}%

\newarxiv{}

\switch{
	\headers{\scriptsize\TheShortTitle}{\TheAuthors}
	\title{%
		\TheTitle%
		\thanks{%
			Submitted to the editors June 15, 2016.%
			\funding{%
				The work of the third author was supported by KU Leuven Research Council BOF/STG-15-043.%
			}%
		}%
	}

	\author{%
		Andreas Themelis\hspace{0.5pt}\footnotemark[3]{\hspace{3.5pt}\textsuperscript,}%
		\thanks{%
			\TheAddressIMT{}
			(\email{\{andreas.themelis,lorenzo.stella\}@imtlucca.it}).%
		}
		\and
		Lorenzo Stella\footnotemark[2]{\hspace{3.5pt}\textsuperscript,}\footnotemark[3]
		\and
		Panagiotis Patrinos\thanks{%
		\TheAddressKU{}
			(+32/16374445 \email{panos.patrinos@esat.kuleuven.be}).%
		}
	}

	\headers{\scriptsize\TheShortTitle}{\TheAuthors}%
}{
	\author{%
		Andreas Themelis,\textsuperscript{1,2}
		Lorenzo Stella,\textsuperscript{1,2}
		Panagiotis Patrinos\textsuperscript{1}%
	}

	\address{%
		\rm \textsuperscript{1}\TheAddressKU\newline
		(panos.patrinos@esat.kuleuven.be).}
	\address{%
		\rm \textsuperscript{2}\TheAddressIMT\newline
		(\{andreas.themelis, lorenzo.stella\}@imtlucca.it).%
	}

	\title[\TheAbbrTitle]{\TheTitle}


	\date\today
}

\ifpdf
	\hypersetup{%
		pdftitle={\TheTitle},
		pdfauthor={\TheAuthors},
	}
\fi

\begin{document}

\siam{%
	\maketitle
}%

\begin{abstract}
	We propose \refZ[], a nonmonotone linesearch algorithm for minimizing the sum of two nonconvex functions, one of which is smooth and the other possibly nonsmooth.
	\refZ[] is the first algorithm that, despite being fit for fully nonconvex problems and requiring only the black-box oracle of forward-backward splitting (FBS) --- namely evaluations of the gradient of the smooth term and of the proximity operator of the nonsmooth one --- achieves superlinear convergence rates under mild assumptions at the limit point when the linesearch directions satisfy a Dennis-Mor\'e condition, and we show that this is the case for quasi-Newton directions.
	Our approach is based on the forward-backward envelope (FBE), an exact and strictly continuous penalty function for the original cost.
	Extending previous results we show that, despite being nonsmooth for fully nonconvex problems, the FBE still enjoys favorable first- and second-order properties which are key for the convergence results of \refZ[].
	Our theoretical results are backed up by promising numerical simulations.
	On large-scale problems, by computing linesearch directions using limited-memory quasi-Newton updates our algorithm greatly outperforms FBS and its accelerated variant (AFBS).
\end{abstract}

\arxiv{%
	\maketitle
}%

\siam{%

}%

\arxivfalse

\section{Introduction}
	In this paper we deal with optimization problems of the form
\begin{equation}\label{eq:Problem}
	\minimize_{x\in\R^n}\ \varphi(x) \equiv f(x) + g(x)
\end{equation}
under the following assumptions, which will be valid throughout the paper without further mention.
\begin{ass}[Basic assumption]\label{ass:fg}%
In problem \eqref{eq:Problem}
\begin{enumerate}
	\item\label{ass:f}%
		\(f\in\cont[1,1](\R^n)\) (differentiable with \(L_f\)-Lipschitz continuous gradient);
	\item\label{ass:g}%
		\(\funcg\) is proper, closed and \(\gamma_g\)-prox-bounded (see \Cref{sec:Notation});
	\item\label{ass:phi}%
		a solution exists, that is, \(\argmin\varphi\neq\emptyset\).
\end{enumerate}
\end{ass}
Both \(f\) and \(g\) are allowed to be nonconvex, making \eqref{eq:Problem} prototypic for a plethora of applications spanning signal and image processing, machine learning, statistics, control and system identification.
A well known algorithm addressing \eqref{eq:Problem} is \emph{forward-backward splitting} (FBS), also known as \emph{proximal gradient method}.
FBS has been thoroughly analyzed under the assumption of \(g\) being convex.
If moreover \(f\) is convex, then FBS is known to converge globally with rate \(O(1/k)\) in terms of objective value, where \(k\) is the iteration count. 
In this case, accelerated variants of FBS, also known as \emph{fast} forward-backward splitting (FFBS), can be derived thanks to the work of Nesterov
\switch{%
	\cite{beck2009fast,nesterov2013gradient},
}{
	\cite{nesterov1983method,tseng2008accelerated,beck2009fast,nesterov2013gradient},
}%
that only require minimal additional computations per iteration but achieve the provably optimal global convergence rate of order \(o(1/k^2)\) \cite{attouch2016rate}.

The work in \cite{patrinos2013proximal} pioneered an alternative acceleration technique.
The method is based on an exact, real-valued penalty function for the original problem \eqref{eq:Problem}, namely the \emph{forward-backward envelope} (FBE), defined as follows
\begin{equation}\label{eq:FBE1}
	\varphi_\gamma(x)
{}={}
	\varphi_\gamma^{f,g}(x)
{}\coloneqq{}
	\inf_{z\in\R^n}
	\set{f(x)+\tinnprod{\nabla f(x)}{z-x} + \tfrac{1}{2\gamma}\|z-x\|^2 + g(z)}
\end{equation}
where \(\gamma>0\) is a given parameter.
We will adopt the simpler notation \(\varphi_\gamma\) without superscript whenever \(f\) and \(g\) are clear from the context.

The name \emph{forward-backward envelope} comes from the fact that \(\varphi_\gamma(x)\) is the value of the minimization problem that defines the \emph{forward-backward} step and alludes to the kinship that it has with the Moreau \emph{envelope}.
These claims will be addressed more in detail in \Cref{sec:FBE}.
When \(f\) is sufficiently smooth and both \(f\) and \(g\) are convex, the FBE was shown to be continuously differentiable and amenable to be minimized with generalized Newton methods.
More recently, \cite{stella2017forward} proposed a linesearch algorithm based on (L-)BFGS quasi-Newton directions for minimizing the FBE.
The curvature information exploited by Newton-like methods acts as an online preconditioner, enabling superlinear rates of convergence under mild assumptions.
However, unlike plain (F)FBS schemes, such methods require accessing second-order information of the smooth term \(f\) (needed for the evaluation of \(\nabla\varphi_\gamma\)), and are well defined only as long as the nonsmooth term \(g\) is convex.
On the contrary, FBS merely requires first-order information on \(f\) and \emph{prox-boundedness} of the nonsmooth term \(g\), in which case all accumulation points are \emph{stationary} for \(\varphi\), \ie they satisfy the first order necessary conditions \cite{attouch2013convergence}.


\switch{%
	\paragraph*{Contributions}
}{%
	\subsection{Contributions}
}%
In this paper we propose \refZ[], a nonmonotone linesearch algorithm that, to the best of our knowledge, is the first that
(1) addresses the same range of problems as FBS,
(2) requires the same black-box oracle as FBS (gradient of one function and \emph{proximity operator} of the other),
(3) yet achieves superlinear rates under mild assumptions (only) at the limit point.
Though related to \minfbe{} algorithm \cite{stella2017forward}, \refZ[] is conceptually different, mainly because it is \emph{gradient-free}, in the sense that it does not require the gradient of the FBE.
Moreover,
\begin{itemize}[
	label={\hspace*{0.5\parindent}\clap\textbullet},
	itemsep=2pt,
	labelsep=0.5\parindent,
	leftmargin=*,
	itemindent=\parindent,
]
\item
	We provide the necessary theoretical background linking the concepts of stationarity of a point for problem \eqref{eq:Problem}, \emph{criticality} and optimality.
	To the best of our knowledge, such an analysis was previously made only for the proximal point algorithm \cite{poliquin1996prox-regular} and for a special case of the projected gradient method \cite{beck2016minimization}.
\item
	The analysis of the FBE, previously studied only in the case of \(f\) being \(C^2(\R^n)\) and \(g\) convex \cite{stella2017forward}, is extended to \(f\) and \(g\) as in \Cref{ass:fg}.
	In particular, we provide mild assumptions on \(f\) and \(g\) that ensure
	(1) continuous differentiabilty of the FBE around critical points,
	(2) (strict) twice differentiability at critical points,
	and (3) equivalence of strong local minimality for the original function and the FBE.
\item%
	Exploiting the investigated properties of the FBE and of critical points we prove that \refZ[] with monotone linesearch converges
	(1) globally if \(\varphi_\gamma\) has the \emph{Kurdyka-{\L}ojasiewicz} property \cite{lojasiewicz1963propriete,lojasiewicz1993geometrie,kurdyka1998gradients},
	and (2) superlinearly when quasi-Newton Broyden
	\arxiv{%
		or BFGS
	}%
	directions are employed, under mild additional requirements at the limit point.
\end{itemize}

\switch{%
	\paragraph*{Organization of the paper}
}{
	\subsection{Organization of the paper}
}%
In \Cref{sec:Preliminaries} we introduce some notation and list some known facts about FBS.
In \Cref{sec:StationaryCritical} we define and explore notions of stationarity and criticality for the investigated problem and relate them with properties of the forward-backward operator.
In \Cref{sec:FBE} we extend the results of \cite{stella2017forward} about the fundamental properties of the FBE to the more general setting addressed in this paper, where \(f\) and \(g\) satisfy \Cref{ass:fg}; for the sake of readability, some of the proofs are deferred to \Cref{sec:FBE_proof}.
\Cref{sec:Algorithm} addresses the core contribution of the paper, \refZ[]; although arbitrary directions can be chosen, we specialize the results on superlinear convergence to \switch{
	a quasi-Newton Broyden method
}{
	quasi-Newton Broyden and BFGS methods
}%
so as to truely maintain the same black-box oracle as FBS.
Some ancillary results needed for the proofs are listed in \Cref{sec:Algorithm_proof}.
Finally, \Cref{sec:Simulations} illustrates numerical results obtained with the proposed method.

\section{Preliminaries}
	\label{sec:Preliminaries}
	\subsection{Notation}
		\label{sec:Notation}
		The identity \(n\times n\) matrix is denoted as \(\id\), and the extended real line as \(\Rinf=\R\cup\set{\infty}\).
The open and closed ball of radius \(r\geq 0\) centered in \(x\in\R^n\) is denoted as \(\ball xr\) and \(\cball xr\), respectively.
Given a set \(E\) and a sequence \(\seq{x^k}\) we write \(\seq{x^k}\subset E\) with the obvious meaning of \(x^k\in E\) for all \(k\in\N\).
The (possibly empty) set of cluster points of \(\seq{x^k}\) is denoted as \(\omega\left(\seq{x^k}\right)\), or simply as \(\omega(x^k)\) whenever the indexing is clear from the context.
We say that \(\seq{x^k}\subset\R^n\) is \DEF{summable} if \(\sum_{k\in\N}\|x^k\|\) is finite, and \DEF{square-summable} if \(\seq{\|x^k\|^2}\) is summable.


Following the terminology of \cite{rockafellar2011variational}, we say that a function \(\funcf\) is \emph{strictly continuous at \(\bar x\)} if
\(
	\limsup_{
		\substack{y,z\to\bar x\\y\neq z}
	}{
		\frac{|f(y)-f(z)|}{\|y-z\|}
	}
\)
is finite, and \DEF{strictly differentiable at \(\bar x\)} if \(\nabla f(\bar x)\) exists and
\(
	\lim_{
		\substack{y,z\to\bar x\\y\neq z}
	}{
		\frac{f(y)-f(z) - \innprod{\nabla f(\bar x)}{y-z}}{\|y-z\|}
	}
{}={}
	0
\).
The set of functions \(\R^n\to\R\) with Lipschitz continuous gradient is denoted as \(\cont[1,1](\R^n)\), and for \(f\in\cont[1,1](\R^n)\) we write \(L_f\) to indicate the Lipschitz modulus of \(\nabla f\).


For a proper, closed function \(\funcg\), a vector \(v\in\partial g(x)\) is a \DEF{subgradient} of \(g\) at \(x\), where the \DEF{subdifferential} \(\partial g(x)\) is considered in the sense of \cite[Def. 8.3]{rockafellar2011variational}
\begin{align*}
	\partial g(x)
{}={} &
	\set{v\in\R^n}[
		\exists\seq{x^k}\to x,
		\seq{v^k\in\hat\partial g(x^k)}\to v
		~\text{s.t.}~
		g(x^k)\to g(x)
	],
\intertext{and \(\hat\partial g(x)\) is the set of \DEF{regular subgradients} of \(g\) at \(x\), namely}
	\hat\partial g(x)
{}={} &
	\set{v\in\R^n}[
		g(z)\geq g(x)
		{}+{}
		\innprod{v}{z-x}
		{}+{}
		o(\|z-x\|),
		~
		\forall z\in\R^n
		\vphantom{\seq{x^k}}
	].
\end{align*}
We have
\(
	\partial\varphi(x)
{}={}
	\nabla f(x)
	{}+{}
	\partial g(x)
\)
and
\(
	\hat\partial\varphi(x)
{}={}
	\nabla f(x)
	{}+{}
	\hat\partial g(x)
\)
\cite[Ex. 8.8(c)]{rockafellar2011variational}.


Given a parameter value \(\gamma>0\), the \DEF{Moreau envelope} function \(g^\gamma\) and the \DEF{proximal mapping} \(\prox_{\gamma g}\) are defined by
\begin{align}
\label{eq:Moreau}
	g^\gamma(x)
{}\coloneqq{} &
	\inf_z
	\set{
		g(z)
		+
		\tfrac{1}{2\gamma}
		\|z-x\|^2
	},
\\
\label{eq:Prox}
	\prox_{\gamma g}(x)
{}\coloneqq{} &
	\argmin_z
	\set{
		g(z)
		+
		\tfrac{1}{2\gamma}
		\|z-x\|^2
	}.
\end{align}
We now summarize some properties of \(g^\gamma\) and \(\prox_{\gamma g}\); the interested reader is referred to \cite{rockafellar2011variational} for a detailed discussion.
A function \(\funcg\) is \DEF{prox-bounded} if there exists \(\gamma>0\) such that \(g+\tfrac{1}{2\gamma}\|{}\cdot{}\|^2\) is bounded below on \(\R^n\).
The supremum of all such \(\gamma\) is the \DEF{threshold \(\gamma_g\) of prox-boundedness} for \(g\).
In particular, if \(g\) is convex or bounded below then \(\gamma_g=\infty\).
In general, for any \(\gamma\in(0,\gamma_g)\) the proximal mapping \(\prox_{\gamma g}\) is nonempty- and compact-valued, and the Moreau envelope \(g^\gamma\) finite \cite[Thm. 1.25]{rockafellar2011variational}.


Given a nonempty closed set \(S\subseteq\R^n\) we let \(\func{\indicator_S}{\R^n}{\Rinf}\) denote its \DEF{indicator function}, namely \(\indicator_S(x)=0\) if \(x\in S\) and \(\indicator_S(x)=\infty\) otherwise, and \(\ffunc{\proj_S}{\R^n}{\R^n}\) the (set-valued) \DEF{projection} \(x\mapsto\argmin_{z\in S}\|z-x\|\).
Proximal mappings can be seen as generalized projections, due to the relation \(\proj_S=\prox_{\gamma{\indicator_S}}\) for any \(\gamma>0\).


For a set-valued mapping \(\ffunc{T}{\R^n}{\R^n}\) we let \(\graph T=\set{(x,y)}[y\in T(x)]\) denote its \DEF{graph}, \(\zer T=\set{x\in\R^n}[0\in T(x)]\) the set of its \DEF{zeros} and \(\fix T=\set{x\in\R^n}[x\in T(x)]\) the set of its \DEF{fixed-points}.

	\subsection{Forward-backward iterations}
Due to the quadratic upper bound
\begin{equation}\label{eq:LipBound}
	f(z)
{}\leq{}
	f(x)+\tinnprod{\nabla f(x)}{z-x}+\tfrac{L_f}{2}\|z-x\|^2
\end{equation}
holding for all \(x,z\in\R^n\) \cite[Prop. A.24]{bertsekas1995nonlinear}, for any \(\gamma\in(0,\nicefrac{1}{L_f})\) the function
\begin{equation}\label{eq:lin}
	\lin zx[f,g]
{}\coloneqq{}
	f(x)+\innprod{\nabla f(x)}{z-x}+\tfrac{1}{2\gamma}\|z-x\|^2+g(z)
\end{equation}
furnishes a majorization model for \(\varphi\), in the sense that
\begin{itemize}
	\def\myVar{\fillwidthof[l]{x}{z}}
	\item
		\(\lin{\myVar}{x}[f,g]\geq\varphi(\myVar)\) for all \(x,z\in\R^n\), and
	\item
		\(\lin xx[f,g]=\varphi(x)\) for all \(x\in\R^n\).
\end{itemize}
Given a point \(x\in\R^n\), one iteration of \DEF{forward-backward splitting} (FBS) for problem \eqref{eq:Problem} consists in the minimization of the majorizing function \(\lin{}{}[f,g]\), namely, in selecting
\begin{equation}\label{eq:FB}
	x^+
{}\in{}
	\T[\gamma][f,g](x)
{}\coloneqq{}
	\fillwidthof[l]{\FB x,}{{\argmin}_z\lin zx[f,g],}
\end{equation}
where \(\gamma\in\bigl(0,\min\set{\gamma_g,\nicefrac{1}{L_f}}\bigr)\) is the stepsize parameter.
The (set-valued) \DEF{forward-backward operator} \(\T[\gamma][f,g]\) can be equivalently expressed as
\begin{subequations}\label{eq:TgammaRgamma}
\begin{align}
\label{eq:T}
	\hphantom{x^+\in{}}
	\T[\gamma][f,g](x)
{}={} &
	\FB x,
\shortintertext{%
	which motivates the bound \(\gamma<\gamma_g\) in \eqref{eq:FB} to ensure the existence of \(x^+\) for any \(x\).
	We also introduce the corresponding (set-valued) \emph{forward-backward residual}, namely%
}
\label{eq:R}
	\hphantom{x^+\in{}}
	\Res[\gamma][f,g](x)
{}\coloneqq{} &
	\fillwidthof[l]{\FB x,}{\tfrac 1\gamma\bigl(x-\T[\gamma][f,g](x)\bigr).}
\end{align}
\end{subequations}
Whenever no ambiguity occurs, we will omit the superscript and write simply \(\lin{}{}\), \(\T\) and \(\Res\) in place of \(\lin{}{}[f,g]\), \(\T[\gamma][f,g]\) and \(\Res[\gamma][f,g]\), respectively.

\eqref{eq:FB} emphasizes that FBS is a \DEF{majorization-minimization} algorithm (MM), a class of methods which has been thoroughly analyzed when the majorizing function is strongly convex in the first argument \cite{bolte2016majorization} (for \(\lin{}{}\), this is the case when \(g\) is convex).
MM algorithms are of interest whenever minimizing the surrogate function \(\lin{{}\cdot{}}{x}\) is significantly easier than directly addressing the non structured minimization of \(\varphi\).
For FBS this translates into simplicity of \(\prox_{\gamma g}\) and \(\nabla f\) operations, cf. \eqref{eq:T}.
Under very mild assumptions FBS iterations \eqref{eq:FB} converge to a \emph{critical point} (see \Cref{sec:StationaryCritical}) independently of the choice of \(x^+\) in the set \(\T(x)\) \cite{attouch2013convergence}.
The key is the following well known \emph{sufficient decrease} property, whose proof can be found in \cite[Lem. 2]{bolte2014proximal}.
\begin{lem}[Sufficient decrease]\label{prop:FBSnonadaptive}%
For any \(\gamma\in(0,\gamma_g)\), \(x\in\R^n\) and \(\bar x\in\T(x)\) it holds that
\(
	\varphi(\bar x)
{}\leq{}
	\varphi(x)
	{}-{}
	\tfrac{1-\gamma L_f}{2\gamma}
	\|x-\bar x\|^2
\).
\end{lem}

\section{Stationary and critical points}
	\label{sec:StationaryCritical}
Unless \(\varphi\) is convex, the \emph{stationarity} condition \(0\in\hat\partial\varphi(x^\star)\) in problem \eqref{eq:Problem} is only necessary for the optimality of \(x^\star\) \cite[Thm. 10.1]{rockafellar2011variational}.
In this section we define different concepts of (sub)optimality and show how they are related for generic functions \(\varphi=f+g\) as in \Cref{ass:fg}. 
\begin{defin}
We say that a point \(x^\star\in\dom\varphi\) is
\begin{enumerate}
	\item
		\DEF{stationary} if
		\(
			0\in\hat\partial\varphi(x^\star)
		\);
	\item\label{def:critical}%
		\DEF{critical} if it is \DEF{\(\gamma\)-critical} for some \(\gamma\in(0,\gamma_g)\), \ie if
		\(
			x^\star\in \T(x^\star)
		\);
	\item
		\DEF{optimal} if
		\(
			x^\star\in\argmin\varphi
		\),
		\ie if it solves \eqref{eq:Problem}.
\end{enumerate}
\end{defin}
\arxiv{%
	In the definition of stationarity we use the \emph{regular} subdifferential \(\hat\partial\varphi\) in order to impose a stronger condition than the `classical' one \(0\in\partial\varphi(x^\star)\), since \(\hat\partial\varphi\subseteq\partial\varphi\).
}%
The notion of criticality was already discussed in \cite{beck2016minimization} under the name of \(L\)-stationarity (\(L\) plays the role of \(\nicefrac 1\gamma\)) for the special case of \(g=\indicator_{B\cap C_s}\), where \(B\) is a convex set and \(C_s\) is the (nonconvex) set of vectors with at most \(s\) nonzero entries.

If \(g\) is convex, then \(\gamma_g=\infty\) and we may talk of criticality without mention of \(\gamma\): in this case, the properties of \(\gamma\)-criticality and stationarity are equivalent regardless of the value of \(\gamma\).
For more general functions \(g\), instead, the value of \(\gamma\) plays a role in determining whether a point is \(\gamma\)-critical or not, which legitimizes the following definition.
\begin{defin}
The \DEF{criticality threshold} is the function \(\func{\Gamma^{f,g}}{\R^n}{[0,\gamma_g]}\)
\begin{equation}\label{eq:Gamma}
	\Gamma^{f,g}(x)
{}\coloneqq{}
	\sup\bigl(\set{\gamma>0}[{x\in\T[\gamma][f,g](x)}]\cup\set{0}\bigr)
\quad
	\text{for \(x\in\R^n\).}
\end{equation}
\end{defin}
As usual, whenever \(f\) and \(g\) are clear from the context we simply write \(\Gamma\) in place of \(\Gamma^{f,g}\).
That \(\Gamma\leq\gamma_g\) is due to the fact that \(\prox_{\gamma g}\) (and consequently \(\T\)) is everywhere empty-valued for \(\gamma>\gamma_g\).
Considering also \(\gamma=0\) forces the set in the definition to be nonempty, and the lower-bound \(\Gamma\geq 0\) in particular; more precisely, observe that, by definition, \(\Gamma(x)>0\) iff \(x\) is a critical point.

\begin{es}\label{es:Gamma}%
Let us consider \(\varphi = f+g\) for \(f(x)=\frac 12x^2\) and \(g = \indicator_C\) where \(C=\set{\pm 1}\).
Clearly, \(\gamma_g = +\infty\) (as \(g\) is lower-bounded), \(L_f = 1\) and \(\pm 1\) are both (unique) optima.
Since \(\hat\partial\varphi(x)=\R\) for \(x\in C\) and \(\hat\partial\varphi\) is clearly empty elsewhere, all points in \(C\) are stationary.
\(\prox_{\gamma g}\) is the (set-valued) projection on \(C\), therefore the forward-backward operator is
\(
	\T(x)
{}={}
	\proj_C((1-\gamma)x)
\).
We have
\[
	\T(-1)
{}={}
	\begin{cases}[@{\!}c@{~~\text{if }} l@{}]
		\set{-1}    & \gamma<1 \\
		\set{\pm 1} & \gamma=1 \\
		\set{1}     & \gamma>1
	\end{cases}
\quad\text{and}\quad
	\T(1)
{}={}
	\begin{cases}[@{\!}c@{~~\text{if }} l@{}]
		\set{1}     & \gamma<1 \\
		\set{\pm 1} & \gamma=1 \\
		\set{-1}    & \gamma>1.
	\end{cases}
\]
In particular, \(\Gamma(1)=\Gamma(-1)=1\).
\end{es}
We now list some properties of critical and optimal points which will be used to derive regularity properties of \(\T\) and \(g^\gamma\).
\begin{thm}[Properties of critical points]\label{prop:GammaCriticality}%
The following properties hold:
\begin{enumerate}
	\item\label{prop:ProxGrad}%
		for \(\gamma\in(0,\gamma_g)\), a point \(x^\star\) is \(\gamma\)-critical iff
		\[
				g(x)
		{}\geq{}
			g(x^\star)
			{}+{}
			\innprod{-\nabla f(x^\star)}{x-x^\star}
			{}-{}
			\tfrac{1}{2\gamma}
			\|x-x^\star\|^2
		\qquad
			\forall x\in\R^n;
		\]
\arxiv{%
	\item\label{prop:CriticalGradf}%
		if \(x^\star\) and \(y^\star\) are \(\gamma\)-critical, then
		\(\mathtight
			\innprod{
				\nabla f(y^\star)-\nabla f(x^\star)
			}{
				y^\star-x^\star
			}
		{}\leq{}
			\tfrac1\gamma
			\|y^\star-x^\star\|^2
		\);
		in particular, if \(f\) is \(\mu\)-strongly convex, then for all \(\gamma\geq\nicefrac1\mu\) there exists at most one \(\gamma\)-critical point;
	\item\label{prop:GammaSmaller}%
		if \(x^\star\) is \(\gamma\)-critical, then it is also \(\gamma'\)-critical for any \(\gamma'\in(0,\gamma]\);
}%
	\item\label{prop:Gamma}%
		if \(x^\star\) is critical, then it is \(\gamma\)-critical for all \(\gamma\in(0,\Gamma(x^\star))\); moreover, \(x^\star\) is also \(\Gamma(x^\star)\)-crit\-i\-cal provided that \(\Gamma(x^\star)<\gamma_g\);
	\item\label{prop:SingleValuedFB}%
		\(\T(x^\star){=}\set{x^\star}\) and \(\Res(x^\star){=}\set{0}\) for any critical point \(x^\star\) and \(\gamma\in(0,\Gamma(x^\star))\).
\end{enumerate}
\begin{proof}
\begin{proofitemize}
\item\ref{prop:ProxGrad}:
	by definition, \(x^\star\) is \(\gamma\)-critical iff \(\lin{x^\star}{x^\star}\leq\lin{x}{x^\star}\) for all \(x\), \ie iff
	\[
		f(x^\star)
		{}+{}
		g(x^\star)
	{}\leq{}
		f(x^\star)
		{}+{}
		\innprod{\nabla f(x^\star)}{x-x^\star}
		{}+{}
		\tfrac{1}{2\gamma}
		\|x-x^\star\|^2
		{}+{}
		g(x)
	\qquad
		\forall x\in\R^n.
	\]
	By suitably rearranging, the claim readily follows.
\arxiv{%
	\item\ref{prop:CriticalGradf}:
		suppose now that \(x^\star\) and \(y^\star\) are \(\gamma\)-critical.
		Plugging \(x=y^\star\) in \ref{prop:ProxGrad}, then interchanging \(x^\star\) and \(y^\star\) and summing the obtained inequalities yields the one in \ref{prop:CriticalGradf}.
		If, additionally, \(f\) is \(\mu\)-strongly convex then
		\[
			\mu
			\|y^\star-x^\star\|^2
		{}\leq{}
			\innprod{
					\nabla f(y^\star)-\nabla f(x^\star)
				}{
					y^\star-x^\star
				}
		{}\leq{}
			\tfrac1\gamma
			\|y^\star-x^\star\|^2
		\]
		where the first inequality is due to \cite[Thm. 2.1.9]{nesterov2004introductory}.
		Therefore, if \(\gamma\geq\nicefrac 1\mu\) then necessarily \(x^\star=y^\star\).
\item\ref{prop:GammaSmaller}:
	follows from the characterization \ref{prop:ProxGrad} of \(\gamma\)-criticality.
}%
\item\ref{prop:Gamma}:
	\siam{%
		due to \ref{prop:ProxGrad}, if \(x^\star\) is \(\gamma\)-critical, apparently it is also \(\gamma'\)-critical for any \(\gamma'\in(0,\gamma]\).
	}%
	From
	\arxiv{\ref{prop:GammaSmaller} and}
	the definition \eqref{eq:Gamma} of the criticality threshold \(\Gamma(x^\star)\), it then follows that \(x^\star\) is \(\gamma\)-critical for any \(\gamma\in(0,\Gamma(x^\star))\).
	Suppose now that \(\Gamma(x^\star)<\gamma_g\).
	Then, due to \ref{prop:ProxGrad} for all \(\gamma\in(0,\Gamma(x^\star))\) we have
	\[
		g(x)
	{}\geq{}
		g(x^\star)
		{}+{}
		\innprod{-\nabla f(x^\star)}{x-x^\star}
		{}-{}
		\tfrac{1}{2\gamma}
		\|x-x^\star\|^2
	\quad
		\forall x\in\R^n.
	\]
	By taking the limit as \(\gamma\nearrow\Gamma(x^\star)\) we obtain that the inequality holds for \(\Gamma(x^\star)\) as well, proving the claim in light of the characterization \ref{prop:ProxGrad}.
\item\ref{prop:SingleValuedFB}:
	let \(x^\star\) be a critical point, and let \(x\in \T(x^\star)\) for some \(\gamma<\Gamma(x^\star)\).
	Fix \(\gamma'\in(\gamma,\Gamma(x^\star))\).
	From \ref{prop:ProxGrad} and \ref{prop:Gamma} it then follows that
	\begin{equation}\label{eq:StrictInequality}
		g(x)
	{}\geq{}
		g(x^\star)
		{}+{}
		\innprod{-\nabla f(x^\star)}{x-x^\star}
		{}-{}
		\tfrac{1}{2\gamma'}
		\|x-x^\star\|^2.
	\end{equation}
	Since \(x,x^\star\in \T(x^\star)\), it holds that \(\lin{x^\star}{x^\star}=\lin{x}{x^\star}\), \ie
	\[
		g(x^\star)
	{}={}
		\innprod{\nabla f(x^\star)}{x-x^\star}
		{}+{}
		\tfrac{1}{2\gamma}
		\|x-x^\star\|^2
		{}+{}
		g(x)
	{}\overrel[\geq]{\eqref{eq:StrictInequality}}{}
		g(x^\star)
		{}+{}
		\bigl(\tfrac{1}{2\gamma}-\tfrac{1}{2\gamma'}\bigr)
		\|x-x^\star\|^2.
	\]
	Since \(\tfrac{1}{2\gamma}-\tfrac{1}{2\gamma'}>0\), necessarily \(x=x^\star\).
	\newarxiv{\qedhere}%
\end{proofitemize}
\end{proof}
\end{thm}
The inequality in \Cref{prop:ProxGrad} can be rephrased as the fact that the vector \(-\nabla f(\bar x)\) is a ``global'' \DEF{proximal subgradient} for \(g\) at \(\bar x\) as in \cite[Def. 8.45]{rockafellar2011variational}, where ``global'' refers to the fact that \(\delta\) can be taken \(+\infty\) in the cited definition.
An interesting consequence is that the definition of criticality depends solely on \(\varphi\) and \emph{not} on the considered decomposition \(f+g\); in fact, it is only the threshold \(\Gamma\) that depends on it.
To see this, let \(\tilde f=f-h\) and \(\tilde g=g+h\) for some \(h\in\cont[1,1](\R^n)\), and consider a point \(x^\star\) which is \(\gamma\)-critical with respect to the decomposition \(f+g\), \ie such that \(x^\star\in\T[\gamma][f,g](x^\star)\).
Combining \Cref{prop:ProxGrad} with the quadratic bound \eqref{eq:LipBound} for \(h\), we obtain
\[
	\let\OLDtilde\tilde
	\renewcommand\tilde[1]{\smash{\OLDtilde{#1}}}
	\tilde g(x)
{}\geq{}
	\tilde g(x^\star)
	{}-{}
	\tinnprod{\nabla \tilde f(x^\star)}{x-x^\star}
	{}-{}
	\smash{\tfrac{1}{2\frac{\gamma}{1+\gamma L_h}}}
	\|x-x^\star\|^2
\quad
	\text{for all \(x\in\R^n\).}
\]
Again from the characterization of \Cref{prop:ProxGrad}, we deduce that
\(
	\smash{
		x^\star\in\T[\tilde\gamma][\tilde f,\tilde g](x^\star)
	}
\),
where
\(
	\tilde\gamma
{}={}
	\frac{\gamma}{1+\gamma L_h}
\).
In particular, considering \(h=-f\) we infer that a point \(x^\star\) is critical iff \(x^\star\in\T[\gamma][0,\varphi](x^\star)=\prox_{\gamma\varphi}(x^\star)\) for some \(\gamma>0\), which legitimizes the notion of criticality without mentioning a specific decomposition.

In the next result we show that criticality is a halfway property between stationarity and optimality.
In light of these relations we shall seek ``suboptimal'' solutions which we characterize as critical points.
\begin{prop}[Optimality, criticality, stationarity]\label{prop:Critical}%
Let
\(
	\bar\gamma
{}\coloneqq{}
	\min\set{\gamma_g,\nicefrac{1}{L_f}}
\).
\begin{enumerate}
	\item\label{prop:FixTinZer}%
		(criticality \(\Rightarrow\) stationarity)~
		\(\fix \T\subseteq\zer\hat\partial\varphi\) for all \(\gamma\in(0,\gamma_g)\);
	\item\label{prop:argmininFixT}%
		(optimality \(\Rightarrow\) criticality)~
		\(\Gamma(x^\star)\geq\bar\gamma\) for all \(x^\star\in\argmin\varphi\); in particular, \(\argmin\varphi\subseteq\fix \T\) for all \(\gamma\in(0,\bar\gamma)\), and also for \(\gamma=\nicefrac{1}{L_f}\) if \(\gamma_g>\nicefrac{1}{L_f}\);
\arxiv{%
	\item\label{prop:argminUniqueT}%
		\(\T(x^\star)=\set{x^\star}\) and \(\Res(x^\star)=\{0\}\) for all \(x^\star\in\argmin\varphi\) and \(\gamma\in(0,\bar\gamma)\).
}%
\end{enumerate}
\begin{proof}
\begin{proofitemize}
\item\ref{prop:FixTinZer}:
	let \(\gamma\in(0,\gamma_g)\) and \(x\in\fix \T\).
	Since \(x\) minimizes
	\(
		g+\tfrac{1}{2\gamma}\|{}\cdot{}-x+\gamma\nabla f(x)\|^2
	\),
	we have
	\(
		0
	{}\in{}
		\hat\partial\bigl[
			g+\tfrac{1}{2\gamma}\|{}\cdot{}-x+\gamma\nabla f(x)\|^2
		\bigr](x)
	{}={}
		\hat\partial g(x)+\nabla f(x)
	{}={}
		\hat\partial\varphi(x)
	\),
	where the first inclusion follows from \cite[Thm. 10.1]{rockafellar2011variational} and the equalities from \cite[Thm. 8.8(c)]{rockafellar2011variational}.
	This proves that \(x\) is stationary.
\item\ref{prop:argmininFixT}:
	Fix \(\gamma\in(0,\bar\gamma)\), \(x^\star\in\argmin\varphi\) and \(y\in \T(x^\star)\).
	Necessarily \(y=x^\star\), otherwise, due to \cref{prop:FBSnonadaptive}, \(\varphi(y)\) would contradict minimality of \(\varphi(x^\star)\).
	Therefore, \(x^\star\) is \(\gamma\)-critical and the claim follows from the arbitrarity of \(\gamma\in(0,\bar\gamma)\).
\arxiv{%
	\item\ref{prop:argminUniqueT}:
		straightforward consequence of \ref{prop:argmininFixT} and \cref{prop:SingleValuedFB}.
}%
\newarxiv{\qedhere}%
\end{proofitemize}
\end{proof}
\end{prop}

As already seen in \Cref{es:Gamma}, the bound \(\Gamma(x^\star)\geq\min\set{\gamma_g,\nicefrac{1}{L_f}}\) at optimal points in \Cref{prop:argmininFixT} is tight, and clearly the implication \emph{``optimality \(\Rightarrow\) criticality''} cannot be reversed (consider, \eg the point \(x^\star=0\) for \(\varphi=\cos\)).
The next example shows that the other implication is also proper.
\begin{es}[Stationarity \(\not\Rightarrow\) criticality]
Let \(f(x)=\frac12x^2\) and \(g(x)=x^{\nicefrac 53}\).
We have \(\gamma_g=+\infty\), \(L_f=1\), and for \(x^\star=0\) it holds that \(\hat\partial\varphi(x^\star)=\set{\nabla\varphi(x^\star)}=\set{0}\).
Therefore, \(x^\star\) is stationary; however, \(\T(x^\star)=\prox_{\gamma g}(0)=\set{-(\nicefrac{5\gamma}{3})^3}\), and in particular \(x^\star\notin \T(x^\star)\) for any \(\gamma>0\), proving \(x^\star\) to be non critical.
\end{es}
\arxiv{%
	We conclude the section showing that criticality entails favorable regularity properties on the Moreau envelope.
	\begin{prop}\label{prop:criticalDifferentiable}%
	If \(\bar x\) is critical, then for all \(\gamma\in(0,\Gamma(\bar x))\) the Moreau envelope \(g^\gamma\) is strictly differentiable at \(\Fw{\bar x}\) with \(\nabla g^\gamma(\Fw{\bar x})=-\nabla f(\bar x)\).
	\begin{proof}
	It follows from \cite[Ex. 10.32]{rockafellar2011variational} that \(g^\gamma\) is strictly continuous with
	\[
		\partial g^\gamma
		\bigl(
			\Fw{\bar x}
		\bigr)
	{}\subseteq{}
		\tfrac 1\gamma
		\bigl[
			\Fw{\bar x}
			{}-{}
			\T(\bar x)
		\bigr]
	{}={}
		\tfrac 1\gamma
		\bigl[
			\Res(\bar x)
			{}-{}
			\gamma\nabla f(\bar x)
		\bigr].
	\]
	From \cref{prop:SingleValuedFB} we then have that
	\(
		\partial g^\gamma
		\bigl(
			\Fw{\bar x}
		\bigr)
	{}\subseteq{}
		\set{-\nabla f(\bar x)}
	\)
	provided that \(\gamma\in(0,\Gamma(\bar x))\), and the proof follows from \cite[Thm. 9.18(b)]{rockafellar2011variational}.
	\end{proof}
	\end{prop}
}%

\section{Forward-backward envelope}
	\label{sec:FBE}
	The FBE \eqref{eq:FBE1} was introduced in \cite{patrinos2013proximal} and further analyzed in \cite{stella2017forward,liu2017further} in the case when \(g\) is convex.
Under such assumption the FBE was shown to be continuously differentiable, which made it possible to derive minimization algorithms based on its gradient.
In the general setting addressed in this paper the FBE might fail to be (continuously) differentiable, and as such we need to resort to \emph{gradient-free} methods.
This task will be addressed in \Cref{sec:Algorithm} where \Cref{alg:zerofpr} will be proposed; other than being applicable to a wider range of problems, the proposed scheme is entirely based on the same oracle of forward-backward iterations, unlike the approaches in \cite{patrinos2013proximal,stella2017forward,liu2017further} which instead require the computation of \(\nabla^2f\).
All this will be possible thanks to continuity properties of the FBE, and to the behavior of \(\varphi_\gamma\) at critical points.
We now focus on its continuity, while the other property will be addressed shortly after in \Cref{prop:EquivFBE}.

\begin{rem}[Alternative expressions for \(\varphi_\gamma\)]
By expanding the square and rearranging the terms in the definition \eqref{eq:FBE1}, \(\varphi_\gamma\) can equivalently be expressed as
\[
	\varphi_\gamma(x)
{}={}
	\inf_{z\in\R^n}
	\set{
		f(x)
		{}-{}
		\tfrac\gamma2
		\|\nabla f(x)\|^2
		{}+{}
		g(z)
		{}+{}
		\tfrac{1}{2\gamma}
		\|z-x+\gamma\nabla f(x)\|^2
	}.
\]
Comparing with \eqref{eq:FB}, it is apparent that the set of minimizers \(z\) in the above expression coincides with \(\T(x)\), the forward-backward operator at \(x\).
Moreover, taking out the constant term
\(
	f(x)
	{}-{}
	\tfrac\gamma2
	\|\nabla f(x)\|^2
\)
from the infimum we immediately obtain the following expression involving the Moreau envelope of \(g\):
\begin{equation}\label{eq:FBEMoreau}
	\newarxiv{\pushQED{\qed}}
	\varphi_\gamma(x)
{}={}
	f(x)
	{}-{}
	\tfrac\gamma2
	\|\nabla f(x)\|^2
	{}+{}
	g^\gamma(x-\gamma\nabla f(x)).
	\newarxiv{\qedhere\popQED}
\end{equation}
\newarxiv{\let\qed\relax}%
\end{rem}
Other than providing an explicit way of computing the FBE, \eqref{eq:FBEMoreau} emphasizes how \(\varphi_\gamma\) inherits the regularity properties of the Moreau envelope of \(g\).
In particular, the next key property follows from the strict continuity of \(g^\gamma\) \cite[Ex. 10.32]{rockafellar2011variational}.
\begin{prop}[Strict continuity of \(\varphi_\gamma\)]\label{prop:FBEC+}%
For any \(\gamma\in(0,\gamma_g)\), the FBE \(\varphi_\gamma\) is a real-valued and strictly continuous function on \(\R^n\).
\end{prop}
	\subsection{Connections with the Moreau envelope}
		\label{subsec:Moreau}
		For the special case \(f=0\), FBS iterations \eqref{eq:FB} reduce to the \emph{proximal point algorithm} (PPA) \(x^+\in\prox_{\gamma\varphi}(x)\), first introduced in \cite{martinet1970breve} for convex functions \(\varphi\) and later generalized for functions with convex majorizing surrogate
\(
	\lin{{}\cdot{}}{x}[\,0,\varphi]
{}={}
	\varphi({}\cdot{})+\frac{1}{2\gamma}\|{}\cdot{}-x\|^2
\),
see \eg \cite{kaplan1998proximal}.
Similarly, the FBE reduces to the Moreau envelope \(\varphi^\gamma=\varphi_\gamma^{0,\varphi}\).
In fact, the FBE extends the connection between PPA and Moreau envelope
\begin{subequations}
\begin{alignat}{5}
	\varphi^\gamma(x)
{}={} &
	{\min}_z\lin zx[\,0,\varphi]
&\quad\leftrightarrow\quad&
	\prox_{\gamma\varphi}(x)
{}={}
	{\argmin}_z\lin zx[\,0,\varphi],
\shortintertext{holding for \(f=0\) in \eqref{eq:lin}, to majorizing functions \(\lin{}{}[f,g]\) with arbitrary \(f\in C^{1,1}(\R^n)\)}
\label{eq:FBET}
	\varphi_\gamma(x)
{}={} &
	{\min}_z\lin zx[f,g]
&\quad\leftrightarrow\quad&
	\hspace*{0pt}
	\fillwidthof[r]{
		{}\prox_{\gamma\varphi}(x)
	}{
		\T(x)
	}
{}={}
	{\argmin}_z\lin zx[f,g].
\end{alignat}
\end{subequations}
In the next section we will see the fundamental qualitative similarities between the FBE and the Moreau envelope.
Namely, for \(\gamma\) small enough both \(\varphi^\gamma\) and \(\varphi_\gamma\) are lower bounds for the original function \(\varphi\) with same minimizers and minimum; in particular the minimization of \(\varphi\) is equivalent to that of \(\varphi^\gamma\) or \(\varphi_\gamma\).
Similarly, the identity
\begin{align*}
	\varphi(\bar x)
{}={} &
	\varphi^\gamma(x)
	{}-{}
	\fillwidthof[l]{
		\tfrac{1-\gamma L_f}{2\gamma}\|x-\bar x\|^2
	}{
		\tfrac{1}{2\gamma}\|x-\bar x\|^2
	}
\quad
	\text{for \(\bar x\in\prox_{\gamma\varphi}(x)\)}
\shortintertext{%
	will be extended to the inequality%
}
	\varphi(\bar x)
{}\leq{} &
	\varphi_\gamma(x)
	{}-{}
	\tfrac{1-\gamma L_f}{2\gamma}\|x-\bar x\|^2
\quad
	\text{for \(\bar x\in \T(x)\).}
\end{align*}

	\subsection{Basic properties}
		\label{subsec:FBEbasic}
		We now provide bounds relating \(\varphi_\gamma\) to the original function \(\varphi\) that extend the well known inequalities involving the Moreau envelope.
\begin{prop}\label{prop:Ineq}%
Let \(\gamma\in(0,\gamma_g)\) be fixed.
Then
\begin{enumerate}
	\item\label{prop:FBEleq}%
		\(
			\varphi_\gamma
		{}\leq{}
			\varphi
		\).
	\item\label{prop:FBEgeqGlobal}%
		\(
			\varphi(\bar x)
		{}\leq{}
			\varphi_\gamma(x)
			{}-{}
			\tfrac{1-\gamma L_f}{2\gamma}
			\|x-\bar x\|^2
		\)
		for all \(x\in\R^n\) and \(\bar x\in \T(x)\).
\end{enumerate}
\begin{proof}
\ref{prop:FBEleq} is obvious from the definition of the FBE (consider \(z=x\) in \eqref{eq:FBE1}).
As to \ref{prop:FBEgeqGlobal}, since the set of minimizers in \eqref{eq:FBE1} is \(\T(x)\) (cf. \eqref{eq:FBET}), \eqref{eq:LipBound} yields
\begin{align*}
	\mathtight
	\varphi_\gamma(x)
{}={} &
		f(x)
		{}+{}
		\innprod{\nabla f(x)}{\bar x-x}
	{}+{}
	g(\bar x)
	{}+{}
	\tfrac{1}{2\gamma}
	\|x-\bar x\|^2
\\
{}\geq{} &
	f(\bar x)
	{}-{}
	\nicefrac{L_f}{2}
	\|\bar x-x\|^2
	{}+{}
	g(\bar x)
	{}+{}
	\tfrac{1}{2\gamma}
	\|x-\bar x\|^2
{}={}
	\varphi(\bar x)
	{}+{}
	\tfrac{1-\gamma L_f}{2\gamma}
	\|x-\bar x\|^2.
	\newarxiv{\qedhere}
\end{align*}
\end{proof}
\end{prop}
With respect to the inequalities holding for convex \(g\) treated in \cite{stella2017forward}, the lower bound in \Cref{prop:Ineq} is weaker, while the upper bound unchanged.
Regardless, an immediate consequence of the result is that the value of \(\varphi\) and \(\varphi_\gamma\) at critical points is the same, and minimizers and infima of the two functions coincide for \(\gamma\) small enough.
\begin{thm}\label{prop:EquivFBE}%
The following hold
\begin{enumerate}
	\item\label{prop:FBEfix}%
		\(\varphi(x)=\varphi_\gamma(x)\)~ for all \(\gamma\in(0,\gamma_g)\) and \(x\in\fix \T\);
	\item\label{prop:FBEinf}%
		\(
			\inf\varphi=\inf\varphi_\gamma
		\)
		~and~
		\(\argmin\varphi=\argmin\varphi_\gamma\)
		~
		for all
		\(
			\gamma
		{}\in{}
			\bigl(
				0,
				\min\set{\nicefrac{1}{L_f},\gamma_g}
			\bigr)
		\).
\end{enumerate}
\end{thm}
The bound \(\gamma<\nicefrac{1}{L_f}\) in \Cref{prop:FBEinf} is tight even when \(f\) and \(g\) are convex, as the counterexample with \(f(x)=\frac12 x^2\) and \(g=\indicator_{\R_+}\) shows (see \cite[Ex. 2.4]{stella2017forward} for details).

Although we will address problem \eqref{eq:Problem} by simply exploiting the \emph{continuity} of the FBE, nevertheless \(\varphi_\gamma\) enjoys favorable properties which are key for the efficacy of the method which will be discussed in \Cref{sec:Algorithm}.
Firstly, observe that, due to \emph{strict} continuity, \(\varphi_\gamma\) is almost everywhere differentiable, as it follows from Rademacher's theorem.
The same applies to the mapping \(x\mapsto x-\gamma\nabla f(x)\), its Jacobian being
\begin{equation}\label{eq:Q}
	Q_\gamma(x)
{}\coloneqq{}
	\id - \gamma\nabla^2\!f(x)
\end{equation}
which is symmetric wherever it exists \cite[Cor. 13.42 and Prop. 13.34]{rockafellar2011variational}.
However, in order to show that the proposed method achieves fast convergence we need additional regularity properties, namely (strict) twice differentiability at critical points and continuous differentiability around.
The rest of the section is dedicated to this task.

	\subsection{Prox-regularity and first-order properties}
		\label{subsec:FBEfirst}
		In the favorable case in which \(g\) is convex and \(f\in C^2(\R^n)\), the FBE enjoys global continuous differentiability \cite{stella2017forward}.
In our setting, \DEF{prox-regularity} acts as a surrogate of convexity; the interested reader is referred to \cite[§13.F]{rockafellar2011variational} for a detailed discussion.
\begin{defin}[Prox-regularity]\label{def:ProxRegularity}%
Function \(g\) is said to be \DEF{prox-regular} at \(x_0\) for \(v_0\in\partial g(x_0)\) if there exist \(\rho,\varepsilon>0\) such that for all \(x'\in \ball{x_0}{\varepsilon}\) and
\[
	(x,v)
{}\in{}
	\graph\partial g
\quad\text{s.t.}\quad
	x\in\ball{x_0}{\varepsilon},~
	v\in\ball{v_0}{\varepsilon},~
	\text{and}~
	g(x)\leq g(x_0)+\varepsilon
\]
it holds that
\(
	g(x')
{}\geq{}
	g(x)
	{}+{}
	\innprod{v}{x'-x}
	{}-{}
	\tfrac \rho2
	\|x'-x\|^2
\).
\end{defin}
\arxiv{%
	To help better visualize this definition we consider the local geometrical property that it entails on the function's epigraph \cite[Cor. 3.4, Thm. 3.5]{poliquin1996prox-regular}.
	If \(g\) is prox-regular at \(x_0\) for \(v_0\) for some constants \(\varepsilon,\rho>0\) as in \Cref{def:ProxRegularity}, then there exists a neighborhood of \((x_0+\nicefrac{v_0}{\rho},g(x_0)-\nicefrac1\rho)\) in which the projection on
	\(
		\epi g
	{}\cap{}
		\bigl(
			\ball{x_0}{\varepsilon}\times\ball{v_0}{\varepsilon}
		\bigr)
	\)
	is single-valued.
}%
Prox-regularity is a mild requirement enjoyed globally and for any subgradient by all convex functions, with \(\varepsilon=+\infty\) and \(\rho=0\).
When \(g\) is prox-regular at \(x_0\) for \(v_0\), then for sufficiently small \(\gamma>0\) the Moreau envelope \(g^\gamma\) is continuously differentiable in a neighborhood of \(x_0+\gamma v_0\) \cite{poliquin1996prox-regular}.
To our purposes, when needed, prox-regularity of \(g\) will be required only at critical points \(x^\star\), and only for the subgradient \(-\nabla f(x^\star)\).
Therefore, with a slight abuse of terminology we define prox-regularity \emph{of} critical points as follows.
\begin{defin}[Prox-regularity \emph{of} critical points]
We say that a critical point \(x^\star\) is \DEF{prox-regular} if \(g\) is prox-regular at \(x^\star\) for \(-\nabla f(x^\star)\).
\end{defin}
\arxiv{%
	Clearly, if \(g\) is convex then any critical point is prox-regular.
	Prox-regularity \emph{of} critical points is really a mild requirement, also considering that the fact of being critical itself entails some regularity properties as shown in \Cref{prop:criticalDifferentiable}.
}%
Examples where
\switch{
	a critical point fails to be prox-regular
}{
	this fails to hold
}%
are of challenging construction; before illustrating a cumbersome such instance in \Cref{es:NonProxRegular}, we first prove an important result that connects prox-regularity with first-order properties of the FBE.
\begin{thm}[Continuous differentiability of \(\varphi_\gamma\)]\label{thm:1stOrder}%
Suppose that \(f\) is of class \(\cont[2]\) around a prox-regular critical point \(x^\star\).
Then, for all \(\gamma\in(0,\Gamma(x^\star))\) there exists a neighborhood \(U_{x^\star}\) of \(x^\star\) on which the following properties hold:
\begin{enumerate}
	\item\label{thm:FBLip}%
		\(\T\) and \(\Res\) are strictly continuous, and in particular single-valued;
	\item\label{thm:FBEC1}%
		\(\varphi_\gamma\in \cont[1]\) with
		\(
			\nabla\varphi_\gamma
		{}={}
			Q_\gamma \Res
		\),
		where \(Q_\gamma\) is as in \eqref{eq:Q}.
\end{enumerate}
\arxiv{%
	Moreover, if \(\gamma\in(0,\min\set{\Gamma(x^\star),\nicefrac{1}{L_f}})\) then
	\begin{enumerate}[resume]
		\item\label{thm:FBELocMin}\hl{REMOVE: will be done later with additional assumption of 2epidiff of g}
			\(x^\star\) is a (strong) local minimum for \(\varphi\) iff it is a (strong) local minimum for \(\varphi_\gamma\).
	\end{enumerate}
}%
\begin{proof}
For \(\gamma'\in(\gamma,\Gamma(x^\star))\), using \cref{{prop:ProxGrad},,{prop:SingleValuedFB}} we obtain that
\begin{equation}\label{eq:ProxGrad'}
	g(x)
{}\geq{}
	g(x^\star)
	{}-{}
	\innprod{\nabla f(x^\star)}{x-x^\star}
	{}-{}
	\tfrac{1}{2\gamma'}
	\|x-x^\star\|^2
\qquad
	\forall x\in\R^n.
\end{equation}
Replacing \(\gamma'\) with \(\gamma\) in the above expression, the inequality is strict for all \(x\neq x^\star\).
From \cite[Thm. 4.4]{poliquin1996prox-regular} applied to the ``tilted'' function
\(
	x
{}\mapsto{}
	g(x+x^\star)
	{}-{}
	g(x^\star)
	{}-{}
	\innprod{\nabla f(x^\star)}{x}
\)
it follows that there is a neighborhood \(V\) of \(\Fw{x^\star}\) in which \(\prox_{\gamma g}\) is strictly continuous and \(g^\gamma\) is of class \(\cont[1+]\) with
\newarxiv{%
	gradient
}%
\switch{
	\(
		\nabla g^\gamma(x)
	{}={}
		\gamma^{-1}
		\left(
			x-\prox_{\gamma g}(x)
		\right)
	\)
	for all \(x\in V\).
}{
	\[
		\nabla g^\gamma(x)
	{}={}
		\gamma^{-1}
		\left(
			x-\prox_{\gamma g}(x)
		\right)
		\qquad\forall x\in V.
	\]
}%
By possibly narrowing \(U_{x^\star}\), we may assume that \(f\in\cont[2](U_{x^\star})\) and \(\Fw x\in V\) for all \(x\in U_{x^\star}\).
\ref{thm:FBEC1} then follows from \eqref{eq:FBEMoreau} and the chain rule of differentiation, and \ref{thm:FBLip} from the fact that strict continuity is preserved by composition.
\arxiv{%

We now show \ref{thm:FBELocMin}.
From \cref{{prop:FBEleq},,{prop:FBEfix}} it follows that
\[
	\varphi(x)
	{}-{}
	\varphi(x^\star)
{}\geq{}
	\varphi_\gamma(x)
	{}-{}
	\varphi_\gamma(x^\star)
\qquad
	\forall x\in\R^n,
\]
and therefore any (strong) local minimum of \(\varphi_\gamma\) is also a (strong) local minimum of \(\varphi\).
Conversely, suppose that \(\gamma<\nicefrac{1}{L_f}\) and that there exist \(\varepsilon\geq 0\) and a neighborhood \(U_{x^\star}\) of \(x^\star\) such that
\[
	\varphi(x)-\varphi(x^\star)
{}\geq{}
	\tfrac \varepsilon 2
	\|x-x^\star\|
\quad
	\text{for all }
	x\in U_{x^\star},
\]
where we let \(\varepsilon\geq0\) so as to include the case of nonstrong minimality of \(x^\star\).
From (Lipschitz-) continuity of \(\T\) and since \(\T(x^\star)=x^\star\), there exists a neighborhood \(V_{x^\star}\) of \(x^\star\) such that \(\T(x)\in U_{x^\star}\) for all \(x\in V_{x^\star}\).
Then, using \cref{{prop:FBEfix},,{prop:FBEgeq}}, for all \(x\in V_{x^\star}\) we have
\begin{align*}
	\varphi_\gamma(x)-\varphi_\gamma(x^\star)
{}\geq{} &
	\varphi(\T(x))
	{}-{}
	\varphi(x^\star)
	{}+{}
	\tfrac{1-\gamma L_f}{2\gamma}
	\|x-\T(x)\|^2
\\
{}\geq{} &
	\tfrac \varepsilon2
	\|\T(x)-x^\star\|^2
	{}+{}
	\tfrac{1-\gamma L_f}{2\gamma}
	\|x-\T(x)\|^2;
\shortintertext{letting
	\(
		\varepsilon'
	{}={}
		\min\set{
			\varepsilon,
			\frac{1-\gamma L_f}{\gamma}
		}
	\),
	which is positive iff \(\varepsilon\) is,
}
{}\geq{} &
	\tfrac{\varepsilon'}{4}
	\|x-x^\star\|^2
\end{align*}
where in the last inequality we used the fact that
\(
	\|u\|^2+\|v\|^2
{}\geq{}
	\frac 12
	\|u-v\|^2
\)
for any \(u,v\in\R^n\).
This concludes the proof of \ref{thm:FBELocMin}.
}%
\end{proof}
\end{thm}

\arxiv{%
	Since prox-regularity is enjoyed globally by convex functions, the following special case is a straightforward consequence.
	\unskip
	\begin{cor}[First-order properties for convex \(g\)]\label{cor:1stOrderConvex}%
	Suppose that \(g\) is convex and that \(f\in\cont[2](\R^n)\) (resp. \(f\in\cont[2+](\R^n)\)).
	Then, for all \(\gamma>0\) all the properties in \Cref{thm:1stOrder} hold globally (\ie for all \(x^\star\in\R^n\) with \(U_{x^\star}=\R^n\)).
	\end{cor}
}%
When \(f=0\), \Cref{thm:1stOrder} restates the known fact that if \(g\) is prox-regular at \(x^\star\) for \(0\in\partial g(x^\star)\), then \(g^\gamma\) is continuosly differentiable around \(x^\star\) with \(\nabla g^\gamma(x)=\frac 1\gamma(x-\prox_{\gamma g}(x))\).
Notice that the bound \(\gamma<\Gamma(x^\star)\) is tight: in general, for \(\gamma=\Gamma(x^\star)\) no continuity of \(\T\) nor continuous differentiability of \(\varphi_\gamma\) around \(x^\star\) can be guaranteed.
In fact, even when \(x^\star\) is \(\Gamma(x^\star)\)-critical, \(\T\) might even fail to be single-valued and \(\varphi_\gamma\) differentiable at \(x^\star\), as the following counterexample shows.
\begin{es}[Why \(\gamma\neq\Gamma(x^\star)\) in first-order properties]
Consider \(f=\frac 12x^2\) and \(g=\indicator_S\) where \(S=\set{0,1}\).
Then, \(L_f=1\), \(\gamma_g=+\infty\), \(\T(x)=\proj_S((1-\gamma)x)\) and the FBE is
\(
	\varphi_\gamma(x)
{}={}
	\tfrac{1-\gamma}{2}
	\|x\|^2
	{}+{}
	\tfrac{1}{2\gamma}\dist((1-\gamma)x,S)^2
\).
At the critical point \(x=1\), which satisfies \(\Gamma(1)=\nicefrac 12\), \(g\) is prox-regular for any subgradient\siam{. }%
\arxiv{%
	(cf. \Cref{prop:DistantDomainedPR} with \(g_1=\indicator_{\set{0}}\) and \(g_2=\indicator_{\set{1}}\)).
}%
For any \(\gamma\in(0,\nicefrac12)\) it is easy to see that \(\varphi_\gamma\) is differentiable in a neighborhood of \(x=1\).
However, for \(\gamma=\nicefrac12\) the distance function has a first-order singularity in \(x=1\), due to the \(2\)-valuedness of \(\T(1)=\proj_S(\nicefrac12)=\set{0,1}\).
\end{es}

\begin{es}[Prox-nonregularity of critical points]\label{es:NonProxRegular}%
Consider \(\varphi=f+g\) where \(f(x)=\tfrac12x^2\), \(g(x)=\indicator_S(x)\) and \(S=\set{\nicefrac 1n}[n\in\N_{\geq 1}]\cup\set{0}\).
For \(x_0=0\) we have \(\Gamma(x_0)=+\infty\), however \(g\) fails to be prox-regular at \(x_0\) for \(v_0=0=-\nabla f(x_0)\).
\switch{
	For any \(\rho>0\) and for any neighborhood \(V\) of \((0,0)\) in \(\graph g\) it is always possible to find a point arbitrarily close to \((0,-\nicefrac1\rho)\) with multi-valued projection on \(V\).
	Specifically, the midpoint
	\(
		P_n
	{}={}
		\bigl(
			\frac12
			(\frac 1n + \frac{1}{n+1})
			,\,
			{}-\nicefrac 1\rho
		\bigr)
	\)
	has 2-valued projection on \(\graph g\) for any \(n\in\N_{\geq 1}\), being it \(\proj_{\graph g}(P_n)=\set{\nicefrac1n,\nicefrac{1}{n+1}}\).
	By considering a large \(n\), \(P_n\) can be made arbitrarily close to \((0,-\nicefrac1\rho)\) and at the same time its projection(s) arbitrarily close to \((0,0)\).
	Therefore, \(g\) cannot be prox-regular at \(0\) for \(0\), for otherwise such projections would be single-valued close enough to \((0,0)\) \cite[Cor. 3.4 and Thm. 3.5]{poliquin1996prox-regular}.
}{
	This can be easily seen with the geometrical interpretation previously described, in that for any \(\rho>0\) and for any neighborhood \(V\) of \((0,0)\) in \(\graph g\) it is always possible to find a point arbitrarily close to \((0,-\nicefrac1\rho)\) with multi-valued projection on \(V\).
	Specifically, the midpoint
	\(
		P_n
	{}={}
		\bigl(
			\frac12
			(\frac 1n + \frac{1}{n+1})
			,\,
			{}-\nicefrac 1\rho
		\bigr)
	\)
	has 2-valued projection on \(\graph g\) for any \(n\in\N_{\geq 1}\), being it \(\proj_{\graph g}(P_n)=\set{\frac1n,\frac{1}{n+1}}\).
	By considering a large \(n\), \(P_n\) can be made arbitrarily close to \((0,-\nicefrac1\rho)\) and at the same time its projection(s) arbitrarily close to \((0,0)\).
}%
As a result, \(g^\gamma(x)=\frac{1}{2\gamma}\dist(x,S)^2\) is not differentiable around \(x=0\), and indeed at each midpoint \(\frac12(\frac1n+\frac{1}{n+1})\) for \(n\in\N_{\geq1}\) it has a nonsmooth spike.
\arxiv{%
	However, it can be easily seen that all the other critical points, \ie the ones in \(S\setminus\set{0}\), are prox-regular, and with basic calculus it can be verified that \(g^\gamma\) and \(\varphi_\gamma\) are differentiable at \(0\) as ensured by \Cref{{prop:criticalDifferentiable},{prop:FBEcriticDiff}}.
}%
\end{es}

To underline how unfortunate the situation depicted in \Cref{es:NonProxRegular} is, notice that adding a linear term \(\lambda x\) to \(f\) for any \(\lambda\neq 0\), yet leaving \(g\) unchanged, restores the desired prox-regularity of each critical point.
Indeed, this is trivially true for any nonzero critical point; besides, \(g\) is prox-regular at \(0\) for any \(\lambda\in(0,-\infty)\), and for \(\lambda<0\) we have that \(0\) is nomore critical.
\arxiv{%
	The reason why prox-regularity fails to hold in the above example is due to the density of isolated points close to \(0\).
	\begin{prop}[Prox-regularity of \emph{convex-cluster-domained} functions]\label{prop:DistantDomainedPR}%
	Let \(g=\min_{i\in I}g_i\) for a family of pairwise disjoint-domained proper, closed, lower semicontinuous and convex functions \(\set{g_i}[i\in I]\).
	Suppose that
	\[
		\mesh(g)
	{}:={}
		\inf\set{\dist(\dom g_i, \dom g_j)}{i,j\in I,~i\neq j}
	{}>{}
		0,
	\]
	then \(g\) is everywhere prox-regular for all subgradients, with constants \(r=0\) and \(\varepsilon=\frac12\mesh(g)\) as in \Cref{def:ProxRegularity}.
	\begin{proof}
	Let \((x_0,v_0)\in\graph\partial g\), then there exists a unique \(i\in I\) such that \((x_0,v_0)\in\graph\partial g_i\).
	By convexity we have
	\[
		g_i(x')
	{}\geq{}
		g_i(x)
		{}+{}
		\innprod{v}{x'-x}
	\qquad
		\forall x,x'\in\R,~v\in\partial g_i(x).
	\]
	The inequality above still holds if we replace \(g_i\) with \(g\) and constrain \(x\) and \(x\) in the ball \(\ball{x_0}{\frac12\mesh(g)}\), which implies \eqref{eq:ProxRegularity} with \(\rho=0\) and \(\varepsilon=\frac12\mesh(g)\).
	\end{proof}
	\end{prop}
	A practical employment of \(g\) as in \Cref{prop:DistantDomainedPR} concerns problems where the optimization variable is constrained in the disjoint union of convex sets with nonvanishing mesh.
	Possibly the most appropriate example regards binary problems where \(g=\indicator_{\set{0,1}}\) or integer programming, in which case \(g=\indicator_\Z\), where no matter what \(f\in\cont[1,1](\R^n)\) is, all critical points of \(\varphi=f+g\) are prox-regular.
}%

	\subsection{Second-order properties}
		\label{subsec:FBEsecond}
		In this section we discuss sufficient conditions for twice-differentiability of the FBE at critical points.
Additionally to prox-regularity, which is needed for local continuous differentiability, we will also need generalized second-order properties of \(g\).
The interested reader is referred to \cite[§13]{rockafellar2011variational} for an extensive discussion on \emph{epi-differentiability}.
\begin{ass}\label{ass:fg2}%
With respect to a given critical point \(x^\star\)
\begin{enumerate}
	\item\label{ass:f2}%
		\(\nabla^2 f\) exists and is (strictly) continuous around \(x^\star\);%
	\item\label{ass:g2}%
		\(g\) is prox-regular and (strictly) twice epi-differentiable at \(x^\star\) for \(-\nabla f(x^\star)\), with its second order epi-derivative being generalized quadratic:
	\begin{equation}\label{eq:GenQuadSecondEpiDer}
		\twiceepi[-\nabla f(x^\star)]{g}{x^\star}[d]
	{}={}
		\innprod{d}{Md} + \delta_S(d),
	\quad
		\forall d\in\R^n
	\end{equation}
	where \(S\subseteq \R^n\) is a linear subspace and \(M\in\R^{n\times n}\).
	Without loss of generality we take \(M\) symmetric, and such that \(\im(M)\subseteq S\) and \(\ker(M)\supseteq S^\perp\).%
	\footnote{%
		This can indeed be done without loss of generality: if \(M\) and \(S\) satisfy \eqref{eq:GenQuadSecondEpiDer}, then it suffices to replace \(M\) with \(M'=\proj_S\frac{M+\trans M}{2}\proj_S\) to ensure the desired properties.%
	}%
\end{enumerate}
We say that the assumptions are ``strictly'' satisfied if the stronger conditions in parenthesis hold.
\end{ass}
\arxiv{%
	The properties of \(M\) in \Cref{ass:g2} cause no loss of generality.
	Indeed, letting \(\Pi_S\) denote the orthogonal projector onto the linear space \(S\) (which is symmetric \cite{bernstein2009matrix}), if \(M\) satisfies \eqref{eq:GenQuadSecondEpiDer}, then so does \(M'=\Pi_S[\tfrac{1}{2}(M+\trans M)]\Pi_S\), which has the wanted properties.
}%

Twice epi-differentiability of \(g\) is a mild requirement, and cases where \(\twiceepi{g}{}\) is generalized quadratic are abundant \cite{rockafellar1988first,rockafellar1989second,poliquin1992nonsmooth,poliquin1995second}.
Moreover, prox-regular and \(C^2\)-partly smooth functions \(g\) (see \cite{lewis2002active,daniilidis2006geometrical}) comprise a wide class of functions that strictly satisfy \Cref{ass:g2} at a critical point \(x^\star\) provided that \emph{strict complementarity} holds, namely if \(-\nabla f(x^\star)\in\relint\partial g(x^\star)\).
In fact, it follows from \cite[Thm. 28]{daniilidis2006geometrical} applied to the \emph{tilted} function \(\tilde g=g+\innprod{\nabla f(x^\star)}{\cdot}\) (which is still \(C^2\)-partly smooth and prox-regular at \(x^\star\) \cite[Cor. 4.6]{lewis2002active}, \cite[Ex. 13.35]{rockafellar2011variational}) that \(\prox_{\gamma\tilde g}\) is continuously differentiable around \(x^\star\) for \(\gamma\) small enough (in fact, for \(\gamma<\Gamma(x^\star)\)).
From \cite[Thm 4.1(g)]{poliquin1996generalized} we then obtain that \(\tilde g\) is strictly twice epi-differentiable at \(x^\star\) with generalized quadratic second-order epiderivative, and the claim follows by \emph{tilting} back to \(g\).

We now show that the quite common properties required in \Cref{ass:fg2} are all is needed for ensuring first-order properties of the proximal mapping and second-order properties of the FBE at critical points.
\begin{thm}[Twice differentiability of \(\varphi_\gamma\)]\label{thm:2ndOrder}%
Suppose that \Cref{ass:fg2} is (strictly) satisfied with respect to a critical point \(x^\star\).
Then, for any \(\gamma\in(0,\Gamma(x^\star))\)
\def\myVar{\displaystyle\jac{\prox_{\gamma g}}{\Fw{x^\star}};}%
\def\myVarOne{\displaystyle\nabla^2\varphi_\gamma(x^\star).}%
\begin{enumerate}
\item\label{lem:JP}%
	\(\prox_{\gamma g}\) is (strictly) differentiable at \(\Fw{x^\star}\) with symmetric and positive semidefinite Jacobian
	\begin{equation}\label{eq:Jprox}
		\fillwidthof[r]{\myVarOne}{
			P_\gamma(x^\star)
		}
	{}\coloneqq{}
		\fillwidthof[l]{\myVar}{
			\jac{\prox_{\gamma g}}{\Fw{x^\star}};
		}
	\end{equation}
\item\label{lem:RStrictDiff}%
	\(\Res\) is (strictly) differentiable at \(x^\star\) with Jacobian
	\begin{equation}\label{eq:JacR}
		\fillwidthof[r]{\myVarOne}{
			\jac{\Res}{x^\star}
		}
	{}={}
		\fillwidthof[l]{\myVar}{
			\tfrac 1\gamma\left[\id-P_\gamma(x^\star) Q_\gamma(x^\star)\right],
		}
	\end{equation}
	where \(Q_\gamma\) is as in \eqref{eq:Q} and \(P_\gamma\) as in \eqref{eq:Jprox};
\item\label{thm:FBEHess}%
	\(\varphi_\gamma\) is (strictly) twice differentiable at \(x^\star\) with symmetric Hessian
	\begin{equation}\label{eq:FBEHess}
		\fillwidthof[r]{\myVarOne}{
			\nabla^2\varphi_\gamma(x^\star)
		}
	{}={}
		\fillwidthof[l]{\myVar}{
			Q_\gamma(x^\star)
			\jac{\Res}{x^\star}.
		}
	\end{equation}
\end{enumerate}
\begin{proof}
See \Cref{proof:thm:2ndOrder}.
\end{proof}
\end{thm}
Again, when \(f\equiv 0\) \Cref{thm:2ndOrder} covers the differentiability properties of the proximal mapping (and consequently the second-order properties of the Moreau envelope, due to the identity \(\nabla g^\gamma(x)=\frac 1\gamma(x-\prox_{\gamma g}(x))\)) as discussed in \cite{poliquin1996generalized}.

We now provide a key result that links nonsingularity of the Jacobian of the forward-backward residual \(\Res\) to strong (local) minimality for the original cost \(\varphi\) and for the FBE \(\varphi_\gamma\), under the generalized second-order properties of \Cref{ass:fg2}.
\begin{thm}[Conditions for strong local minimality]\label{thm:StrongMinimality}%
Suppose that \Cref{ass:fg2} is satisfied with respect to a critical point \(x^\star\), and let
\(
	\gamma
{}\in{}
	(0,\min\set{\Gamma(x^\star),\nicefrac{1}{L_f}})
\).
The following are equivalent:
\begin{enumerateq}
	\item\label{thm:StrLocMin}\label{thm:StrongMinimPsi}%
		\(x^\star\) is a strong local minimum for \(\varphi\);
	\item\label{thm:LocMin}\label{thm:JRnonsing}%
		\(x^\star\) is a local minimum for \(\varphi\) and \(J\Res(x^\star)\) is nonsingular;
	\item\label{thm:FBEHessPosDef}\label{thm:HessDef+}%
		the (symmetric) matrix \(\nabla^2 \varphi_\gamma(x^\star)\) is positive definite;
	\item\label{thm:StrLocMinFBE}\label{thm:StrongMinimFBE}%
		\(x^\star\) is a strong local minimum for \(\varphi_\gamma\);
	\item\label{thm:LocMinFBE}%
		\(x^\star\) is a local minimum for \(\varphi_\gamma\) and \(J\Res(x^\star)\) is nonsingular.
\end{enumerateq}
\begin{proof}
See \Cref{proof:thm:StrongMinimality}.
\end{proof}
\end{thm}

\section{\zerofpr{} algorithm}
	\label{sec:Algorithm}
	The first algorithmic framework exploiting the FBE for solving composite minimization problems was studied in \cite{patrinos2013proximal}, and other schemes have been recently investigated in \cite{stella2017forward,liu2017further}.
All such methods tackle the problem by looking for a (local) minimizer of the FBE, exploting the equivalence of (local) minimality for the original function \(\varphi\) and for the FBE \(\varphi_\gamma\), for \(\gamma\) small enough.
To do so, they all employ the concept of directions of descent, thus requiring the gradient of the FBE to be well defined everywhere.
In the more general framework addressed in this paper, such basic requirement is not met, which is why we approach the problem from a different perspective.
This leads to \refZ[], the first algorithm, to the best of our knowledge, that despite requiring only the black-box oracle of FBS and being suited for fully nonconvex problems it achieves superlinear convergence rates.
{\renewcommand\thealgorithm{\zerofpr}%
	\begin{algorithm*}[tb]%
		\algcaption{generalized forward-backward with nonmonotone linesearch}%
		\label{alg:zerofpr}%
		\begin{algorithmic}[1]
\Require{%
	\(\gamma\in(0,\min\set{\nicefrac{1}{L_f},\gamma_g})\),~
	\(\beta,p_{\rm min}\in(0,1)\),~
	\(\sigma\in(0,\gamma\frac{1-\gamma L_f}{2})\),~
	\(x^0\in\R^n\).%
}%
\Initialize{%
	\(\bar\Phi_0=\varphi_{\gamma}(x^0)\),~
	\(k=0\).%
}%
\STATE\label{state:zerofpr:Initial}%
	Select \(\bar x^k\in T_{\gamma}(x^k)\) ~and set~ \(r^k=\tfrac 1\gamma(x^k-\bar x^k)\)%
\STATE\label{state:zerofpr:Stop}%
	{\bf if}~ \(\|r^k\|=0\), ~{\bf then}~ stop; ~{\bf end if}%
\STATE\label{state:zerofpr:d}%
	Select a direction \(d^k\in\R^n\)%
\STATE\label{state:zerofpr:LS}%
	Let \(\tau_k\in\set{\beta^m}[m\in\N]\) be the smallest such that
	\(
		x^{k+1}
	{}={}
		\bar x^k + \tau_kd^k
	\)
	\begin{equation}\label{eq:LS}
			\varphi_{\gamma}(x^{k+1})
		{}\leq{}
			\bar\Phi_k
			{}-{}
			\sigma
			\|r^k\|^2
	\end{equation}
	\vspace{-1.25\baselineskip}%
\renewcommand\myVar{\bar\Phi_{k+1}}%
\STATE\label{state:zerofpr:nm}%
	\(
		\fillwidthof[l]\myVar{\bar\Phi_{k+1}}
	{}={}
		(1-p_k)\bar\Phi_k
		{}+{}
		p_k\varphi_{\gamma}(x^{k+1})
	\)~
	for some \(p_k\in[p_{\rm min},1]\)%
\item[]%
	\(k\gets k+1\)~ and go to \cref{state:zerofpr:Initial}.%
\end{algorithmic}

	\end{algorithm*}%
}%
	\subsection{Overview}
		Instead of directly addressing the minimization of \(\varphi\) or \(\varphi_\gamma\), we seek solutions of the following nonlinear inclusion (\emph{generalized equation})
\begin{equation}\label{eq:ProblemR}
	\text{find }~
	x^\star\in\R^n
\quad\text{such that}\quad
	0
{}\in{}
	\Res(x^\star).
\end{equation}
By doing so we address the problem from the same perspective of FBS, that is, finding fixed points of the forward-backward operator \(\T\) or, equivalently, zeros of its residual \(\Res\).
Despite \(\Res\) might be quite irregular when \(g\) is nonconvex, it enjoys favorable properties at the very solutions to \eqref{eq:ProblemR} --- \ie at \(\gamma\)-critical points --- starting from single-valuedness, cf. \Cref{prop:SingleValuedFB}.
If mild assumptions are met, \(\Res\) turns out to be continuous around and even differentiable at critical points (cf. \Cref{thm:1stOrder,thm:2ndOrder}), and as a consequence the \emph{inclusion} problem \eqref{eq:ProblemR} reduces to a well behaved system of \emph{equations}, as opposed to \emph{generalized equations}, when close to solutions.

This motivates addressing problem \eqref{eq:ProblemR} with fast methods for nonlinear equations.
Newton-like schemes are iterative methods that prescribe updates of the form
\begin{equation}\label{eq:Fastd}
	x^+=x-H\Res(x)
\end{equation}
which essentially amount to selecting \(H=H(x)\), a linear operator that ideally carries information of the geometry of \(\Res\) around \(x\), in the attempt to yield an optimal iterate \(x^+\).
For instance, when \(\Res\) is sufficiently regular Newton method corresponds to selecting \(H\) as the inverse of an element of the generalized Jacobian of \(\Res\) at \(x\), enabling fast convergence when close to a solution under some assumptions.
However, selecting \(H\) as in Newton method would require information additional to the forward-backward oracle \(\T\), and as such it goes beyond the scope of the paper.
For this reason we focus instead on \emph{quasi-Newton} schemes, in which \(H\) are linear operators recursively defined with low-rank updates that satisfy the \emph{(inverse) secant condition}
\begin{equation}\label{eq:Secant}
	H^+y = s,
\quad\text{where}\quad
	s=x^+-x
	~~\text{and}~~
	y\in \Res(x^+)-\Res(x).
\end{equation}
A famous result \cite{dennis1974characterization} states that, under mild assumptions and starting sufficiently close to a solution \(x^\star\), updates as in \eqref{eq:Fastd} are superlinearly convergent to \(x^\star\) iff the \emph{Dennis-Moré condition} holds, namely the limit
\(
	\frac{\|(H^{-1}-J\Res(x^\star))s\|}{\|s\|}
{}\to{}
	0
\).
More recently, in \cite{dontchev2012generalizations} the result was extended to generalized equations of the form \(f(x)+G(x)\ni 0\), where \(f\) is smooth and \(G\) possibly set-valued.
The study focuses on Josephy-Newton methods where the update \(x^+\) is the solution of the inner problem \(f(x)-Bx\in Bx^++G(x^+)\), where \(B=H^{-1}\), which can be interpreted as a forward-backward step in the metric induced by \(B\).
In particular, differently from the here proposed \refZ[], the method in \cite{dontchev2012generalizations} has the crucial limitation that, unless the operator \(B\) has a very particular structure, the \emph{backward} step \((B+G)^{-1}\) may be prohibitely challenging.
%

		\subsubsection{Globalization strategy}
			Quasi-Newton schemes are
\Newswitch{
	extremely handy and widely used methods.
}{
	powerful and widely used methods.
}%
However, it is well known that they are effective only when close enough to a solution and might even diverge otherwise.
To cope with this crucial downside there comes the need of a globalization strategy; this is usually addressed by means of a linesearch over a suitable \emph{merit function} \(\psi\), along directions of descent for \(\psi\) so as to ensure sufficient decrease for small enough stepsizes.
Unfortunately, the potential choice \(\psi(x)=\frac12\|\Res(x)\|^2\) is not regular enough for a `direction of descent' to be everywhere defined.
The proposed \Cref{alg:zerofpr} bypasses this limitation by exploiting the favorable properties of the FBE.

Globalizing the convergence of any fast local method is the core contribution of \refZ[], an algorithm that exploits the favorable properties of the FBE, and that requires exactly the same oracle of FBS.
Conceptually, \refZ[] is really elementary; for simplicity, let us first consider the monotone case, \ie with \(p_k\equiv 1\) so that \(\bar\Phi_k=\varphi_\gamma(x^k)\) (cf. \cref{state:zerofpr:nm}).
The following steps are executed for updating
\newsiam{the }%
iterate \(x^k\):
\begin{enumerate}
\item
	first, at \cref{state:zerofpr:Initial} a nominal forward-backward call yields an element \(\bar x^k\in \T(x^k)\) that decreases the value of \(\varphi_\gamma\) by at least \(\gamma\frac{1-\gamma L_f}{2}\|r^k\|^2\) (\cref{prop:FBEleq});
\item
	then, at \cref{state:zerofpr:d} an update direction \(d^k\) at \(\bar x^k\) (\emph{not} at \(x^k\)!) is selected;
\item
	because of the sufficient decrease \(x^k\mapsto\bar x^k\) on \(\varphi_\gamma\) and the continuity of \(\varphi_\gamma\), at \cref{state:zerofpr:LS} a stepsize \(\tau_k\) can be found with finite many backtrackings \(\tau_k\gets\beta\tau_k\) that ensures a decrease for \(\varphi_\gamma\) of at least \(\sigma\|r^k\|^2\) in the update \(x^k\mapsto\bar x^k+\tau_kd^k\), for any \(\sigma<\frac{1-\gamma L_f}{2}\).
\end{enumerate}
In order to reduce the number of backtrackings, \(p_k<1\) can be selected resulting in a \emph{nonmonotone} linesearch.
The sufficient decrease is enforced with respect to a parameter \(\bar\Phi_k\geq\varphi_\gamma(x^k)\) (cf. \cref{lem:tau}), namely a convex combination of \(\set{\varphi_\gamma(x^i)}_{i=0}^k\).
For the sake of convergence, \(\seq{p_k}\) can be selected arbitrarily in \((0,1]\) as long as it is bounded away from \(0\), hence the role of the user-set lower bound \(p_{\rm min}\).
Consequently, small values of \(\sigma\) and \(p_k\) concur in reducing conservatism in the linesearch by favoring larger stepsizes. \begin{lem}[{Nonmonotone linesearch globalization}]\label{lem:tau}%
For all \(k\in\N\) the iterates generated by \refZ[] satisfy
\begin{equation}
\label{eq:barPhi}
	\varphi_{\gamma}(\bar x^k)
{}\leq{}
	\varphi(\bar x^k)
{}\leq{}
	\varphi_{\gamma}(x^k)
{}\leq{}
	\bar\Phi_k
\end{equation}
and there exists \(\bar\tau_k>0\) such that
\begin{equation}\label{eq:NMLS}
	\varphi_\gamma(\bar x^k+\tau d^k)
{}\leq{}
	\bar\Phi_k
	{}-{}
	\sigma
	\|x^k-\bar x^k\|^2
\qquad
	\forall\tau\in[0,\bar\tau_k].
\end{equation}
In particular, the number of backtrackings at \cref{state:zerofpr:LS} is finite.
\begin{proof}
The first two inequalities in \eqref{eq:barPhi} are due to \cref{prop:FBEgeqGlobal}.
Moreover,
\[
	\bar\Phi_{k+1}
{}={}
	(1-p_k)\bar\Phi_k
	{}+{}
	p_k\varphi_{\gamma}(x^{k+1})
{}\geq{}
	(1-p_k)\varphi_{\gamma}(x^{k+1})
	{}+{}
	p_k\varphi_{\gamma}(x^{k+1})
{}={}
	\varphi_{\gamma}(x^{k+1}),
\]
where the inequality follows by the linesearch condition \eqref{eq:LS}; this proves the last inequality in \eqref{eq:barPhi}.
As to \eqref{eq:NMLS}, let \(k\) be fixed and contrary to the claim suppose that for all \(\varepsilon>0\) there exists \(\tau_\varepsilon\in[0,\varepsilon]\) such that the point
\(
	x_\varepsilon
{}={}
	\bar x^k+\tau_\varepsilon d^k
\)
satisfies
\(
	\varphi_{\gamma}(x_\varepsilon)
{}>{}
	\varphi_{\gamma}(x^k)
	{}-{}
	\sigma
	\|x^k-\bar x^k\|^2
\).
Taking the limit for \(\varepsilon\to 0^+\), continuity of \(\varphi_\gamma\) as ensured by \cref{prop:FBEC+} yields
\[
	\varphi_{\gamma}(\bar x^k)
{}\geq{}
	\varphi_{\gamma}(x^k)
	{}-{}
	\sigma
	\|x^k-\bar x^k\|^2
{}>{}
	\varphi_{\gamma}(x^k)
	{}-{}
	\gamma\tfrac{1-\gamma L_f}{2}
	\|x^k-\bar x^k\|^2
\]
where the last inequality is due to the fact that \(x^k\neq\bar x^k\).
This contradicts \cref{prop:FBEgeqGlobal}; therefore, there exists \(\bar\tau_k>0\) such that
\(
	\varphi_{\gamma}(\bar x^k+\tau d^k)
{}\leq{}
	\varphi_{\gamma}(x^k)
	{}-{}
	\sigma
	\|x^k-\bar x^k\|^2
\)
for all \(\tau\in[0,\bar\tau_k]\).
By combining this with \eqref{eq:barPhi} the claim follows.
\end{proof}
\end{lem}
\Cref{lem:tau} ensures that regardless of the choice of \(d^k\), \refZ[] does not get stuck in infinite loops.
In \Cref{sec:Convergence} we will also show that the algorithm returns solutions of problem \eqref{eq:ProblemR}, and that under mild assumptions at the limit point the convergence rate is superlinear when \emph{good} directions are selected at \cref{state:zerofpr:d}.
Before going through the technicalities, we briefly anticipate what such \emph{good} directions are.

		\subsubsection[Choice of the directions]{Choice of the directions: quasi-Newton methods}
			\label{sec:Directions}
			As already emphasized, fast convergence of \refZ[] will be obtained thanks to the employment of Newton-like directions \(d^k\).
Differently from the classical Newton-like step \eqref{eq:Fastd}, when stepsize \(1\) is accepted, the update in \refZ[] is of the form \(x^+=\bar x+d\) rather than \(x^+=x+d\), where \(\bar x\) is an element of \(\T(x)\).
Therefore, \(d\) needs to be a Newton-like direction at \(\bar x\), and \emph{not} at \(x\), namely
\begin{equation}\label{eq:-Hr}
	d^k=-H_k\bar r^k
\quad\text{for some }
	\bar r^k\in \Res(\bar x^k)
\end{equation}
(as opposed to \(\bar r^k\in \Res(x^k)\)).
\arxiv{%
	In the following subsections we propose Broyden's and BFGS quasi-Newton methods.
	Due to some technicalities, we will need some modifications to the recipe \eqref{eq:-Hr} in order to prove that BFGS enables the 	desired superlinear convergence.
}%

			\paragraph{Broyden's method}
				We consider a \emph{modified} Broyden's scheme \cite{powell1970numerical} that performs rank-one updates of the form
\begin{subequations}\label{subeq:Broyden}
\begin{equation}
	H_{k+1}
{}={}
	H_k
	{}+{}
	\frac{
		s_k-H_ky_k
	}{
		\tinnprod{s_k}{(\nicefrac{1}{\vartheta_k}-1)s_k+H_ky_k}
	}
	\trans{s_k}H_k
\quad\text{with}\quad
	\begin{cases}
		s_k={} & x^{k+1}-\bar x^k\\
		y_k={} & r^{k+1}-\bar r^k,
	\end{cases}
\end{equation}
for a sequence \(\seq{\vartheta_k}\subset(0,2]\).
The original Broyden formula \cite{broyden1965class} corresponds to selecting \(\vartheta_k\equiv 1\), whereas for other values of \(\vartheta_k\) the secant condition \eqref{eq:Secant} is \emph{drifted} to \(H^+\tilde y=s\), where \(\tilde  y=(1-\vartheta)H^{-1}s+\vartheta y\).
In particular, \cite{powell1970numerical} suggests
\begin{equation}\label{eq:PowellTheta}
	\vartheta_k
{}\coloneqq{}
	\begin{cases}[l@{~~\text{if }}l]
			1
		&
			|\gamma_k| \geq \bar\vartheta
		\\
			\frac{
				1-\sign(\gamma_k)
				\bar\vartheta
			}{
				1-\gamma_k
			}
		&
			|\gamma_k| < \bar\vartheta
	\end{cases}
\qquad\text{where}\quad
	\gamma_k
{}\coloneqq{}
	\frac{
		\tinnprod{H_ky^k}{s^k}
	}{
		\|s^k\|^2
	}
\end{equation}
\end{subequations}
and \(\bar\vartheta\in(0,1)\) is a fixed parameter, with the convention \(\sign 0=1\).
Starting from an invertible matrix \(H_0\) this
\newsiam{specific }%
selection ensures that all matrices \(H_k\) are invertible.

			\paragraph[BFGS method]{BFGS method}
				\switch{
	BFGS method consists
}{
	A direct application of BFGS method results
}%
in the following update rule for matrices \(H_k\) in \eqref{eq:-Hr}: starting from a symmetric and positive definite \(H_0\),
\begin{equation}\label{eq:BFGS}
	H_{k+1}
{}={}
	\bigl(I-\rho_ks_k\trans y_k\bigr)
	H_k
	\bigl(I-\rho_ky_k\trans s_k\bigr)
	{}+{}
	\rho_ks_k\trans s_k
\Newswitch{,\quad}{\!,~~}
	\rho_k
{}={}
	\begin{cases}[c @{~~} l]
		\tfrac{1}{\tinnprod{s_k}{y_k}} & \text{if } \tinnprod{s_k}{y_k}>0 \\
		0 & \text{otherwise,}
	\end{cases}
\end{equation}
with \(s_k=x^{k+1}-\bar x^k\) and \(y_k=r^{k+1}-\bar r^k\), see \eg \cite[§6.1]{nocedal2006numerical}.
BFGS is the most popular quasi-Newton scheme; it is based on rank-two updates that, additionally to the secant condition, enforce also symmetricity.
In fact, BFGS is guaranteed to satisfy the Dennis-Moré condition only provided that the Jacobian of the nonlinear system at the limit point is symmetric \cite{byrd1989tool}.
Although this is not the case for \(J\Res(x^\star)\), we observed in practice that BFGS directions \eqref{eq:BFGS} perform extremely well.
\arxiv{%
	For the sake of the theory, in the following section we propose a modification of \eqref{eq:-Hr} that provably enables superlinear convergence of \refZ[] with BFGS method.
	Such modification comes at the cost of requiring one extra evaluation of \(\nabla f\) per update; in fact, despite not supported by our theory extensive numerical evidence seems to confirm that BFGS directions as proposed in \eqref{eq:BFGS} are more efficient.
}%

			\arxiv{%
				\paragraph[Symmetrized BFGS method]{``Symmetrized'' BFGS method}
					In order to obtain a certificate of superlinear convergence with BFGS, we may \emph{symmetrize} \(\Res\) into the operator
	\begin{equation}\label{eq:tildeR}
		\tilde \Res(x)
	{}\coloneqq{}
		\tfrac 1\gamma
		\bigl[
			(\id-\gamma\nabla f)(x)
			{}-{}
			(\id-\gamma\nabla f)(\T(x))
		\bigr]
	{}\subseteq{}
		\Res(x)
		{}+{}
		\nabla f(\T(x))
		{}-{}
		\nabla f(x)
	\end{equation}
	and replace problem \eqref{eq:ProblemR} with the equivalent generalized equation
	\begin{equation}\label{eq:ProblemTildeR}
		\text{find }~
		x^\star\in\R^n
	\quad\text{such that}\quad
		0
	{}\in{}
		\tilde \Res(x^\star).
	\end{equation}
	Since \(\func{\id-\gamma\nabla f}{\R^n}{\R^n}\) is invertible for \(\gamma<\nicefrac{1}{L_f}\), it is apparent that \eqref{eq:ProblemR} and \eqref{eq:ProblemTildeR} are equivalent problems, in that \(\zer\Res=\zer\tilde\Res\) for all \(\gamma\in(0,\nicefrac{1}{L_f})\).
	Moreover, if \Cref{ass:fg2} is satisfied at a critical point \(x^\star\), then from the chain rule and \Cref{thm:2ndOrder} we have that \(\tilde \Res\) is differentiable at \(x^\star\) with Jacobian
	\[
		J\tilde \Res(x^\star)
	{}={}
		J\Res(x^\star)
		{}+{}
		\nabla^2f(x^\star)
		\bigl[
			J\T(x^\star)-I
		\bigr]
	{}={}
		Q_\gamma(x^\star)
		J\Res(x^\star)
	{}={}
		\nabla^2\varphi_\gamma(x^\star)
	\]
	which is indeed \emph{symmetric}, as it was shown in \Cref{thm:FBEHess}.

	For the iterates of \refZ[],
	\(
		r^k+\nabla f(\bar x^k)-\nabla f(x^k)
	\)
	is an element of
	\(
		\tilde \Res(x^k)
	\);
	similarly,
	\(
		\bar r^k+\nabla f(\bar x^k-\gamma\bar r^k)-\nabla f(\bar x^k)
	\)
	belongs to \(\tilde \Res(\bar x^k)\) for any \(\bar r^k\in \Res(x^k)\).
	We may then tackle \eqref{eq:ProblemTildeR} by selecting the following directions in \refZ[]
	\begin{subequations}\label{subeq:tildeBFGS}
	\begin{equation}
		d^k
	{}={}
		-H_k
		\bigl[
			\bar r^k+\nabla f(\bar x^k-\gamma\bar r^k)-\nabla f(\bar x^k)
		\bigr]
	\quad
		\text{for some \(\bar r^k\in \Res(\bar x^k)\)}
	\end{equation}
	where, starting from a symmetric and invertible \(H_0\), matrices \(H_k\) are recursively defined with the BFGS rank-two (inverse) update \eqref{eq:BFGS} with
	\begin{equation}
		\begin{cases}
			s_k={} & x^{k+1}-\bar x^k\\
			y_k={} & r^{k+1}+\nabla f(\bar x^{k+1})-\nabla f(x^{k+1})
				{}-{}
					\bigl(\bar r^k+\nabla f(\bar x^k-\gamma\bar r^k)-\nabla f(\bar x^k)\bigr).
		\end{cases}
	\end{equation}
	\end{subequations}
	Compared to the ``non symmetrized'' BFGS directions \eqref{eq:BFGS}, this selections requires the additional evaluation of \(\nabla f(\bar x^k-\gamma\bar r^k)\).

			}%
			\paragraph{Limited-memory variants}
				Ultimately, instead of storing and operating on dense \(m\times m\) matrices, \emph{limited-memory} variants of quasi-Newton schemes keep in memory only a few (usually \(3\) to \(20\)) most recent pairs \((s^k,y^k)\) implicitly representing the approximate inverse Jacobian.
Their employment considerably reduces storage and computations over the full-memory counterparts, and as such they are the methods of choice for large-scale problems.
The most popular limited-memory method is L-BFGS: based on BFGS, it efficiently computes matrix-vector products with the approximate inverse Jacobian using a \emph{two-loop recursion} procedure \cite{liu1989limited,nocedal1980updating,nocedal2006numerical}.

	\subsection{Connections with other methods}
		The first algorithmic framework exploiting the FBE was studied in \cite{patrinos2013proximal}, where two semismooth Newton methods were analyzed for convex \(f\) and \(g\) with \(f\in\cont[2,1](\R^n)\) (twice continuously differentiable with Lipschitz continuous gradient).
A generalization of the scheme was then studied in \cite{stella2017forward} under less restrictive assumptions, with particular attention to quasi-Newton directions in place of semismooth Newton methods.
The proposed algorithm interleaves descent steps over the FBE with forward-backward steps.
\cite{liu2017further} then analyzed global and linear convergence properties of a generic linesearch algorithmic framework for minimizing the FBE based on gradient-related directions, for analytic \(f\) and subanalytic, convex, and lower bounded \(g\).

Though apparently closely related,
the approach that we provide in this paper presents major conceptual differences from any of the ones above.
Apart from the significantly less restrictive assumptions, the crucial distinction is that our method is \emph{derivative-free}, \ie it does not require the gradient of the FBE.
As a consequence, no computation nor the existence of \(\nabla^2 f\) is required, resulting in a method that, differently from the others, truly relies on the very same oracle information of the forward-backward operator \(\T\).
\arxiv{%
	Furthermore, in our setting the update directions \(d^k\) can be \emph{arbitrary}, as no concept of ``descent'' is required.
	\Cref{table:Comparisons} offers a schematic overview of some differences between FBS and these three FBE-based generalizations.
	\begin{table}[h]
		\centering%
		\newcommand\unknown{{\bf\color{red!50!yellow}?}}%
		{\footnotesize%
\def\myVar#1{\textfillwidthof[c]{\minfbe{} \cite{stella2017forward}}{#1}}%
\setlength\tabcolsep{2pt}%
\begin{tabular}{@{}|l|c|c|c|c|@{}}
	\cline{2-5}
	\multicolumn{1}{l|}{\vphantom{\(\Big|\)}} & \myVar{FBS} & \myVar{Alg. \cite{liu2017further}} & \myVar{\minfbe{} \cite{stella2017forward}} & \myVar{\refZ[]}\\\hline
	\rowcolor{pink!50}
	\textbf{Superlinear convergence} & \xmark & \unknown & \cmark & \cmark\\
	\quad Nonmonotone linesearch & --- & \xmark & \xmark & \cmark\\
	\quad No need of \(d^k\) bounded {\scriptsize (Newton-like dir.)} & --- & \unknown & \cmark & \xmark \\[2pt]
	%
	\rowcolor{gray!10}
	Arbitrary update direction \(d^k\) & \xmark & \xmark & \xmark & \cmark\\
	Support for nonconvex \(g\) & \cmark & \xmark & \xmark & \cmark\\
	\rowcolor{gray!10}
	Support for \(f\notin C^2\) & \cmark & \xmark & \xmark & \cmark\\
	Global \(O(1/k)\) rate {\scriptsize (in convex case)} & \cmark & \unknown & \cmark & \xmark\\
	\rowcolor{gray!10}
	Only \(\prox_{\gamma g}\) and \(\nabla f\) evaluations & \cmark & \xmark & \xmark & \cmark\\
	Conditions on \(\varphi\) for global/linear converg. & \multirow{ 2}{*}{\cmark} & \multirow{ 2}{*}{\xmark} & \multirow{ 2}{*}{\cmark} & \multirow{ 2}{*}{\xmark}\\[-2pt]
	\textfillwidthof[c]{Conditions on \(\varphi\) for global/linear converg.}{\scriptsize (as opposed to conditions on \(\varphi_\gamma\))} &&&& \\
	\rowcolor{gray!10}
	Descent algorithm & \cmark & \cmark  & \cmark & \xmark\\
	\hline
\end{tabular}
}%

		\caption{%
			Comparison between features of plain FBS, the generic algorithm proposed in \cite{liu2017further}, \minfbe{} \cite{stella2017forward}, and \refZ[].
			Symbols `---' indicate a property that does not apply to the algorithm, whereas `\unknown{}' a feature that has not been investigated.%
		}%
		\label{table:Comparisons}%
	\end{table}
}%

	\subsection{Main remarks}
		In this section we list a few observations that come in handy when implementing \refZ[].
\begin{rem}[Adaptive variant when \(L_f\) is unknown]\label{rem:Adaptive}%
In practice, no prior knowledge of the global Lipschitz constant \(L_f\) is required for \refZ[].
In fact, replacing \(L_f\) with an initial estimate \(L>0\) and fixing a backtracking ratio \(\alpha\in(0,1)\), after \cref{state:zerofpr:Stop} the following instruction can be added:
\begin{algorithmic}[1]
\makeatletter%
	\renewcommand{\ALC@lno}{%
		\ALC@linenosize\arabic{ALC@line}bis\ALC@linenodelimiter%
	}%
\stepcounter{ALC@line}%
\makeatother%
\IF{
	\(f(\bar x^k) > f(x^k)-\tinnprod{\nabla f(x^k)}{x^k-\bar x^k}+\frac L2\|x^k-\bar x^k\|^2\)
}%
	\item[]%
		\(\gamma\gets\alpha\gamma\),~
		\(L\gets\nicefrac L\alpha\),~
		\(\sigma\gets\alpha\sigma\),~
		\(\bar\Phi_k \gets \varphi_\gamma(x^k)\)~
		and go to \cref{state:zerofpr:Initial}.%
\ENDIF{}%
\end{algorithmic}
Whenever the quadratic bound \eqref{eq:LipBound} is violated with \(L\) in place of \(L_f\), the estimated Lipschitz constant \(L\) is increased and \(\gamma\) decreased accordingly; as a consequence, the FBE \(\varphi_\gamma\) changes and the nonmonotone linesearch is restarted.
Since replacing \(L_f\) with any \(L\geq L_f\) still satisfies \eqref{eq:LipBound}, it follows that \(L\) is incremented only a finite number of times.
Therefore, there exists an iteration \(k_0\) starting from which \(\gamma\) and \(\sigma\) are constant; in particular, all the results of the paper remain valid starting from iteration \(k_0\), at latest.
\arxiv{%
	Notice that, due to \Cref{prop:Gamma}, it makes sense to insert the instruction after \cref{state:zerofpr:Stop} rather than after \cref{state:zerofpr:Initial}, since having \(x^k\) \(\gamma\)-critical automatically ensures it to be \(\gamma'\)-critical for any \(\gamma'<\gamma\).
}%
\end{rem}

\begin{rem}[Support for locally Lipschitz \(\nabla f\)]\label{rem:C1+}%
If \(\dom g\) is bounded and, as it is reasonable, the directions \(\seq{d^k}\) selected at \cref{state:zerofpr:d} do not diverge, then \Cref{ass:f} on \(f\) can be relaxed to \(\nabla f\) being \emph{locally} Lipschitz.

In fact, it follows from the definition of proximal mapping that \(\seq{\bar x^k}\subseteq\dom g\), and if the directions are bounded then there exists a compact domain \(\Omega\supseteq\dom g\) such that \(\seq{x^k}\subseteq\Omega\).
Then, all results of the paper apply by replacing \(L_f\) with \(\lip_\Omega\nabla f\), the (finite) Lipschitz constant of \(\nabla f\) on \(\Omega\).
\end{rem}

\begin{rem}[Cost per iteration]\label{rem:cost}%
Evaluating \(\varphi_\gamma\) essentially amounts to one evaluation of \(\T\); this is evident from the expression \eqref{eq:FBEMoreau}, together with the observation that
\(
	g^\gamma(\Fw x)
{}={}
	g(\bar x)+\tfrac{1}{2\gamma}\|\Fw x-\bar x\|^2
\)
for any \(\bar x\in \T(x)\).
Therefore, computing \(\varphi_\gamma(\bar x^k + \tau_k d^k)\) at \cref{state:zerofpr:LS} yields an element \(\bar x^{k+1}\in \T(x^{k+1})\) required in \cref{state:zerofpr:Initial}, since \(x^{k+1} = \bar x^k + \tau_k d^k\) at every iteration.
In general, one evaluation of \(\T\) per backtracking step is required.
If the directions \(d^k\) are computed with Broyden or BFGS methods \eqref{subeq:Broyden} and \eqref{eq:BFGS}, then one additional evaluation of \(\T\) is required for retrieving \(d^k\); in the best case of \(\tau_k = 1\) being accepted, which asymptotically happens under mild assumptions (cf. \cref{thm:BroydenBFGS}), the algorithm then requires exactly \emph{two} evaluations of \(\T\) per iteration.
\arxiv{%
	If instead ``symmetrized'' BFGS directions \eqref{subeq:tildeBFGS} are employed, then one additional evaluation of \(\T\) and one of \(\nabla f\) are required for retrieving \(d^k\), resulting in exactly \emph{two} evaluations of \(\T\) and one additional of \(\nabla f\) per iteration in case \(\tau_k = 1\).
}%
\end{rem}

\begin{rem}[Extension of FBS]
Observe that by selecting \(d^k\equiv 0\) \refZ[] reduces to the classical FBS algorithm.
\Cref{prop:FBEgeqGlobal} combined with the relation
\(
	\varphi_\gamma(x^k)\leq\bar\Phi_k
\)
due to \eqref{eq:barPhi} shows that the condition at
\switch{
	\cref{state:zerofpr:LS}
}{
	\cref{state:zerofpr_nonadaptive:LS}
}%
is always statisfied (with \(\tau_k=1\)).
Therefore,
\(
	x^{k+1}
{}={}
	\bar x^k+d^k
{}={}
	\bar x^k
{}\in{}
	\T(x^k)
\)
for all \(k\), which is FBS, cf. \eqref{eq:FB}.
\end{rem}

	\subsection{Convergence results}
		\label{sec:Convergence}
		In this section we analyze the properties of cluster points of the iterates generated by \refZ[].
Specifically,
\begin{itemize}[%
	itemsep=2pt,
	labelsep=5pt,
	leftmargin=20pt,
]%
\item
	every cluster point of \(\seq{x^k}\) and \(\seq{\bar x^k}\) solves problem \eqref{eq:ProblemR} (\Cref{thm:Critical});
\item
	if the linesearch is (eventually) monotone, then global and linear convergence are achieved under mild assumptions (\Cref{thm:Global,thm:Linear});
\item
	directions satisfying the Dennis-Moré condition, such as Broyden's,
	\arxiv{%
		and symmetrized BFGS,
	}%
	enable superlinear rates under mild assumptions (\Cref{thm:DM,thm:BroydenBFGS}).
\end{itemize}
In what follows, we exclude the trivial case in which the optimality condition \(r^k=0\) is achieved in a finite number of iterations, and therefore
\newsiam{we }%
assume \(r^k\neq 0\) for all \(k\).
\begin{thm}[Criticality of cluster points]\label{thm:Critical}%
The following hold for the iterates generated by \refZ[]:
\begin{enumerate}
	\item\label{thm:Critical:r}%
		\(r^k\to 0\) square-summably, and all cluster points of \(\seq{x^k}\) and \(\seq{\bar x^k}\) are critical; more precisely,
		\(
			\omega(x^k)=\omega(\bar x^k)\subseteq\fix \T
		\);
	\item\label{thm:Critical:Lim}%
		\(\seq{\varphi_\gamma(x^k)}\) converges to a (finite) value \(\varphi_\star\), and so does \(\seq{\varphi(\bar x^k)}\) if \(\seq{x^k}\) is bounded.
\end{enumerate}
\begin{proof}
\begin{proofitemize}
\item\ref{thm:Critical:r}:
	For all iterates \(k\) we have
	\begin{equation}\label{eq:SD}
		\bar\Phi_{k+1}
	{}={}
		(1-p_k)
		\bar\Phi_k
		{}+{}
		p_k
		\varphi_\gamma(x^{k+1})
	\smash{{}\overrel[\leq]{\eqref{eq:LS}}{}}
		\bar\Phi_k - \sigma p_k\|r^k\|^2
	{}\leq{}
		\bar\Phi_k - \sigma p_{\rm min}\|r^k\|^2.
	\end{equation}
	By telescoping the above inequality and using \eqref{eq:barPhi}, we obtain
	\begin{equation}\label{eq:Telescopic}
		\textstyle
		\bar\Phi_k
		{}-{}
		\inf\varphi
	{}\geq{}
		\bar\Phi_0
		{}-{}
		\bar\Phi_{k+1}
	{}={}
		\sum_{i=0}^k
		\left[
			\bar\Phi_i
			{}-{}
			\bar\Phi_{i+1}
		\right]
	{}\geq{}
		\sigma p_{\rm min}
		\sum_{i=0}^k
		{
			\|r^i\|^2
		},
	\end{equation}
	proving \(r^k\to 0\) square-summably.
	Suppose now that \(\seq{x^k}[k\in K]\to x'\) for some \(x'\in\R^n\) and \(K\subseteq\N\).
	Then, since
	\(
		\|\bar x^k-x^k\|
	{}={}
		\gamma\|r^k\|
	{}\to{}
		0
	\),
	in particular \(\seq{\bar x^k}[k\in K]\to x'\) as well.
	Due to the arbitrarity of the cluster point \(x'\) it follows that \(\omega(x^k)\subseteq\omega(\bar x^k)\), and a similar reasoning proves the converse inclusion, hence \(\omega(x^k)=\omega(\bar x^k)\).
	Moreover, we have
	\(
		x^k
	{}\in{}
		\cball{\bar x^k}{\gamma\|r^k\|}
	{}\subseteq{}
		\FB{x^k}
		{}+{}
		\cball{0}{\gamma\|r^k\|}
	\)
	and since
	\(
		\seq{
			x^k
			{}-{}
			\gamma\nabla f(x^k)
		}[k\in K]
	{}\to{}
		x'
		{}-{}
		\gamma\nabla f(x')
	\),
	from the outer semicontinuity of \(\prox_{\gamma g}\) \cite[Ex. 5.23(b)]{rockafellar2011variational} it follows that \(x'\in\FB{x'}\), \ie \(x'\in\fix \T\).
\item\ref{thm:Critical:Lim}:
	from \eqref{eq:SD} it follows that \(\seq{\bar\Phi_k}\) is decreasing, and in particular its limit exists, be it \(\varphi_\star\).
	Due to \eqref{eq:barPhi}, necessarily \(\varphi_\star\geq\inf\varphi>-\infty\), therefore
	\[
		0
	{}\leftarrow{}
		\bar\Phi_k-\bar\Phi_{k+1}
	{}={}
		p_k
		\left(
			\bar\Phi_k-\varphi_\gamma(x^{k+1})
		\right)
	{}\geq{}
		p_{\rm min}
		\left(
			\bar\Phi_{k+1}-\varphi_\gamma(x^{k+1})
		\right)
	{}\geq{}
		0,
	\]
	proving that \(\varphi_\gamma(x^{k+1})\to\varphi_\star\).
	If \(\seq{x^k}\) is bounded, then so is \(\seq{\bar x^k}\) due to compact-valuedness of \(\prox_{\gamma g}\) \cite[Thm. 1.25]{rockafellar2011variational}.
	Due to \cref{prop:FBEC+} \(\varphi_\gamma\) is \(L\)-Lipschitz continuous on a compact set containing \(\seq{x^k}\) and \(\seq{\bar x^k}\) for some \(L>0\).
	Then,
	\[
		0
	{}\leq{}
			\varphi_\gamma(x^k)-\varphi(\bar x^k)
	{}\leq{}
			\varphi_\gamma(x^k)-\varphi_\gamma(\bar x^k)
	{}\leq{}
		L\gamma\|r^k\|
	{}\to{}
		0
	\]
	where the inequalities follow from \cref{prop:Ineq}.
	Consequently, \(\seq{\varphi(\bar x^k)}\to\varphi_\star\) as well.
	\newarxiv{\qedhere}
\end{proofitemize}
\end{proof}
\end{thm}

		\subsubsection{Global and linear convergence}
If follows from \eqref{eq:barPhi} and the fact that \(\seq{\bar\Phi_k}\) is a decreasing sequence (cf. \eqref{eq:SD}), that the iterates of \refZ[] satisfy \(\varphi(\bar x^k)\leq\bar\Phi_0=\varphi(\bar x^0)\).
As a consequence, a sufficient condition for ensuring that the sequence \(\seq{\bar x^k}\) does not diverge --- and consequently nor does \(\seq{x^k}\) provided that the sequence of directions \(\seq{d^k}\) is bounded --- is that the level set \(\set{\varphi\leq\varphi(\bar x^0)}\) is compact.
In the adaptive variant discussed in \Cref{rem:Adaptive}, this translates to boundedness of the level set \(\set{\varphi\leq\varphi(\bar x^{k_0})}\), where \(k_0\) denotes the iteration starting from which \(\gamma\) is constant.
Since such point is unknown a priori, the sufficient condition needs be strengthened to \(\varphi\) having bounded level sets.

We now show that if \(\varphi_\gamma\) is well-behaved at cluster points, then the whole sequence generated by \refZ[] is convergent.
Good behavior involves the existence of a \emph{desingularizing function}, that is, \(\varphi_\gamma\) needs to possess the \emph{Kurdyka-\L ojasiewicz} property, a mild requirement that we restate here for the reader's convenience.
\begin{defin}[KL property]\label{def:KL}%
A proper and lower semicontinuous function \(\func{h}{\R^n}{\Rinf}\) has the Kurdyka-\L ojasiewicz property (KL property) at \(x^\star\in\dom\partial h\) if there exist a concave \emph{desingularizing function} (or \emph{KL function}) \(\func{\psi}{[0,\eta]}{[0,+\infty)}\) for some \(\eta>0\) and a neighborhood \(U_{x^\star}\) of \(x^\star\), such that
\begin{enumerate}
	\item\label{def:KL1}%
		\(\psi(0) = 0\);
	\item\label{def:KL2}%
		\(\psi\) is \(\cont[1]\) with \(\psi'>0\) on \((0,\eta)\);
	\item\label{def:KL3}%
		for all \(x\in U_{x^\star}\) s.t. \(h(x^\star) < h(x) < h(x^\star) + \eta\) it holds that
		\begin{equation}\label{eq:KL3}
			\psi'\bigl(h(x)-h(x^\star)\bigr)
			\dist\bigl(0,\partial h(x)\bigr)
		{}\geq{}
			1.
		\end{equation}
\end{enumerate}
\end{defin}
The KL property is a mild requirement enjoyed by semi-algebraic functions and by subanalytic functions which are continuous on their domain \cite{bolte2007lojasiewicz,bochnak1998real} see also \cite{lojasiewicz1963propriete,lojasiewicz1993geometrie,kurdyka1998gradients}.
Moreover, since semi-algebraic functions are closed under parametric minimization, from the expression \eqref{eq:FBE1} it is apparent that \(\varphi_\gamma\) is semi-algebraic provided that \(f\) and \(g\) are.
More precisely, in all such cases the desingularizing function can be taken of the form \(\psi(s)=\rho s^\theta\) for some \(\rho>0\) and \(\theta\in(0,1]\), in which case it is usually referred to as a \L ojasiewicz function.
This property has been extensively exploited to provide convergence rates of optimization algorithms such as FBS, see \cite{attouch2009convergence,attouch2010proximal,attouch2013convergence,bolte2014proximal,frankel2015splitting,ochs2014ipiano}.
Further properties of \(f\) and \(g\) that ensure \(\varphi_\gamma\) to satisfy such requirement are discussed in \cite{liu2017further}.

We first show how the KL property on \(\varphi_\gamma\) ensures global convergence of the iterates of \refZ[] if the linesearch is eventually monotone, \ie if \(p_k=1\) for \(k\) sufficiently large, and then show that linear convergence is attained when the KL function is actually a \L ojasiewicz function with large enough exponent.
\begin{thm}[Global convergence (monotone LS)]\label{thm:Global}%
Consider the iterates generated by \refZ[] with \(p_k=1\) for \(k\)'s large enough, and with directions satisfying
\begin{equation}\label{eq:dr}
	\|d^k\|\leq D\|r^k\|
\quad
	\text{for all }
	k
\end{equation}
for some \(D\geq 0\).
Suppose that \(\seq{x^k}\) remains bounded, that \(\varphi_\gamma\) has the KL property on \(\omega(x^k)\), and that every cluster point is prox-regular.
If \(f\) is of class \(C^2\) in a neighborhood of \(\omega(x^k)\), then \(\seq{x^k}\) and \(\seq{\bar x^k}\) are convergent to (the same point) \(x^\star\), and the sequence of residuals \(\seq{r^k}\) is summable.
\begin{proof}
	From \cref{lem:omega} we know that \(\varphi_\gamma\) is constant on the (nonempty) compact set \(\omega(x^k)\).
	It then follows from \cite[Lem. 6]{bolte2014proximal} that there exist \(\eta,\varepsilon>0\) and a \emph{uniformized KL function}, namely a function \(\psi\) satisfying \cref{{def:KL1},,{def:KL2},,{def:KL3}} for all \(x^\star\in\omega(x^k)\) and \(x\) such that
	\(
		\dist(x,\omega(x^k))<\varepsilon
	\)
	and
	\(
		\varphi(x^\star) < \varphi(x) < \varphi(x^\star) + \eta
	\).
	Let
	\(
		\varphi_\star
	{}\coloneqq{}
		\lim_{k\to\infty}\varphi_\gamma(x^k)
	\),
	which exists and is finite (cf. \cref{thm:Critical}), and let \(k_1\in\N\) be such that \(p_k=1\) for all \(k\geq k_1\).
	Then we have (cf. \cref{state:zerofpr:nm} and \eqref{eq:LS})
	\begin{equation}\label{eq:p1}
		\bar\Phi_k=\varphi_\gamma(x^k)
	\quad\text{and}\quad
		\varphi_\gamma(x^k)>\varphi_\gamma(x^{k+1})>\varphi_\star
	\qquad
		\forall k\geq k_1.
	\end{equation}
	By possibly restricting \(\varepsilon\), from \cref{thm:FBEC1} and since \(\omega(x^k)\) is compact it follows that \(\varphi_\gamma\) is differentiable in an \(\varepsilon\)-enlargement of \(\omega(x^k)\).
	\cref{lem:omega} ensures that there exists \(k_2\geq k_1\) such that for all \(k\geq k_2\) we have
	\(
		\varphi_\star
	{}<{}
		\varphi_\gamma(x^k)
	{}<{}
		\varphi_\star+\eta
	\)
	and
	\(
		\dist(x^k,\omega(x^k))
	{}<{}
		\varepsilon
	\).
	For all such \(k\), by \cref{thm:FBEC1} we have
	\(
		\nabla\varphi_\gamma(x^k)
	{}={}
		Q_\gamma(x^k)
		\Res(x^k)
	{}={}
		\bigl[I-\gamma\nabla^2f(x^k)\bigr]
		r^k
	\)
	and the uniformized KL property yields
	\begin{equation}\label{eq:Global:KL1}
		\textstyle
		\psi'
		\bigl(
			\varphi_\gamma(x^k)-\varphi_\star
		\bigr)
	{}\geq{}
		\frac{1}{\|\nabla\varphi_\gamma(x^k)\|}
	{}\geq{}
		\frac{1}{(1+\gamma L_f)\|r^k\|}.
	\end{equation}
	Letting
	\(
		\Delta_k
	{}\coloneqq{}
		\psi
		\bigl(
			\varphi_\gamma(x^k)-\varphi_\star
		\bigr)
	{}>{}
		0
	\),
	by concavity of \(\psi\) and \eqref{eq:p1} it follows that
	\begin{align}
	\nonumber
		\Delta_k
		{}-{}
		\Delta_{k+1}
	{}\geq{} &
		\psi'
		\bigl(
			\varphi_\gamma(x^k)-\varphi_\star
		\bigr)
		\bigl(
			\varphi_\gamma(x^k) - \varphi_\gamma(x^{k+1})
		\bigr)
	\\
	{}\geq{} &
		\frac{
			\varphi_\gamma(x^k) - \varphi_\gamma(x^{k+1})
		}{
			(1+\gamma L_f)\|r^k\|
		}
	{}={}
		\frac{
			\bar\Phi_k - \bar\Phi_{k+1}
		}{
			(1+\gamma L_f)\|r^k\|
		}
	{}\overrel[\geq]{\eqref{eq:SD}}{}
		\frac{
			\sigma p_{\rm min}
		}{
			1+\gamma L_f
		}
		\|r^k\|.
	\label{eq:Global:KL2}
	\end{align}
	By telescoping the inequality it follows that \(\seq{\|r^k\|}\) is summable, hence, due to \cref{lem:Deltaxr}, also \(\seq{\|x^{k+1}-x^k\|}\) is.
	Therefore, \(\seq{x^k}\) is a Cauchy sequence and as such it admits a limit, this being also the limit of \(\seq{\bar x^k}\) in light of \cref{thm:Critical:r} (and the fact that \(\seq{\bar x^k}\) is also bounded).
\end{proof}
\end{thm}
\begin{thm}[Linear convergence (monotone LS)]\label{thm:Linear}%
Consider the iterates generated by \refZ{}.
Suppose that the hypothesis of \Cref{thm:Global} are satisfied, and that the KL function can be taken of the form \(\psi(s)=\rho s^\theta\) for some \(\theta\in[\nicefrac12,1]\).
Then, \(\seq{x^k}\) and \(\seq{\bar x^k}\) are \(R\)-linearly convergent.
\begin{proof}
	From \cref{thm:Global} we know that \(\seq{x^k}\) and \(\seq{\bar x^k}\) converge to the same (\(\gamma\)-critical) point, be it \(x^\star\).
	Defining
	\(
		B_k
	{}\coloneqq{}
		\sum_{i\geq k}\|r^i\|
	\),
	from \cref{{lem:Deltaxr},{lem:DeltaBarxr}} we have
	\[
		\textstyle
		\|x^k-x^\star\|
	{}\leq{}
		\sum_{i\geq k}\|x^{i+1}-x^i\|
	{}\leq{}
		(\gamma+D)B_k
	\quad\text{and}\quad
		\|\bar x^k-x^\star\|
	{}\leq{}
		(3\gamma+D)B_k.
	\]
	Therefore, the proof reduces to showing that \(\seq{B_k}\) converges with asymptotic \(Q\)-linear rate.
	Inequality \eqref{eq:Global:KL1} reads
	\(
		\varphi_\gamma(x^k)
		{}-{}
		\varphi_\star
	{}\leq{}
		\bigl[
			(1+\gamma L_f)
			\rho\theta
			\|r^k\|
		\bigr]^{\frac{1}{1-\theta}}
	\),
	and since \(r^k\to 0\) for large enough \(k\) we have
	\[
		\newsiam{\Delta_k{}\coloneqq{}}
		\psi
		\bigl(
			\varphi_\gamma(x^k)-\varphi_\star
		\bigr)
	{}={}
		\rho
		[
			\varphi_\gamma(x^k)
			{}-{}
			\varphi_\star
		]^\theta
	{}\leq{}
		\rho
		\bigl[
			\smashoverbrace{
				(1+\gamma L_f)
				\rho
				\theta
				\|r^k\|
			}{
				<1~\text{for large }k
			}
		\bigr]^{\frac{\theta}{1-\theta}}
	{}\leq{}
		\rho^2
		(1+\gamma L_f)
		\|r^k\|.
	\]
	Therefore, eventually
	\(
		\Delta_k
		\newarxiv{{}\coloneqq{}\psi\bigl(\varphi_\gamma(x^k)-\varphi_\star\bigr)}
	{}<{}
		1
	\)
	and from \eqref{eq:Global:KL2} we get
	\[
		\textstyle
		B_k
	{}={}
		\sum_{i\geq k}\|r^i\|
	\smash{{}\overrel[\leq]{\eqref{eq:Global:KL2}}{}}
		\frac{1+\gamma L_f}{\sigma p_{\rm min}}
		\sum_{i\geq k}(\Delta_i-\Delta_{i+1})
	{}\leq{}
		\frac{1+\gamma L_f}{\sigma p_{\rm min}}
		\Delta_k
	{}\leq{}
		\frac{
			\rho^2
			(1+\gamma L_f)^2
		}{
			\sigma p_{\rm min}
		}
		\|r^k\|
	{}={}
		C\|r^k\|
	\]
	for some \(C>0\).
	Therefore, 	for large enough \(k\) we have
	\(
		B_k
	{}\leq{}
		C\|r^k\|
	{}={}
		C(B_k-B_{k+1})
	\),
	\ie
	\(
		B_{k+1}
	{}\leq{}
		(1-\nicefrac 1C)B_k
	\),
	proving asymptotic \(Q\)-linear convergence of \(B_k\).
\end{proof}
\end{thm}

		\subsubsection{Superlinear convergence}
			\label{sec:Superlinear}
			In the next result we show that under mild assumptions \refZ[] exhibits superlinear rates of convergence if the directions satisfy a Dennis-Mor\'e condition.
Then, we show that the Broyden scheme \eqref{subeq:Broyden} produces directions that satisfy such condition, and that due to the acceptance of unit stepsize \(\tau_k=1\), eventually each iteration of \refZ[] will require only two evaluations of \(\T\) (cf. \cref{rem:cost}).
We remind that
\switch{
	a sequence \(\seq{x^k}\) such that
}{
	\(\seq{x^k}\) with
}%
\(x^k\neq x^\star\) for all \(k\) is said to be \DEF{superlinearly convergent} to \(x^\star\) if
\switch{
	\(
			\|x^{k+1}-x^\star\|/
			\|x^k-x^\star\|
	{}\to{}
		0
	\)
	as \(k\to\infty\).
}{
	\[
		\lim_{k\to\infty}{
			\frac{
				\|x^{k+1}-x^\star\|
			}{
				\|x^k-x^\star\|
			}
		}
	{}={}
		0.
	\]
}%
\arxiv{
The following definition is taken from \cite[§7.5]{facchinei2007finite}.
\begin{defin}[Superlinearly convergent directions]%
Let \(\seq{x^k}\subset\R^n\) be a sequence converging to a point \(x^\star\).
We say that a sequence \(\seq{d^k}\subset\R^n\) is \DEF{superlinearly convergent with respect to \(\seq{x^k}\)} if
\[
	\lim_{k\to\infty}{
		\frac{
			\|x^k+d^k-x^\star\|
		}{
			\|x^k-x^\star\|
		}
	}
{}={}
	0.
\]
\end{defin}
}

\begin{thm}[Superlinear convergence under Dennis-Moré condition]\label{thm:DM}%
Suppose that \Cref{ass:fg2} is strictly satisfied at a strong local minimum \(x^\star\) of \(\varphi\), and consider the iterates generated by \refZ[].
Suppose that \(\seq{x^k}\) converges to \(x^\star\) and that the directions \(\seq{d^k}\) satisfy the Dennis-Moré condition
\begin{equation}\label{eq:DM}
	\lim_{k\to\infty}{
		\frac{
			\|
				\bar r^k
				{}+{}
				J\Res(x^\star)d^k
			\|
		}{
			\|d^k\|
		}
	}
{}={}
	0
\qquad
	\text{where \(\bar r^k\in \Res(\bar x^k)\).}
\end{equation}
Then, eventually stepsize \(\tau_k=1\) is always accepted and the sequences
\(\seq{x^k}\),
\(\seq{\bar x^k}\),
and \(\seq{r^k}\),
converge with superlinear rate.
\arxiv{%
	Moreover, all the claims still hold replacing \(\Res\) with \(\tilde \Res\).
}%
\begin{proof}
	From \cref{thm:1stOrder,lem:RStrictDiff,thm:StrongMinimality,thm:FBEHess} we know that \(\nabla\varphi_\gamma\) and \(\Res\) are strictly differentiable at \(x^\star\), with
	\(
		G_\star
	{}\coloneqq{}
		\nabla^2
		\varphi_\gamma(x^\star)
	{}={}
		Q_\gamma(x^\star)J\Res(x^\star)
	{}\succ{}
		0,
	\)
	and that there exists a neighborhood \(U_{x^\star}\) of \(x^\star\) in which \(\varphi_\gamma\) is differentiable and \(\Res\) Lipschitz continuous.
	Since \(\bar x^k=x^k-\gamma r^k\to x^\star\) due to \cref{thm:Critical:r}, it holds that \(x^k,\bar x^k\in U_{x^\star}\) for all \(k\) large enough.
	By single-valuedness of \(\Res\), for all such \(k\) we may write \(\Res(x^k)\) and \(\Res(\bar x^k)\) in place of \(r^k\) and \(\bar r^k\), respectively.
	In particular, since \(x^\star\in\fix \T\) (cf. \cref{thm:Critical:r}), necessarily \(\Res(\bar x^k)\to 0\).
	In turn, due to \eqref{eq:DM} it also holds that \(d^k\to 0\).
	Let \(x^{k+1}_0\coloneqq\bar x^k+d^k\); then,
	\newsiam{from \eqref{eq:DM} we have}
	\[
		\mathtight
		\newsiam{0{}\leftarrow{}}
		\frac{
			\Res(\bar x^k)
			{}+{}
			J\Res(x^\star)
			d^k
		}{
			\|d^k\|
		}
	{}={}
		\frac{
			\Res(\bar x^k)
			{}+{}
			J\Res(x^\star)
			(x^{k+1}_0-\bar x^k)
			{}-{}
			\Res(x_0^{k+1})
		}{
			\|x^{k+1}_0-\bar x^k\|
		}
		{}+{}
		\frac{
			\Res(x^{k+1}_0)
		}{
			\|x^{k+1}_0-\bar x^k\|
		}
	\]
	and since \(x^{k+1}_0-\bar x^k=d^k\to 0\), from
	\newarxiv{\eqref{eq:DM} and }%
	strict differentiability of \(\Res\) at \(x^\star\) applied on the first term on the right-hand side it follows that
	\begin{equation}\label{eq:Rzero}
		\lim_{k\to\infty}{
			\nicefrac{
				\|\Res(x^{k+1}_0)\|
			}{
				\|x^{k+1}_0-\bar x^k\|
			}
		}
	{}={}
		0.
	\end{equation}
	By possibly restricting \(U_{x^\star}\), nonsingularity of \(J\Res(x^\star)\) ensures the existence of a constant \(\alpha>0\) such that
	\(
		\|\Res(x)\|
	{}\geq{}
		\alpha
		\|x-x^\star\|
	\)
	for all \(x\in U_{x^\star}\).
	Since \(\bar x^k+d^k\to x^\star\), eventually \(x^{k+1}_0\in U_{x^\star}\).
	\Newswitch{From \eqref{eq:Rzero} we obtain}{We have}
	\[
		\newsiam{0{}\leftarrow{}}
		\frac{
			\|\Res(x^{k+1}_0)\|
		}{
			\|x^{k+1}_0-\bar x^k\|
		}
	{}\geq{}
		\alpha
		\frac{
			\|x^{k+1}_0-x^\star\|
		}{
			\|x^{k+1}_0-\bar x^k\|
		}
	{}\geq{}
		\alpha
		\frac{
			\|x^{k+1}_0-x^\star\|
		}{
			\|x^{k+1}_0-x^\star\|
			{}+{}
			\|\bar x^k-x^\star\|
		}
	{}={}
		\alpha
		\smash{
			\frac{
				\frac{
					\|x^{k+1}_0-x^\star\|
				}{
					\|\bar x^k-x^\star\|
				}
			}{
				1
				{}+{}
				\frac{
					\|x^{k+1}_0-x^\star\|
				}{
					\|\bar x^k-x^\star\|
				}
			}
		}
	\]
	and
	\Newswitch{therefore}{due to \eqref{eq:Rzero}}%
	\begin{equation}\label{eq:lim2}
		\lim_{k\to\infty}{
			\frac{
				\smash{\|x^{k+1}_0-x^\star\|}
			}{
				\|\bar x^k-x^\star\|
			}
		}
	{}={}
		\lim_{k\to\infty}{
			\frac{
				\smash{\|\bar x^k+d^k-x^\star\|}
			}{
				\|\bar x^k-x^\star\|
			}
		}
	{}={}
		0.
	\end{equation}
	A second-order expansion of \(\varphi_\gamma\) at \(x^\star\) yields
	{\mathtight%
	\begin{align*}
		\varphi_\gamma(\bar x^k)
	{}={} &
		\varphi_\gamma(x^\star)
		{}+{}
		\tfrac12
		\innprod{G_\star(\bar x^k-x^\star)}{\bar x^k-x^\star}
		{}+{}
		o(\|\bar x^k-x^\star\|^2)
	\shortintertext{and}
		\varphi_\gamma(\bar x^k+d^k)
	{}={} &
		\varphi_\gamma(x^\star)
		{}+{}
		\tfrac12
		\innprod{G_\star(\bar x^k+d^k-x^\star)}{\bar x^k+d^k-x^\star}
		{}+{}
		o(\|\bar x^k+d^k-x^\star\|^2)
	\\
	{}={} &
		\varphi_\gamma(x^\star)
		{}+{}
		o(\|\bar x^k-x^\star\|^2)
	\end{align*}%
	}%
	where the last equality is due to \eqref{eq:lim2}.
	Substracting,
	\begin{align*}
		\varphi_\gamma(\bar x^k+d^k)
		{}-{}
		\varphi_\gamma(\bar x^k)
	{}={} &
		{}-{}
		\tfrac12
		\innprod{G_\star(\bar x^k-x^\star)}{\bar x^k-x^\star}
		{}+{}
		o(\|\bar x^k-x^\star\|^2)
	\\
	{}\leq{} &
		-\beta\|\bar x^k-x^\star\|^2
		{}+{}
		o(\|\bar x^k-x^\star\|^2),
	\end{align*}
	where \(\beta=\frac 12\lambda_{\rm min}(G_\star)>0\).
	Therefore, there exists \(k_0\in\N\) such that
	\(
		\varphi_\gamma(\bar x^k+d^k)
	{}\leq{}
		\varphi_\gamma(\bar x^k)
	\)
	for all \(k\geq k_0\); in particular, for all such \(k\)
	\[
		\varphi_\gamma(\bar x^k+d^k)
	{}\leq{}
		\varphi_\gamma(\bar x^k)
	{}\leq{}
		\varphi_\gamma(x^k)
		{}-{}
		\gamma\tfrac{1-\gamma L_f}{2}\|r^k\|^2
	{}\leq{}
		\bar\Phi_k
		{}-{}
		\sigma\|r^k\|^2,
	\]
	where the second inequality follows from \cref{prop:FBEgeqGlobal}, and the last one from \eqref{eq:barPhi} and the fact that \(\sigma<\gamma\frac{1-\gamma L_f}{2}\).
	Therefore, for \(k\geq k_0\) the linesearch condition \eqref{eq:LS} holds with \(\tau_k=1\), and unitary stepsize is always accepted.
	In particular, the limit \eqref{eq:lim2} reads
	\(
		\lim_{k\to\infty}{
				\|x^{k+1}-x^\star\|/
				\|\bar x^k-x^\star\|
		}
	{}={}
		0,
	\)
	and from the inequality
	\begin{align*}
		\|\bar x^k-x^\star\|
	{}={} &
		\|\bar x^k-x^k+x^k-x^\star\|
	{}\leq{}
		\gamma\|\Res(x^k)\|
		{}+{}
		\|x^k-x^\star\|
	\\
	{}={} &
		\gamma\|\Res(x^k)-\Res(x^\star)\|
		{}+{}
		\|x^k-x^\star\|
	{}\leq{}
		(\gamma L_R+1)
		\|x^k-x^\star\|
	\end{align*}
	superlinear convergence of \(\seq{x^k}\) follows.
	Since
	\(
		\|r^k\|
	{}\leq{}
		L_R\|x^k-x^\star\|
	\),
	then also \(\seq{r^k}\) converges superlinearly, and in turn, since
	\(
		\|\bar x^k-x^\star\|
	{}\leq{}
		\gamma\|r^k\|
		{}+{}
		\|x^k-x^\star\|
	\),
	also \(\seq{\bar x^k}\) does.
	\arxiv{%

		The same proof applies replacing \(\Res\) with \(\tilde \Res\), by simply observing that \(\tilde \Res\) is strictly differentiable at \(x^\star\) and its Jacobian coincides with \(\nabla^2\varphi_\gamma(x^\star)\).
		In fact, close to \(x^\star\) we may write
		\(
			\tilde \Res(x)
		{}={}
			\Res(x)+\nabla f(\T(x))-\nabla f(x)
		\),
		therefore
		\begin{align}
		\nonumber
		J\tilde \Res(x^\star)
		{}={} &
			J\Res(x^\star)
			{}+{}
			\nabla^2f(\T(x^\star))J\T(x^\star)
			{}-{}
			\nabla^2f(x^\star)
		\\
		\label{eq:JtildeR}
		{}={} &
			Q_\gamma(x^\star)J\Res(x^\star)
		{}={}
			\nabla^2\varphi_\gamma(x^\star),
		\end{align}
		where the second equality comes from the fact that \(x^\star=\T(x^\star)\) and \(J\T=I-\gamma J\Res\).
	}%
\end{proof}
\end{thm}

We conclude the section showing that employing 
\switch{%
	Broyden directions \eqref{subeq:Broyden}
}{
	Broyden \eqref{subeq:Broyden} or BFGS \eqref{subeq:tildeBFGS} directions
}%
in \refZ[] enables superlinear convergence rates, provided that \(\Res\) is Lipschitz continuously \emph{semidifferentiable} at the limit point (see \cite{ip1992local}).
\begin{thm}[Superlinear convergence with Broyden \arxiv{and BFGS }directions]\label{thm:BroydenBFGS}%
Suppose that \Cref{ass:fg2} is (strictly) satisfied at a strong local minimum \(x^\star\) of \(\varphi\) at which \(\Res\) is Lipschitz-continuously semidifferentiable.
Consider the iterates generated by \refZ[] with directions \(d^k\) selected with
\switch{
	Broyden method \eqref{subeq:Broyden},
}{
	either Broyden method \eqref{subeq:Broyden} or BFGS scheme \eqref{subeq:tildeBFGS},
}%
and suppose that \(x^k\to x^\star\).

Then, the Dennis-Moré condition \eqref{eq:DM} is satisfied, and in particular all the claims of \Cref{thm:DM} hold.
\begin{proof}
	\switch{
		From \cref{thm:2ndOrder} we know that \(\Res\) is strictly differentiable at
		\newarxiv{the critical point }%
		\(x^\star\) and Lipschitz-continuously semidifferentiable there.
		Denoting \(G_\star=J\Res(x_\star)\),
	}{
		From the assumptions, \eqref{eq:JtildeR} and \cref{thm:2ndOrder}, we know that \(\Res\) and \(\tilde \Res\) are strictly differentiable at \(x^\star\), and Lipschitz-continuously semidifferentiable there.
		Suppose first that the directions are selected according to Broyden's method, and let \(G_\star=J\Res(x_\star)\).
		We have
	}%
	\[
		\frac{
			\|y^k-G_\star s^k\|
		}{
			\|s^k\|
		}
	{}={}
		\frac{
			\|\Res(x^{k+1})-\Res(\bar x^k)-G_\star(x^{k+1}-\bar x^k)\|
		}{
			\|x^{k+1}-\bar x^k\|
		}
	\]
	and since \(x^k,\bar x^k\to x^\star\), due to \cite[Lem. 2.2]{ip1992local} there exists \(L>0\) such that
	\Newswitch{\(}{\[}
		\frac{
			\|y^k-G_\star s^k\|
		}{
			\|s^k\|
		}
	{}\leq{}
		L
		\max{
			\set{
				\|x^{k+1}-x^\star\|
				{},{}
				\|\bar x^k-x^\star\|
			}
		}
	\Newswitch{\) for \(k\) large enough. }{\quad\text{for \(k\) large enough.}\]}%
	In particular, due to \cref{thm:Linear,lem:Lojasiewicz},
	\(
		\frac{
			\|y_k-G_\star s_k\|
		}{
			\|s_k\|
		}
	\)
	is summable.
	\arxiv{%
		Replacing \(\Res\) with \(\tilde \Res\), the same proof applies to the case in which BFGS directions are selected.
		We now differentiate the two cases.
		\begin{proofitemize}
		\item\emph{Broyden's directions.}
	}%
		Let \(E_k=B_k-G_\star\) and let \(\|{}\cdot{}\|_F\) denote the Frobenius norm.
		With a simple modification of the proofs of \cite[Thm. 4.1]{ip1992local} and \cite[Lem. 4.4]{artacho2014local} that takes into account the scalar \(\vartheta_k\in[\bar\vartheta,2-\bar\vartheta]\) we obtain
		\[
			\left\|E_{k+1}\right\|_F
		{}\leq{}
			\left\|
				E_k
				\left(
					\newarxiv{\!}
					\id-\vartheta_k\tfrac{s_k\trans{(s_k)}}{\|s_k\|^2}
				\right)
			\right\|_F
			\newarxiv{\!}
			{}+{}
			\vartheta_k
			\frac{\|y_k-G_\star s_k\|}{\|s_k\|}
		{}\leq{}
			\left\|E_k\right\|_F
			{}-{}
			\frac{\bar\vartheta(2-\bar\vartheta)}{2\|E_k\|_F}
			\frac{\|E_ks_k\|^2}{\|s_k\|^2}.
		\]
		The last term on the right-hand side, be it \(\sigma_k\), is summable and therefore the sequence \(\seq{E_k}\) is bounded.
		Therefore,
		\Newswitch{\(}{\[}%
			\|E_{k+1}\|_F
			{}-{}
			\|E_k\|_F
		{}\leq{}
			\sigma_k
			{}-{}
			\tfrac{\bar\vartheta(2-\bar\vartheta)}{2\bar E}
			\left(
				\tfrac{\|(B_k-G_\star)s_k\|}{\|s_k\|}
			\right)^{\!2}
		\Newswitch{\)}{\]}%
		where \(\bar E\coloneqq\sup\seq{\|E_k\|_F}\).
		Telescoping the inequality, summability of \(\sigma_k\) ensures that of
		\(
			\frac{\|(B_k-G_\star)s_k\|^2}{\|s_k\|^2}
		\)
		proving in particular the claimed Dennis-Moré condition \eqref{eq:DM}.
	\arxiv{%
		\item\emph{BFGS directions.}
			Since \(J\tilde \Res(x^\star)=\nabla^2\varphi_\gamma(x^\star)\) is symmetric and positive definite, and since \(x^k\to x^\star\), eventually \(\innprod{y^k}{s^k}>0\).
			We may now invoke \cite[Thm. 3.2]{byrd1989tool} to infer that the Dennis-Moré condition \eqref{eq:DM} holds.
		\end{proofitemize}
	}%
\end{proof}
\end{thm}

\section{Simulations}
	\label{sec:Simulations}
	We now present numerical results with the proposed method.
In \refZ[] we set $\beta = \nicefrac{1}{2}$, and for the nonmonotone linesearch we used the sequence $p_k = (\eta Q_k + 1)^{-1}$ where $Q_0 = 1$, $Q_{k+1} = \eta Q_k + 1$, $\eta = 0.85$: in this way $\seq{p_k}$ is computed as in \cite{zhang2004nonmonotone,li2015accelerated}.

We performed experiments with different choices of \(d^k\) in step \ref{state:zerofpr:d}.
In particular,
\begin{itemize}[
	leftmargin=*,
	labelindent=.5\parindent,
]
  \item \zerofpr(Broyden): \(d^k = -H_k\bar r^k\), and \(H_k\) obtained by the Broyden method \eqref{subeq:Broyden} with \(\bar\vartheta=10^{-4}\);
  \item \zerofpr(BFGS): \(d^k = -H_k\bar r^k\), where \(H_k\) is computed using BFGS updates \eqref{eq:BFGS};
  \item \zerofpr(L-BFGS): \(d^k\) is computed using L-BFGS \cite[Alg. 7.4]{nocedal2006numerical} with memory \(10\).
\end{itemize}
We only show the results with full quasi-Newton updates (Broyden, BFGS) for one of the examples: for the other experiments we focus on L-BFGS, which is better suited for large-scale problems.
Although $J\Res$ is nonsymmetric at the critical points in general, we observed that the symmetric updates of BFGS and L-BFGS perform very well in practice and outperform the Broyden method.

We compared \refZ[] with the forward-backward splitting algorithm (denoted FBS), that is \eqref{eq:FB}, the inertial FBS (denoted IFBS) proposed in \cite[Eq. (7)]{bot2016inertial} (with parameter \(\beta = 0.2\)), and the nonmonotone accelerated FBS (denoted AFBS) proposed in \cite[Alg. 2]{li2015accelerated} for fully nonconvex problems. 
All experiments were performed in MATLAB. The implementation of the methods used in the
tests are available online.\footnote{\url{http://github.com/kul-forbes/ForBES}}

	\subsection{Nonconvex sparse approximation}
		\label{sec:NoncvxSparse}
\begin{figure}
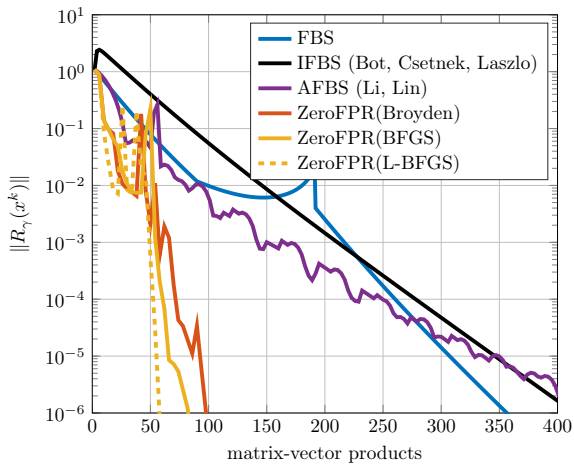

	\begin{minipage}[t]{0.6\linewidth}
		\vspace*{0pt}
		{{%
		\pgfkeys{/pgf/images/include external/.code={\includegraphics[width=\linewidth]{#width=\linewidth}}}%
		\tikzsetnextfilename{NoncvxSparseApprox/mu_3e-2_ops_vs_fpr}%
		\input{./TeX/Tikz/NoncvxSparseApprox/mu_3e-2_ops_vs_fpr.tex}%
	}}
	\end{minipage}
	\begin{minipage}[t]{0.39\linewidth}
		\vspace*{0pt}
		\caption{%
			\emph{Nonconvex sparse approximation.}
			Convergence of fixed-point residual in FBS and \refZ[], for different choices of the search directions, in the case \(n=1500\), \(\lambda = 0.03\).
			Since \refZ[] employs a linesearch, and therefore the complexity of each iteration is unknown \emph{a priori}, we recorded the number of matrix-vector products by \(A\) and \(\trans{A}\) performed during the iterations, and display it on the horizontal axis.%
		}%
		\label{fig:NoncvxSparseApprox1}%
	\end{minipage}
\end{figure}

		Here we consider the problem of finding a sparse solution to a least-squares problem. As discussed in \cite{xu2012regularization}, this is achieved by solving the following nonconvex problem:
\begin{equation}\label{eq:NoncvxSparse}
	\minimize\ \tfrac{1}{2}\|Ax-b\|_2^2 + \lambda\|x\|_{1/2}^{1/2},
\end{equation}
where \(\lambda>0\) is a regularization parameter, and \(\|x\|_{1/2} = \left(\sum_{i=1}^n |x_i|^{1/2}\right)^2\) is the \(\ell_{1/2}\) quasi-norm, a nonconvex regularizer whose role is that of inducing sparsity in the solution of \eqref{eq:NoncvxSparse}.
Function \(\|x\|_{1/2}^{1/2}\) is separable, and its proximal mapping can be computed in closed form as follows, see \cite[Thm. 1]{xu2012regularization}: for \(i=1,\ldots,n\)
\[
	\Bigl[
		\prox_{\gamma\|\cdot\|_{1/2}^{1/2}}(x)
	\Bigr]_i
{}={}
	\tfrac{2x_i}{3}\left(1+\cos\left(\tfrac{2\pi}{3}-\tfrac{2p_\gamma(x_i)}{3}\right)\right),
\]
where
\(
	p_\gamma(x_i)
{}={}
	\arccos{}
	\bigl(
		\nicefrac{\gamma}{8}
		(\nicefrac{|x_i|}{3})^{-\nicefrac 32}
	\bigr)
\).
We performed experiments using the setting of \cite[Sec. 8.2]{daubechies2010iteratively}: matrix \(A\in\R^{m\times n}\) has \(m=n/5\) rows and was generated with random Gaussian entries, with zero mean and variance \(1/m\).
Vector \(b\) was generated as \(b = Ax_{\mathrm{orig}} + v\) where \(x_{\mathrm{orig}}\in\R^n\) was randomly generated with \(k=5\) nonzero normally distributed entries, and \(v\) is a noise vector with zero mean and variance \(1/m\).
Then we solved problem \eqref{eq:NoncvxSparse} using \(x^0 = 0\) as starting iterate for all algorithms.
We computed the average and worst-case performance of the algorithms in a variety of scenarios, generating \(100\) random problems for each combination \((n, \lambda)\).
The results are illustrated in \Cref{tbl:NoncvxSparseApprox1}: \refZ[] finds local minima significantly faster than FBS, IFBS and AFBS.
\Cref{fig:NoncvxSparseApprox1} shows the behavior of the considered algorithms in one of the generated problems, where fast asymptotic convergence of \refZ[] is apparent.

\begin{table}[tbp]
	\centering
	\footnotesize
	\setlength{\tabcolsep}{10pt}%
	\newarxiv{\setlength{\tabcolsep}{8pt}}%
	\begin{tabular}{|r|r|r|r|r|r|}
	\hline
	$n$ & $\lambda$ & \multicolumn{1}{c|}{FBS} & \multicolumn{1}{c|}{IFBS} & \multicolumn{1}{c|}{AFBS} & \multicolumn{1}{c|}{ZeroFPR(L-BFGS)} \\
		&           & avg/max (s) & avg/max (s) & avg/max (s) & avg/max (s)\\
	\hline
	500 & 0.10 & 0.141/0.405 & 0.159/0.449 & 0.135/0.221 & 0.037/0.088 \\
		& 0.03 & 0.498/2.548 & 0.688/3.962 & 0.274/0.430 & 0.084/0.126 \\
		& 0.01 & 1.305/5.445 & 1.721/4.942 & 0.570/1.157 & 0.152/0.560 \\
	\hline
	1000 & 0.10 & 0.176/0.287 & 0.231/0.659 & 0.228/0.483 & 0.021/0.077 \\
		& 0.03 & 0.576/2.756 & 0.645/4.165 & 0.382/0.841 & 0.091/0.275 \\
		& 0.01 & 1.864/9.740 & 2.391/8.311 & 0.795/1.446 & 0.222/0.438 \\
	\hline
	2000 & 0.10 & 0.291/0.599 & 0.392/0.719 & 0.393/0.640 & 0.025/0.055 \\
		& 0.03 & 0.553/1.841 & 0.602/3.270 & 0.464/0.702 & 0.088/0.198 \\
		& 0.01 & 2.108/10.934 & 2.439/8.010 & 0.979/1.411 & 0.271/0.464 \\
	\hline
	\end{tabular}
	\caption{%
		Nonconvex sparse approximation.
		Performance of FBS, IFBS, AFBS and \refZ[] on problems with different values of \(n\) and \(\lambda\).
		The table shows average and maximum CPU time required to reach \(\|R_\gamma(x^k)\|\leq 10^{-6}\) in \(100\) random experiments.
		Each algorithm was run on the same set of randomly generated problems, with \(x^0=0\).%
	}%
	\label{tbl:NoncvxSparseApprox1}%
\end{table}


	\subsection{Dictionary learning}
		\label{sec:DictionaryLearning}
		Given a collection of \(m\) signals of dimension \(n\), collected as columns in a matrix \(Y\in\R^{n\times m}\), we seek for a sparse representation of each of them as combination of a set of \(k\) vectors \(\set{d_1,\ldots,d_k}\), called \emph{dictionary atoms}. To do so, we solve the following problem
\begin{equation}\label{eq:DictionaryLearning}
	\minimize_{D, C}\ \tfrac12\|Y - DC\|_F^2
\quad\stt~~
	\begin{array}[t]{l@{~~}rl}
		\|d_i\|_2 = 1        & i={} & 1,\ldots,k, \\
		\|c_j\|_0\leq N      & j={} & 1,\ldots,m, \\
		\|c_j\|_\infty\leq T & j={} & 1,\ldots,m,
	\end{array}
\end{equation}
where \(D=(d_1,\ldots,d_k)\in\R^{n\times k}\), \(C=(c_1,\ldots,c_n)\in\R^{k\times m}\), \(\|\cdot\|_0\) is the \(\ell_0\) pseudo-norm, \ie the number of nonzero coefficents, while \(N\in\N\) and \(T>0\) are parameters.
This is similar to the problem considered in \cite{aharon2006ksvd}.
Here we bound the set of feasible points by means of the \(\ell_\infty\)-norm constraint: this has the effect of making the domain of problem \eqref{eq:DictionaryLearning} compact and \(\nabla f\), as a consequence, Lipschitz continuous over the problem domain (cf. \cref{rem:C1+}).
Furthermore, here we explicitly constrain the norm of the dictionary atoms: in fact, the objective value of \eqref{eq:DictionaryLearning} is unchanged if an atom \(d_i\) is scaled by a factor, and the corresponding row of \(C\) is scaled by the inverse factor.

Problem \eqref{eq:DictionaryLearning} takes the form \eqref{eq:Problem} by letting \(f(D, C) = \tfrac12\|Y-DC\|^2_F\) and \(g(D,C) = \delta_S(D,C)\), where
\(
	S
{}={}
	\overbrace{S_D \times \ldots \times S_D}^{\text{\(k\) times}}
	{}\times{}
	\overbrace{S_C \times \ldots \times S_C}^{\text{\(m\) times}},
\)
with
\[
  S_D {}={} \set{d\in\R^n}[\|d\|_2 = 1], \quad
  S_C {}={} \set{\smash{c\in\R^k}}[\|c\|_0 \leq N, \|c\|_\infty \leq T],
\]
\ie set \(S\) is the product of Euclidean spheres and box-constrained \(\ell_0\) level sets.
Both \(f\) and \(g\) are nonconvex in this case.
Projection onto \(S_D\) is simply a matter of scaling, while that onto \(S_C\) simply amounts to projecting the \(N\) largest coefficients (in absolute value) onto the box \([-T, T]\) and setting to zero the remaining ones.

\begin{figure}[tbp]
	\begin{minipage}[t]{0.6\linewidth}
		\vspace*{0pt}
		{{%
		\pgfkeys{/pgf/images/include external/.code={\includegraphics[width=\linewidth]{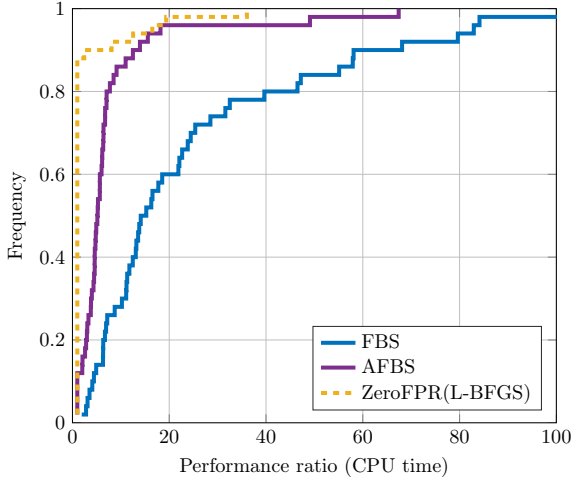}}}%
		\tikzsetnextfilename{DictLearning/matfact_perf_cputime}%
%
%
\definecolor{mycolor1}{rgb}{1.00000,0.00000,1.00000}%
\begin{tikzpicture}

\begin{axis}[%
grid,
at={(1.011in,0.642in)},
scale only axis,
xmin=0,
xmax=100.0,
ymin=0,
ymax=1,
axis background/.style={fill=white},
title style={font=\bfseries},
xlabel={Performance ratio (CPU time)},
ylabel={Frequency},
legend style={legend cell align=left,align=left,draw=white!15!black},
legend pos=south east
]
\addplot [color=color_FBS,solid,line width=2.0pt,
]
  table[row sep=crcr]{%
  1.91567944697685	0.02\\
  2.81250108712113	0.02\\
  2.81250108712113	0.04\\
  3.11921005735204	0.04\\
  3.11921005735204	0.06\\
  3.52472217719163	0.06\\
  3.52472217719163	0.08\\
  4.10391242256737	0.08\\
  4.10391242256737	0.1\\
  4.43850262304769	0.1\\
  4.43850262304769	0.12\\
  4.90709908938579	0.12\\
  4.90709908938579	0.14\\
  6.3112021748218	0.14\\
  6.3112021748218	0.16\\
  6.3344588635766	0.16\\
  6.3344588635766	0.18\\
  6.3633819065776	0.18\\
  6.3633819065776	0.2\\
  6.7148769982079	0.2\\
  6.7148769982079	0.22\\
  6.95110950825902	0.22\\
  6.95110950825902	0.24\\
  7.18596251263235	0.24\\
  7.18596251263235	0.26\\
  8.7291506806401	0.26\\
  8.7291506806401	0.28\\
  10.2210476365986	0.28\\
  10.2210476365986	0.3\\
  11.1015114860911	0.3\\
  11.1015114860911	0.32\\
  11.1700583895906	0.32\\
  11.1700583895906	0.34\\
  11.3132385152794	0.34\\
  11.3132385152794	0.36\\
  11.7659202184527	0.36\\
  11.7659202184527	0.38\\
  12.4685659593654	0.38\\
  12.4685659593654	0.4\\
  13.0884793969968	0.4\\
  13.0884793969968	0.42\\
  13.2101141061727	0.42\\
  13.2101141061727	0.44\\
  13.5709713600855	0.44\\
  13.5709713600855	0.46\\
  13.8005376232451	0.46\\
  13.8005376232451	0.48\\
  14.1124223454658	0.48\\
  14.1124223454658	0.5\\
  15.2582021259549	0.5\\
  15.2582021259549	0.52\\
  16.2912493197452	0.52\\
  16.2912493197452	0.54\\
  16.5361337268762	0.54\\
  16.5361337268762	0.56\\
  17.7439542789797	0.56\\
  17.7439542789797	0.58\\
  18.5458412102046	0.58\\
  18.5458412102046	0.6\\
  21.8948569769454	0.6\\
  21.8948569769454	0.62\\
  22.0941524370904	0.62\\
  22.0941524370904	0.64\\
  22.6689359719396	0.64\\
  22.6689359719396	0.66\\
  23.8290293649467	0.66\\
  23.8290293649467	0.68\\
  24.4437756639771	0.68\\
  24.4437756639771	0.7\\
  25.3760022252135	0.7\\
  25.3760022252135	0.72\\
  28.5104420101124	0.72\\
  28.5104420101124	0.74\\
  31.5993058793729	0.74\\
  31.5993058793729	0.76\\
  32.5179910958846	0.76\\
  32.5179910958846	0.78\\
  39.6516530048859	0.78\\
  39.6516530048859	0.8\\
  46.5095500984769	0.8\\
  46.5095500984769	0.82\\
  47.192173711758	0.82\\
  47.192173711758	0.84\\
  55.1111407701651	0.84\\
  55.1111407701651	0.86\\
  57.9257999395547	0.86\\
  57.9257999395547	0.88\\
  58.0991806935582	0.88\\
  58.0991806935582	0.9\\
  68.1197164466936	0.9\\
  68.1197164466936	0.92\\
  79.6399250327896	0.92\\
  79.6399250327896	0.94\\
  82.9385636451353	0.94\\
  82.9385636451353	0.96\\
  84.1178179446665	0.96\\
  84.1178179446665	0.98\\
  106.711702642423	0.98\\
  106.711702642423	1\\
};
\addlegendentry{FBS};

\addplot [color=color_AFBS,solid,line width=2.0pt,
]
  table[row sep=crcr]{%
  1	0.02\\
  1	0.02\\
  1	0.04\\
  1	0.04\\
  1	0.06\\
  1	0.06\\
  1	0.08\\
  1	0.08\\
  1	0.1\\
  1	0.1\\
  1	0.12\\
  2.0290313565877	0.12\\
  2.0290313565877	0.14\\
  2.07565043004577	0.14\\
  2.07565043004577	0.16\\
  2.60548028667929	0.16\\
  2.60548028667929	0.18\\
  2.82370037639502	0.18\\
  2.82370037639502	0.2\\
  2.98565055943808	0.2\\
  2.98565055943808	0.22\\
  3.06729686657146	0.22\\
  3.06729686657146	0.24\\
  3.30596846081398	0.24\\
  3.30596846081398	0.26\\
  3.79416396311413	0.26\\
  3.79416396311413	0.28\\
  3.8336176177495	0.28\\
  3.8336176177495	0.3\\
  4.00962342048357	0.3\\
  4.00962342048357	0.32\\
  4.34980151031009	0.32\\
  4.34980151031009	0.34\\
  4.47396146484014	0.34\\
  4.47396146484014	0.36\\
  4.59667269830089	0.36\\
  4.59667269830089	0.38\\
  4.59844718296093	0.38\\
  4.59844718296093	0.4\\
  4.61020964713299	0.4\\
  4.61020964713299	0.42\\
  4.75212511274206	0.42\\
  4.75212511274206	0.44\\
  4.81854750712483	0.44\\
  4.81854750712483	0.46\\
  4.92129062408832	0.46\\
  4.92129062408832	0.48\\
  5.1417114086978	0.48\\
  5.1417114086978	0.5\\
  5.21990382269349	0.5\\
  5.21990382269349	0.52\\
  5.27339657137621	0.52\\
  5.27339657137621	0.54\\
  5.65522870532505	0.54\\
  5.65522870532505	0.56\\
  5.66394883986475	0.56\\
  5.66394883986475	0.58\\
  5.66422632260331	0.58\\
  5.66422632260331	0.6\\
  6.01475129071327	0.6\\
  6.01475129071327	0.62\\
  6.12902512279619	0.62\\
  6.12902512279619	0.64\\
  6.16231466613065	0.64\\
  6.16231466613065	0.66\\
  6.35694952110141	0.66\\
  6.35694952110141	0.68\\
  6.40727642147161	0.68\\
  6.40727642147161	0.7\\
  6.41619335347321	0.7\\
  6.41619335347321	0.72\\
  6.74659342014477	0.72\\
  6.74659342014477	0.74\\
  6.75437347769891	0.74\\
  6.75437347769891	0.76\\
  7.0377211778802	0.76\\
  7.0377211778802	0.78\\
  7.06109417772636	0.78\\
  7.06109417772636	0.8\\
  7.73934835247104	0.8\\
  7.73934835247104	0.82\\
  8.50094361959408	0.82\\
  8.50094361959408	0.84\\
  9.11719453275933	0.84\\
  9.11719453275933	0.86\\
  10.9894660029494	0.86\\
  10.9894660029494	0.88\\
  12.4926272794372	0.88\\
  12.4926272794372	0.9\\
  13.9467122357913	0.9\\
  13.9467122357913	0.92\\
  15.5770014732925	0.92\\
  15.5770014732925	0.94\\
  18.1957801946851	0.94\\
  18.1957801946851	0.96\\
  49.0801470466952	0.96\\
  49.0801470466952	0.98\\
  67.4166255800897	0.98\\
  67.4166255800897	1\\
};
\addlegendentry{AFBS};

\addplot [color=color_BFGS,\linestylelimited,line width=2.0pt,
]
  table[row sep=crcr]{%
  1	0.02\\
  1	0.02\\
  1	0.04\\
  1	0.04\\
  1	0.06\\
  1	0.06\\
  1	0.08\\
  1	0.08\\
  1	0.1\\
  1	0.1\\
  1	0.12\\
  1	0.12\\
  1	0.14\\
  1	0.14\\
  1	0.16\\
  1	0.16\\
  1	0.18\\
  1	0.18\\
  1	0.2\\
  1	0.2\\
  1	0.22\\
  1	0.22\\
  1	0.24\\
  1	0.24\\
  1	0.26\\
  1	0.26\\
  1	0.28\\
  1	0.28\\
  1	0.3\\
  1	0.3\\
  1	0.32\\
  1	0.32\\
  1	0.34\\
  1	0.34\\
  1	0.36\\
  1	0.36\\
  1	0.38\\
  1	0.38\\
  1	0.4\\
  1	0.4\\
  1	0.42\\
  1	0.42\\
  1	0.44\\
  1	0.44\\
  1	0.46\\
  1	0.46\\
  1	0.48\\
  1	0.48\\
  1	0.5\\
  1	0.5\\
  1	0.52\\
  1	0.52\\
  1	0.54\\
  1	0.54\\
  1	0.56\\
  1	0.56\\
  1	0.58\\
  1	0.58\\
  1	0.6\\
  1	0.6\\
  1	0.62\\
  1	0.62\\
  1	0.64\\
  1	0.64\\
  1	0.66\\
  1	0.66\\
  1	0.68\\
  1	0.68\\
  1	0.7\\
  1	0.7\\
  1	0.72\\
  1	0.72\\
  1	0.74\\
  1	0.74\\
  1	0.76\\
  1	0.76\\
  1	0.78\\
  1	0.78\\
  1	0.8\\
  1	0.8\\
  1	0.82\\
  1	0.82\\
  1	0.84\\
  1	0.84\\
  1	0.86\\
  1	0.86\\
  1	0.88\\
  2.42684970129721	0.88\\
  2.42684970129721	0.9\\
  8.07685727699042	0.9\\
  8.07685727699042	0.92\\
  12.5435382915583	0.92\\
  12.5435382915583	0.94\\
  16.5012947824111	0.94\\
  16.5012947824111	0.96\\
  19.3777632534125	0.96\\
  19.3777632534125	0.98\\
  36.0439251676772	0.98\\
  36.0439251676772	1\\
};
\addlegendentry{ZeroFPR(L-BFGS)};

\end{axis}
\end{tikzpicture}%
%
	}}
	\end{minipage}
	\hfill
	\begin{minipage}[t]{0.39\linewidth}
		\vspace*{0pt}
		\centering
		\caption{%
			\emph{Dictionary learning.}
			Performance profiles of FBS, AFBS and \refZ[](L-BFGS) when applied to \(50\) randomly generated problems with \(n=20\), \(m=500\), \(k=50\), \(T=10^6\) and \(N=3\).
			The algorithms are executed until tolerance \(\|R_\gamma(x^k)\| \leq 10^{-4}\) is reached.
			In the great majority of cases, \refZ[](L-BFGS) reaches a critical point significantly faster than FBS.
		}%
		\label{fig:DictLearning1}%
	\end{minipage}
\end{figure}

We have tested the proposed algorithm on a sequence of problems generated according to \cite[§V.A]{aharon2006ksvd}.
We set \(n=20\), \(m=500\), \(k=50\), therefore in this case the problem has \(26000\) variables.
A generating dictionary \(D_{\mathrm{gen}}\in\R^{20\times 50}\) was selected randomly generated with normal entries, and each column was normalized to one.
Then a random matrix \(C_{\mathrm{gen}}\in\R^{50\times 500}\) was selected with \(3\) normally distributed nonzero coefficient per column.
Then we set \(Y=C_{\mathrm{gen}}D_{\mathrm{gen}} + V\), where \(V\in\R^{n\times m}\) is normally distributed with variance \(10^{-2}\).
We generated \(50\) random problems according to this procedure, and applied the algorithms to problem \eqref{eq:DictionaryLearning} with \(N=3\) and \(T=10^6\).
In this case
IFBS is not applicable since the Lipschitz constant of \(\nabla f\) over the problem domain is unknown, and an appropriate stepsize-selection rule is not provided for the algorithm in \cite{bot2016inertial}.
We compared FBS, AFBS and \refZ[], using the backtracking procedure discussed in \Cref{rem:Adaptive} to adaptively adjust the stepsize \(\gamma\).
We used \((D^0, C^0) = (0, 0)\) as initial iterate, while the algorithms were stopped as soon as \(\|R_\gamma(x^k)\| \leq 10^{-4}\).
The results are shown in \Cref{fig:DictLearning1}:
in most of the cases, \refZ[](L-BFGS) exhibited a speedup of a factor \(5\)-to-\(100\) with respect to FBS, and \(3\)-to-\(60\) with respect to AFBS, at reaching a critical point.

	\subsection{Matrix decomposition}
		\label{sec:MatrixDecomposition}
		We consider the problem of approximating a given matrix $A\in \R^{m\times n}$ as the sum of a low-rank and a sparse component, by solving
\begin{align}\label{eq:MatDecomp}
	\minimize_{X_L, X_S\in\R^{m\times n}} \tfrac{1}{2}\|A - X_L - X_S\|_F^2 + \lambda \|X_S\|_0
\quad
	\stt\ \rank(X_L) \leq r.
\end{align}
This problem has application, for example, in the analysis of video imagery,
specifically the separation of the background (fixed over time) scenery
from the foreground (moving) objects in a series of video frames.
In this case, matrix $A$ contains $n$ video frames (columns), each
consisting of $m$ pixels, and $X_L$, $X_S$ will respectively contain the background scenery
and foreground objects identified in each frame.
Therefore here $f(X_L, X_S) = \tfrac{1}{2}\|A - X_L - X_S\|_F^2$ and $L_f = 2$,
while $g(X_L, X_S) = \indicator_{\rank \leq r}(X_L) + \lambda\|X_S\|_0$.
The proximal mapping of $g$ is given by
$$ \prox_{\gamma g}(X_L, X_S) = (\proj_{\rank \leq r}(X_L), \prox_{\gamma \lambda \|\cdot\|_0}(X_S)). $$
Here, $\prox_{\gamma \lambda \|\cdot\|_0}$ is the hard-thresholding operation, defined componentwise as
\[
\left[\prox_{\gamma \lambda \|\cdot\|_0} (X_S)\right]_{ij} =
 \begin{cases}
 	0 &\quad \mbox{if\ } \bigl(X_S\bigr)_{ij} \leq \sqrt{2\gamma\lambda} \\
 	\bigl(X_S\bigr)_{ij} &\quad\mbox{otherwise}
 \end{cases}
\quad
	i,j=1,\ldots,m.
\]
The set of matrices of rank at most $r$ is nonconvex and closed, and the projection onto it is given by
\(
	\proj_{\rank \leq r}(X) {}={} U_r \diag(\sigma_1, \ldots, \sigma_r) V_r^T
\),
where $\sigma_1\ldots \sigma_r$ are $r$ largest singular values of $X$, and $U_r, V_r$ are the matrices of left and right singular vectors, respectively.
Each computation of $\Pi_{\rank \leq r}$ requires a partial SVD which is, from the computational perspective, the most significantly expensive operation in this case.

\begin{figure}[tbp]
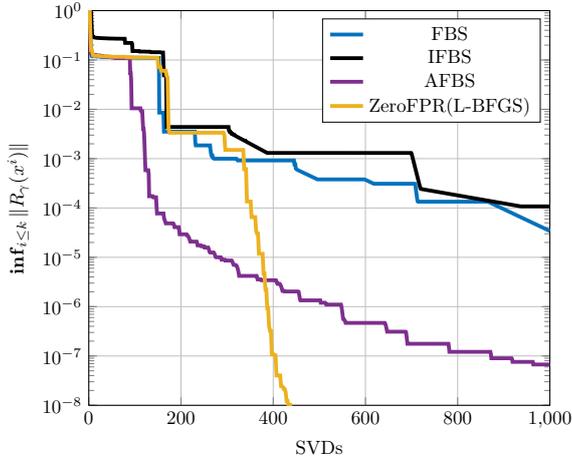

	\begin{minipage}[t]{0.6\linewidth}
		\vspace*{0pt}
		\centering
		{{%
		\pgfkeys{/pgf/images/include external/.code={\includegraphics[width=\linewidth]{#width=\linewidth}}}%
		\tikzsetnextfilename{MatDecomp/svds_vs_fpr_50frames_lam_3e-3}%
		\input{./TeX/Tikz/MatDecomp/svds_vs_fpr_50frames_lam_3e-3.tex}%
	}}
	\end{minipage}
	\hfill
	\begin{minipage}[t]{0.39\linewidth}
		\vspace*{0pt}%
		\caption{%
			\emph{Matrix decomposition.}
			Convergence of the fixed-point residual of \refZ[] using L-BFGS directions compared to FBS and IFBS on the \emph{ShoppingMall} dataset.\newline
		}%
		\label{fig:BackFore1}%
	\end{minipage}
\end{figure}

\begin{figure}[tbp]
	\centering
%
	\includegraphics[width=0.25\textwidth]{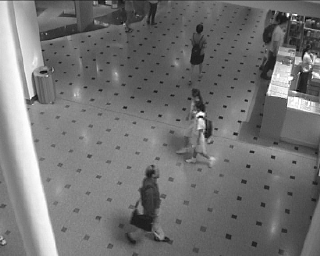}
	\hspace{10pt}
	\includegraphics[width=0.25\textwidth]{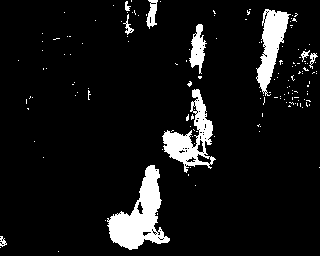}
	\hspace{10pt}
	\includegraphics[width=0.25\textwidth]{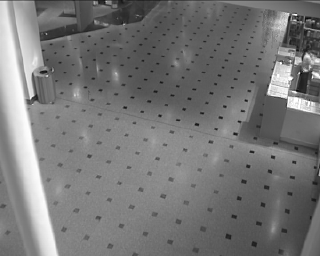}

	\vspace{5pt}
	\includegraphics[width=0.25\textwidth]{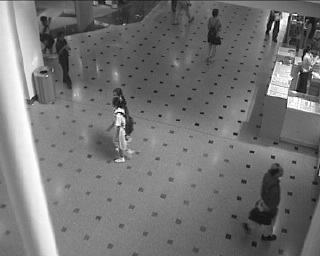}
	\hspace{10pt}
	\includegraphics[width=0.25\textwidth]{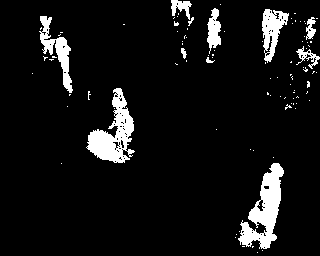}
	\hspace{10pt}
	\includegraphics[width=0.25\textwidth]{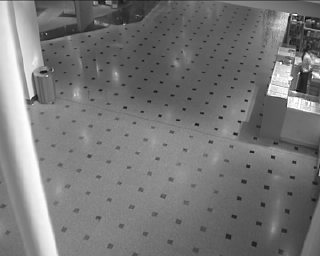}
	\caption{%
		Matrix decomposition: separation of the background scenery from the foreground moving objects in the \emph{ShoppingMall} dataset.%
	}%
	\label{fig:BackFore2}%
\end{figure}

We applied this technique to a sequence of $n = 50$ frames coming from the
\emph{ShoppingMall} dataset.\footnote{\url{http://perception.i2r.a-star.edu.sg/bk_model/bk_index.html}}
The footage consists of \(m = 320\times 256\) grayscale pixels frames, therefore the problem has \(8192000\) variables in total.
In problem \eqref{eq:MatDecomp} we used $r = 1$ and $\lambda = 3\cdot 10^{-3}$.
The results are shown in \Cref{fig:BackFore1,fig:BackFore2}. Also in this case, the fast asymptotic convergence of \refZ[](L-BFGS) is apparent.

\section{Conclusions}
	\label{sec:Conclusions}
The forward-backward envelope is a valuable tool for deriving efficient algorithms tackling nonsmooth and nonconvex problems of the form $\varphi = f+g$, as it can be used as a merit function to devise globally convergent linesearch methods solving the system of nonlinear equations defining the stationary points of $\varphi$.

\refZ[] implements this idea, and we proved that it globally converges to a stationary point under the assumption that $\varphi_\gamma$ has the Kurdyka-{\L}ojasiewicz property. Furthermore, if the linesearch directions satisfy the Dennis-Mor\'e condition (for example, if they are determined according to the Broyden method), the convergence rate at strong local minima is superlinear.

Numerical simulations with the proposed method on convex and nonconvex problems confirm our theoretical results.
Using Broyden method, BFGS (in the case of small-scale problems) and L-BFGS (for large-scale problems) to compute directions in \refZ[] greatly outperform FBS and its accelerated variant.
It is our belief that the surprising efficacy of \mbox{(L-)BFGS} is due to the fact that, under the appropriate assumptions, the Jacobian of \(\Res\) at strong local minima is similar to a symmetric and positive definite matrix. Future investigation may better explain the effectiveness of symmetric update formulas in this framework.

\arxivtrue

	\phantomsection
	\switch{
		\bibliographystyle{siamplain}
	}{
		\bibliographystyle{plain}
	}%
	\bibliography{TeX/Bibliography.bib}

\begin{thebibliography}{10}

\bibitem{aharon2006ksvd}
Michal Aharon, Michael Elad, and Alfred Bruckstein.
\newblock {K-SVD}: An algorithm for designing overcomplete dictionaries for
  sparse representation.
\newblock {\em IEEE Transactions on signal processing}, 54(11):4311--4322,
  2006.

\bibitem{artacho2014local}
F.~J. Arag{\'o}n~Artacho, A.~Belyakov, A.~L. Dontchev, and M.~L{\'o}pez.
\newblock Local convergence of quasi-{N}ewton methods under metric regularity.
\newblock {\em Computational Optimization and Applications}, 58(1):225--247,
  2014.

\bibitem{attouch2009convergence}
Hedy Attouch and J\'er\^ome Bolte.
\newblock On the convergence of the proximal algorithm for nonsmooth functions
  involving analytic features.
\newblock {\em Mathematical Programming}, 116(1-2):5--16, 2009.

\bibitem{attouch2010proximal}
H{\'e}dy Attouch, J{\'e}r{\^o}me Bolte, Patrick Redont, and Antoine Soubeyran.
\newblock Proximal alternating minimization and projection methods for
  nonconvex problems: An approach based on the {K}urdyka-{\l}o\-ja\-sie\-wicz
  inequality.
\newblock {\em Mathematics of Operations Research}, 35(2):438--457, 2010.

\bibitem{attouch2013convergence}
Hedy Attouch, J{\'e}r{\^o}me Bolte, and Benar~Fux Svaiter.
\newblock Convergence of descent methods for semi-algebraic and tame problems:
  proximal algorithms, forward--backward splitting, and regularized
  gauss--seidel methods.
\newblock {\em Mathematical Programming}, 137(1):91--129, 2013.

\bibitem{attouch2016rate}
Hedy Attouch and Juan Peypouquet.
\newblock The rate of convergence of {N}esterov's accelerated forward-backward
  method is actually faster than \(1/k^2\).
\newblock {\em SIAM Journal on Optimization}, 26(3):1824--1834, 2016.

\bibitem{beck2016minimization}
Amir Beck and Nadav Hallak.
\newblock On the minimization over sparse symmetric sets: Projections,
  optimality conditions, and algorithms.
\newblock {\em Math. Oper. Res.}, 41(1):196--223, 2016.

\bibitem{beck2009fast}
Amir Beck and Marc Teboulle.
\newblock {A Fast Iterative Shrinkage-Thresholding Algorithm for Linear Inverse
  Problems}.
\newblock {\em SIAM Journal on Imaging Sciences}, 2(1):183--202, 2009.

\bibitem{bertsekas1995nonlinear}
Dimitri~P. Bertsekas.
\newblock {\em Nonlinear Programming}.
\newblock Athena Scientific, 1995.

\bibitem{bochnak1998real}
J.~Bochnak, M.~Coste, and M-F. Roy.
\newblock {\em Real Algebraic Geometry}.
\newblock Springer, 1998.

\bibitem{bolte2007lojasiewicz}
J{\'{e}}r{\^{o}}me Bolte, Aris Daniilidis, and Adrian Lewis.
\newblock The {{\L}}o\-ja\-sie\-wicz inequality for nonsmooth subanalytic
  functions with applications to subgradient dynamical systems.
\newblock {\em SIAM Journal on Optimization}, 17(4):1205--1223, 2007.

\bibitem{bolte2014proximal}
J{\'e}r{\^o}me Bolte, Shoham Sabach, and Marc Teboulle.
\newblock Proximal alternating linearized minimization for nonconvex and
  nonsmooth problems.
\newblock {\em Mathematical Programming}, 146(1-2):459--494, 2014.

\bibitem{bolte2016majorization}
Jérôme Bolte and Edouard Pauwels.
\newblock Majorization-minimization procedures and convergence of {SQP} methods
  for semi-algebraic and tame programs.
\newblock {\em Mathematics of Operations Research}, 41(2):442--465, 2016.

\bibitem{bot2016inertial}
Radu~Ioan Bo{\c{t}}, Ern{\"o}~Robert Csetnek, and Szil{\'a}rd~Csaba
  L{\'a}szl{\'o}.
\newblock An inertial forward--backward algorithm for the minimization of the
  sum of two nonconvex functions.
\newblock {\em EURO Journal on Computational Optimization}, 4(1):3--25, 2016.

\bibitem{broyden1965class}
Charles~G. Broyden.
\newblock A class of methods for solving nonlinear simultaneous equations.
\newblock {\em Mathematics of Computation}, 19(92):577--593, 1965.

\bibitem{byrd1989tool}
Richard~H. Byrd and Jorge Nocedal.
\newblock A tool for the analysis of quasi-{N}ewton methods with application to
  unconstrained minimization.
\newblock {\em SIAM J. Numer. Anal.}, 26(3):727--739, June 1989.

\bibitem{daniilidis2006geometrical}
Aris Daniilidis, Warren Hare, and Jérôme Malick.
\newblock Geometrical interpretation of the predictor-corrector type algorithms
  in structured optimization problems.
\newblock {\em Optimization}, 55(5-6):481--503, 2006.

\bibitem{daubechies2010iteratively}
Ingrid Daubechies, Ronald DeVore, Massimo Fornasier, and C.~Sinan
  G{\"u}nt{\"u}rk.
\newblock Iteratively reweighted least squares minimization for sparse
  recovery.
\newblock {\em Communications on Pure and Applied Mathematics}, 63(1):1--38,
  2010.

\bibitem{dennis1974characterization}
John~E. Dennis and Jorge~J. Mor{\'e}.
\newblock A characterization of superlinear convergence and its application to
  quasi-{N}ewton methods.
\newblock {\em Mathematics of computation}, 28(126):549--560, 1974.

\bibitem{dontchev2012generalizations}
Asen Dontchev.
\newblock Generalizations of the dennis--moré theorem.
\newblock {\em SIAM Journal on Optimization}, 22(3):821--830, 2012.

\bibitem{frankel2015splitting}
Pierre Frankel, Guillaume Garrigos, and Juan Peypouquet.
\newblock Splitting methods with variable metric for kurdyka--{\l}ojasiewicz
  functions and general convergence rates.
\newblock {\em Journal of Optimization Theory and Applications},
  165(3):874--900, 2015.

\bibitem{ip1992local}
Chi-Ming Ip and Jerzy Kyparisis.
\newblock Local convergence of quasi-{N}ewton methods for {B}-differentiable
  equations.
\newblock {\em Mathematical Programming}, 56(1-3):71--89, 1992.

\bibitem{kaplan1998proximal}
A.~Kaplan and R.~Tichatschke.
\newblock Proximal point methods and nonconvex optimization.
\newblock {\em Journal of Global Optimization}, 13(4):389--406, 1998.

\bibitem{kurdyka1998gradients}
Krzysztof Kurdyka.
\newblock On gradients of functions definable in o-minimal structures.
\newblock {\em Annales de l'institut Fourier}, 48(3):769--783, 1998.

\bibitem{lewis2002active}
A.~S. Lewis.
\newblock Active sets, nonsmoothness, and sensitivity.
\newblock {\em SIAM Journal on Optimization}, 13(3):702--725, 2002.

\bibitem{li2015accelerated}
Huan Li and Zhouchen Lin.
\newblock Accelerated proximal gradient methods for nonconvex programming.
\newblock In {\em Advances in neural information processing systems}, pages
  379--387, 2015.

\bibitem{liu1989limited}
Dong~C. Liu and Jorge Nocedal.
\newblock On the limited memory {BFGS} method for large scale optimization.
\newblock {\em Mathematical Programming}, 45(1-3):503--528, 1989.

\bibitem{liu2017further}
Tianxiang Liu and Ting~Kei Pong.
\newblock Further properties of the forward--backward envelope with
  applications to difference-of-convex programming.
\newblock {\em Computational Optimization and Applications}, pages 1--32, 2017.

\bibitem{lojasiewicz1963propriete}
Stanislaw {\L}ojasiewicz.
\newblock Une propri{\'e}t{\'e} topologique des sous-ensembles analytiques
  r{\'e}els.
\newblock {\em Les {\'e}quations aux d{\'e}riv{\'e}es partielles}, pages
  87--89, 1963.

\bibitem{lojasiewicz1993geometrie}
Stanislaw {\L}ojasiewicz.
\newblock Sur la g{\'e}om{\'e}trie semi- et sous- analytique.
\newblock {\em Annales de l'institut Fourier}, 43(5):1575--1595, 1993.

\bibitem{martinet1970breve}
B.~Martinet.
\newblock Brève communication. {R}égularisation d'inéquations
  variationnelles par approximations successives.
\newblock {\em ESAIM: Modélisation Mathématique et Analyse Numérique},
  4(R3):154--158, 1970.

\bibitem{nesterov2013gradient}
Yuri Nesterov.
\newblock {Gradient Methods for Minimizing Composite Functions}.
\newblock {\em Mathematical Programming}, 140(1):125--161, 2013.

\bibitem{nocedal1980updating}
Jorge Nocedal.
\newblock Updating quasi-{N}ewton matrices with limited storage.
\newblock {\em Mathematics of computation}, 35(151):773--782, 1980.

\bibitem{nocedal2006numerical}
Jorge Nocedal and Stephen Wright.
\newblock {\em Numerical Optimization}.
\newblock Springer, New York, 2nd edition edition, August 2006.

\bibitem{ochs2014ipiano}
Peter Ochs, Yunjin Chen, Thomas Brox, and Thomas Pock.
\newblock i{P}iano: Inertial proximal algorithm for nonconvex optimization.
\newblock {\em SIAM Journal on Imaging Sciences}, 7(2):1388--1419, 2014.

\bibitem{patrinos2013proximal}
Panagiotis Patrinos and Alberto Bemporad.
\newblock Proximal {N}ewton methods for convex composite optimization.
\newblock In {\em IEEE Conference on Decision and Control}, pages 2358--2363,
  2013.

\bibitem{poliquin1996generalized}
RA~Poliquin and RT~Rockafellar.
\newblock Generalized {H}essian properties of regularized nonsmooth functions.
\newblock {\em SIAM Journal on Optimization}, 6(4):1121--1137, 1996.

\bibitem{poliquin1992nonsmooth}
Ren\'e~A. Poliquin and R.~Tyrrell Rockafellar.
\newblock Amenable functions in optimization.
\newblock {\em Nonsmooth optimization: methods and applications}, pages
  338--353, 1992.

\bibitem{poliquin1995second}
Ren{\'e}~A Poliquin and R.~Tyrrell Rockafellar.
\newblock Second-order nonsmooth analysis in nonlinear programming.
\newblock {\em Recent advances in nonsmooth optimization}, pages 322--349,
  1995.

\bibitem{poliquin1996prox-regular}
Ren{\'e}~A Poliquin and R.~Tyrrell Rockafellar.
\newblock Prox-regular functions in variational analysis.
\newblock {\em Transactions of the American Mathematical Society},
  348(5):1805--1838, 1996.

\bibitem{powell1970numerical}
M.~Powell.
\newblock A hybrid method for nonlinear equations.
\newblock {\em Numerical Methods for Nonlinear Algebraic Equations}, pages
  87--144, 1970.

\bibitem{rockafellar1988first}
R.~Tyrrell Rockafellar.
\newblock First- and second-order epi-differentiability in nonlinear
  programming.
\newblock {\em Transactions of the American Mathematical Society},
  307(1):75--108, 1988.

\bibitem{rockafellar1989second}
R.~Tyrrell Rockafellar.
\newblock Second-order optimality conditions in nonlinear programming obtained
  by way of epi-derivatives.
\newblock {\em Mathematics of Operations Research}, 14(3):462--484, 1989.

\bibitem{rockafellar2011variational}
R.~Tyrrell Rockafellar and Roger~J.B. Wets.
\newblock {\em Variational analysis}, volume 317.
\newblock Springer, 2011.

\bibitem{stella2017forward}
Lorenzo Stella, Andreas Themelis, and Panagiotis Patrinos.
\newblock Forward--backward quasi-{N}ewton methods for nonsmooth optimization
  problems.
\newblock {\em Computational Optimization and Applications}, pages 1--45, 2017.

\bibitem{xu2012regularization}
Zongben Xu, Xiangyu Chang, Fengmin Xu, and Hai Zhang.
\newblock \(l_{1/2}\) regularization: a thresholding representation theory and
  a fast solver.
\newblock {\em IEEE Transactions on neural networks and learning systems},
  23(7):1013--1027, 2012.

\bibitem{zhang2004nonmonotone}
Hongchao Zhang and William~W. Hager.
\newblock A nonmonotone line search technique and its application to
  unconstrained optimization.
\newblock {\em SIAM Journal on Optimization}, 14(4):1043--1056, 2004.

\end{thebibliography}

\arxivfalse

\begin{appendix}
	\proofsection{sec:FBE}
		\label{sec:FBE_proof}
		\arxiv{%
	\begin{appendixproof}{prop:FBE1stOrder}
	\begin{proofitemize}
	\item\ref{prop:FBESubgrad}:
		from \cite[Ex. 10.32]{rockafellar2011variational} it follows that
		\(
			\partial g^\gamma
		{}\subseteq{}
			\frac 1\gamma(\id-\prox_{\gamma g})
		\).
		Moreover, since \(\gamma<\nicefrac{1}{L_f}\) and \(x\) is (strictly) twice differentiable at \(x\), matrix \(Q_\gamma(x)\) is well defined and nonsingular.
		In particular, from \cite[Thm. 10.49 and Ex. 9.25(c)]{rockafellar2011variational} it follows that
		\[
			\partial\left[g^\gamma(x-\gamma\nabla f(x))\right]
		{}\subseteq{}
			\tfrac 1\gamma Q_\gamma(x)
			\left(x-\gamma\nabla f(x)-\T(x)\right)
		{}={}
			\tfrac 1\gamma Q_\gamma(x)\Res(x)
			{}-{}
			Q_\gamma(x)\nabla f(x),
		\]
		and the claim follows from the expression \eqref{eq:FBE2} for \(\varphi_\gamma\).
	\item\ref{prop:FBEcriticDiff}:
		suppose that \(x\) is critical and let \(\gamma<\Gamma(x)\).
		Due to \cref{prop:SingleValuedFB} \(\Res(x)=\set{0}\), and \ref{prop:FBESubgrad} then implies that \(\partial\varphi_\gamma(x)\subseteq\set{0}\).
		The claim follows from \cite[Thm. 9.18(b)]{rockafellar2011variational}.
	\end{proofitemize}
	\end{appendixproof}
}%
\begin{appendixproof}{thm:2ndOrder}
\begin{proofitemize}
\item\ref{lem:JP}:
	It follows from \cite[Thm.s 3.8 and 4.1]{poliquin1996generalized} that \(\prox_{\gamma g}\) is (strictly) differentiable at \(x^\star-\gamma\nabla f(x^\star)\) iff \(g\) (strictly) satisfies \cref{ass:g2}.
	Consequently, if \(f\) is of class \(C^2\) around \(x^\star\) (and in particular strictly differentiable at \(x^\star\) \cite[Cor. 9.19]{rockafellar2011variational}), \(\Res(x) = x-\FB x\) is (strictly) differentiable at \(x^\star\) with Jacobian as in \eqref{eq:JacR} due to the chain rule of differentiation (and the fact that strict differentiability is preserved by composition).
	For \(\gamma'\in(\gamma,\Gamma(x^\star))\) and \(w\in\R^n\) we have
	\[
		\mathtight
		\twiceepi[-\nabla f(x^\star)]{g}{x^\star}[w]
	{}={}
		\liminf_{
			\substack{
				w'\to\fillwidthof[l]{0^+}{w}\\
				\fillwidthof[r]{w'}{\tau}\to 0^+
			}
		}{
			\frac{
				g(x^\star+\tau w')
				{}-{}
				g(x^\star)
				{}+{}
				\tau
				\innprod{\nabla f(x^\star)}{w}
			}{
				\nicefrac{\tau^2}{2}
			}
		}
	~{}\overrel[\geq]{\eqref{eq:ProxGrad'}}{}
		-\tfrac{1}{\gamma'}\|w\|^2.
	\]
	The expression \eqref{eq:GenQuadSecondEpiDer} of the second-order epi-derivative then implies
	\(
		\innprod{Mw}{w}
	{}\geq{}
		-\frac{1}{\gamma'}\|w\|^2
	\)
	for all \(w\in\R^n\) (since \(Mw=0\) for \(w\in S^\perp\)).
	Therefore,
	\(
		\lambda_{\rm min}(M)
	{}\geq{}
		-\nicefrac{1}{\gamma'}
	{}>{}
		-\nicefrac1\gamma
	\),
	proving \(I+\gamma M\) to be positive definite, and in particular invertible.
\proofswitch{
	To obtain an expression for \(P_\gamma(x^\star) = \jac{\prox_{\gamma g}}{x^\star-\gamma\nabla f(x^\star)}\) we can apply \cite[Ex. 13.45]{rockafellar2011variational} to the \emph{tilted} function \(g+\innprod{\nabla f(x^\star)}{{}\cdot{}}\) so that, letting
	\(
		\twiceepi g{}=\twiceepi[-\nabla f(x^\star)]g{x^\star}[{}\cdot{}]
	\)
	and \(\proj_S\) the idempotent and symmetric projection matrix on \(S\),
	\begin{align}
	\nonumber
		P_\gamma(x^\star)d
	{}={}&
		\prox_{(\gamma/2)\twiceepi g{}}(d)
	\\
	\nonumber
	{}={}&
		\argmin_{d'\in S}{
			\set{
				\tfrac12\innprod{d'}{Md'} + \tfrac{1}{2\gamma}\|d' - d\|^2
			}
		}
	\\
	\nonumber
	{}={}&
		\proj_S
		\argmin_{d'\in\R^n}{
			\set{
				\tfrac12\innprod{\proj_S d'}{M \proj_S d'} + \tfrac{1}{2\gamma}\|\proj_S d' - d\|^2
			}
		}
	\\
	\nonumber
	{}={}&
		\proj_S
		{\bigl(
			\proj_S[I+\gamma M]\proj_S
		\bigr)}^\dagger
		\proj_S d
	\\
	\label{eq:JacProx}
	{}={}&
		\proj_S
		[I+\gamma M]^{-1}
		\proj_S
	\end{align}
	where \(^\dagger\) indicates the pseudo-inverse, and last equality is due to \cite[Facts 6.4.12(i)-(ii) and 6.1.6(xxxii)]{bernstein2009matrix}.
}{
	We may now trace the proof of \cite[Lem. 2.9]{stella2017forward} to infer that
	\(
		JP_\gamma(x^\star)
	{}={}
		\proj_S
		[I+\gamma M]^{-1}
		\proj_S
	\).
}
	Apparently, \(JP_\gamma(x^\star)\) is symmetric and positive semidefinite.
\item\ref{lem:RStrictDiff}:
	\proofswitch{%
		With basic calculus rules it can be easily verified that, since \(\Res(x^\star)=0\), \(\nabla\varphi_\gamma=Q_\gamma \Res\) is (strictly) differentiable at \(x^\star\) provided that \(Q_\gamma\) is (strictly) continuous at \(x^\star\) and \(\Res\) is (strictly) differentiable at \(x^\star\).
	}{%
		Since \(Q_\gamma\) is (strictly) continuous at \(x^\star\) and \(\Res\) is (strictly) differentiable at \(x^\star\), from \cite[Lem. 6.2]{stella2017forward} we have that \(\nabla\varphi_\gamma=Q_\gamma \Res\) is (strictly) differentiable at \(x^\star\), and \eqref{eq:JacR} follows by the chain rule.
	}%
\item\ref{thm:FBEHess}:
	A simple application of the chain rule proves \eqref{eq:FBEHess}; moreover, combined with \eqref{eq:JacR} we obtain
	\(
		\nabla^2\varphi_\gamma(x^\star)
	{}={}
		\tfrac1\gamma
		\left[
			Q_\gamma(x^\star)
			{}-{}
			Q_\gamma(x^\star)P_\gamma(x^\star)Q_\gamma(x^\star)
		\right]
	\),
	and since both \(Q_\gamma(x^\star)\) and \(P_\gamma(x^\star)\) are symmetric, then so is \(\nabla^2\varphi(x^\star)\).
\end{proofitemize}
\end{appendixproof}

\begin{appendixproof}{thm:StrongMinimality}
	We will show that all conditions are equivalent to either one of the following
	\begin{enumerate}[%
		start=6,
		label={\textit{(\alph*)}},
		ref={\protect\OLDref*{thm:StrongMinimality}\textit{(\alph*)}},
	]
	\item\label{thm:2epiderPositive}\label{thm:2epiDef+}%
		\textit{\(
			\tinnprod{d}{(\nabla^2 f(x^\star)+M)d}>0
		\)
		\(\forall d\in S\), where \(M\) and \(S\) are as in \cref{ass:fg2};}%
	\item\label{thm:JRposdef}%
		\textit{\(J\Res(x^\star)\) is similar to a symmetric and positive definite matrix.}%
	\end{enumerate}
	\begin{proofitemize}
	\item\ref{thm:HessDef+} \(\Leftrightarrow\) \ref{thm:StrongMinimFBE}:
		trivial, since \(\nabla^2\varphi_\gamma(x^\star)\) exists as shown in \cref{thm:FBEHess}.
	
	\item\ref{thm:StrongMinimPsi} \(\Leftrightarrow\) \ref{thm:2epiDef+}:
		follows from \cite[Thm. 13.24(c)]{rockafellar2011variational}, since
		\[
			\twiceepi[0]{\varphi}{x^\star}[d]
		{}={}
			\tinnprod{d}{\nabla^2 f(x^\star)d}
			{}+{}
			\twiceepi[-\nabla f(x^\star)]{g}{x^\star}[d]
		{}={}
			\tinnprod{d}{(\nabla^2 f(x^\star)+M)d} + \delta_S(d).
		\]

	\item\ref{thm:HessDef+} \(\Leftrightarrow\) \ref{thm:LocMinFBE}:
		if \(\nabla^2\varphi_\gamma(x^\star)\succ 0\), then \(x^\star\) is a (strong) local minimum for \(\varphi_\gamma\) and, due to \eqref{eq:FBEHess}, necessarily \(J\Res(x^\star)\) is invertible.
		Conversely, if \(x^\star\) is a local minimum for \(\varphi_\gamma\), then \(\nabla^2\varphi_\gamma(x^\star)\succeq 0\).
		If, additionally, \(J\Res(x^\star)\) is invertible, then due to \eqref{eq:FBEHess} \(\nabla^2\varphi_\gamma(x^\star)\) is also invertible, and therefore positive definite.

	\item\ref{thm:HessDef+} \(\Leftrightarrow\) \ref{thm:JRposdef}:
		by comparing \eqref{eq:JacR} and \eqref{eq:FBEHess} we observe that \(J\Res(x^\star)\) is similar to the (symmetric) matrix
		\(
			Q_\gamma(x^\star)^{-\nicefrac12}
			\nabla^2\varphi_\gamma(x^\star)
			Q_\gamma(x^\star)^{-\nicefrac12}
		\),
		which is positive definite iff \(\nabla^2\varphi_\gamma(x^\star)\) is.
	
	\item\ref{thm:2epiDef+} \(\Leftrightarrow\) \ref{thm:JRposdef}:
		\proofswitch{
			from the proof of the implication above we already know that \(\jac{\Res}{x^\star}\) is similar to a symmetric matrix; therefore, \ref{thm:JRposdef} is equivalent to the condition \(\lambda_{\rm min}(J\Res(x^\star))>0\).
			Let \(P = P_\gamma(x^\star)\) and \(Q = Q_\gamma(x^\star)\), so that \(J\Res(x^\star)=\frac1\gamma(I-QP)\).
			From \cite[Thm. 7.7.3(a)]{horn2012matrix} it follows that \(\lambda_{\rm min}(J\Res(x^\star)) > 0\) iff \(Q^{-1} \succ P\).
			For all \(d\in S\), using \eqref{eq:JacProx} we have
			\begin{align*}
				\innprod{d}{(Q^{-1}-P)d}
			{}={} &
				\innprod{d}{Q^{-1}d}
				{}-{}
				\innprod{d}{
					\proj_S
					[I+\gamma M]^{-1}
					\proj_Sd
				}
			\\
			{}={} &
				\innprod{d}{Q^{-1}d}
				{}-{}
				\innprod{\proj_Sd}{
					[I+\gamma M]^{-1}
					\proj_Sd
				}
			\\
			{}={} &
				\innprod{d}{Q^{-1}d}
				{}-{}
				\innprod{d}{
					[I+\gamma M]^{-1}
					d
				}
			\\
			{}={} &
				\innprod{d}{Q^{-1}d}
				{}-{}
				\innprod{d}{
					[I+\gamma M]^{-1}
					d
				}
			\end{align*}
			and last quantity is positive iff \(I+\gamma M\succ Q\) on \(S\).
			By definition of \(Q\), we then have that this holds iff
			\(
				\nabla^2f(x^\star)+M\succ 0
			\)
			on \(S\), which is \ref{thm:2epiDef+}.
			For \(d\in S^\bot\) the computation trivializes to
			\(
				\innprod{d}{(Q^{-1}-P)d}
			{}={}
				\innprod{d}{Q^{-1}d}
			{}>{}
				0
			\)
			regardless.
		}{
			the proof is the same as that of \cite[Thm. 2.11\textit{(b)}\(\Leftrightarrow\)\textit{(c)}]{stella2017forward}.
		}%
	
%
	\item\ref{thm:LocMin} \(\Rightarrow\) \ref{thm:JRposdef}:
		with similar reasonings as in the proof of the implications ``\ref{thm:StrongMinimPsi} \(\Leftrightarrow\) \ref{thm:2epiDef+} \(\Leftrightarrow\) \ref{thm:JRposdef}'', we conclude that local minimality of \(x^\star\) for \(\varphi\) entails \(J\Res(x^\star)\) being similar to a symmetric and positive \emph{semi}definite matrix.
		Therefore, if \(J\Res(x^\star)\) is nonsingular, then it is similar to a symmetric and positive definite matrix.
	\item\ref{thm:LocMinFBE} \(\Rightarrow\) \ref{thm:LocMin}:
		trivial, since \(\varphi_\gamma\leq\varphi\) and \(\varphi_\gamma(x^\star)=\varphi(x^\star)\) (cf. \cref{prop:FBEleq,prop:FBEfix}).
	\end{proofitemize}
\end{appendixproof}

	\section{Additional results for \texorpdfstring{\Cref{sec:Algorithm}}{\S\ref{sec:Algorithm}}}
		\label{sec:Algorithm_proof}
\begin{appendixlem}\label{lem:IteratesDifference}%
Consider the iterates generated by \refZ[] and suppose that the directions \(\seq{d^k}\) are selected so as to satisfy \eqref{eq:dr}.
Then,
\begin{enumerate}
	\item\label{lem:Deltaxr}%
		\(
			\|x^{k+1} - x^k\|
		{}\leq{}
			(\gamma+D)\|r^k\|
		\)
	\item\label{lem:DeltaBarxr}%
		\(
			\|\bar x^{k+1} - \bar x^k\|
		{}\leq{}
			\gamma\|r^{k+1}\|
			{}+{}
			(2\gamma+D)
			\|r^k\|
		\)
	\item\label{lem:DifferenceTo0}%
		in particular, \(\|x^{k+1}-x^k\|\) and \(\|\bar x^{k+1} - \bar x^k\|\) converge to 0.
\end{enumerate}
\begin{proof}
For all \(k\) we have
\[
	\|x^{k+1}-x^k\|
{}={}
	\|\bar x^k + \tau_kd^k - x^k\|
{}={}
	\|\tau_kd^k-\gamma r^k\|
{}\leq{}
	\gamma\|r^k\|
	{}+{}
	\tau_k
	\|d^k\|
{}\leq{}
	(\gamma+D)
	\|r^k\|
\]
where in the last inequality we used the fact that \(\tau_k\in(0,1]\).
This proves \ref{lem:Deltaxr}, and \ref{lem:DeltaBarxr} trivially follows by the triangular inequality
\(
	\|\bar x^{k+1} - \bar x^k\|
	{}\leq{}
	\|x^{k+1} - x^k\|
	{}+{}
	\gamma\|r^{k+1}\|
	{}+{}
	\gamma\|r^k\|
\).
Using this, \ref{lem:DifferenceTo0} follows from \cref{thm:Critical:r}.
\end{proof}
\end{appendixlem}

\begin{lem}\label{lem:omega}%
Consider the iterates generated by \refZ[].
Suppose that \eqref{eq:dr} is satisfied and that the sequence \(\seq{x^k}\) is bounded.
Then, \(\omega(x^k)=\omega(\bar x^k)\) are nonempty compact and connected sets over which \(\varphi\) and \(\varphi_\gamma\) are constant and coincide.
Moreover,
\begin{equation}\label{eq:omegaDist}
	\lim_{k\to\infty}{\dist(x^k,\omega(x^k))}
{}={}
	\lim_{k\to\infty}{\dist(\bar x^k,\omega(x^k))}
{}={}
	0.
\end{equation}
\begin{proof}
The sets of cluster points are nonempty because of boundedness of the sequences; in turn, connectedness and compactness as well as \eqref{eq:omegaDist} are shown in \cite[Rem. 5]{bolte2014proximal}, which applies since \(\|x^{k+1}-x^k\|\) and \(\|\bar x^{k+1} - \bar x^k\|\) converge to \(0\) (cf. \cref{lem:DifferenceTo0}).
Moreover, since \(\seq{\varphi_\gamma(x^k)}\) converges to some value \(\varphi_\star\in\R\) and \(\omega(x^k)=\omega(\bar x^k)\subseteq\fix\T\) as shown in \cref{thm:Critical}, it follows \Cref{prop:FBEfix} that \(\varphi\) and \(\varphi_\gamma\) coincide on \(\omega(x^k)\) (and equal \(\varphi_\star\)).
\end{proof}
\end{lem}

\begin{lem}\label{lem:Lojasiewicz}%
Suppose that \Cref{ass:fg2} is satisfied at a strong local minimum \(x^\star\) of \(\varphi\).
Then, for any \(\gamma\in(0,\nicefrac{1}{L_f})\) the FBE \(\varphi_\gamma\) possesses the KL property at \(x^\star\), and the desingularizing function \(\psi\) can be taken of the form
\(
	\psi(s)
{}={}
	\rho s^{\nicefrac 12}
\)
for some \(\rho>0\).
\begin{proof}
	From \cref{thm:FBEHessPosDef} it follows that \(x^\star\) is a strong local minimum for \(\varphi_\gamma\) at which \(\varphi_\gamma\) is twice differentiable with \(H_\star\coloneqq\nabla^2\varphi_\gamma(x^\star)\succ 0\).
	Let
	\(
		\lambda
	{}\coloneqq{}
		\lambda_{\rm min}(H_\star)
	\)
	and
	\(
		\Lambda
	{}\coloneqq{}
		\lambda_{\rm max}(H_\star)
	\).
	Since \(\nabla\varphi_\gamma(x^\star)=0\), from a second-order expansion of \(\varphi_\gamma\) and a first-order expansion of \(\nabla\varphi_\gamma\) we obtain that there exists a neighborhood \(U_{x^\star}\) of \(x^\star\) such that, for all \(x\in U_{x^\star}\),
	\(
		\varphi_\gamma(x)-\varphi_\gamma(x^\star)
	{}\leq{}
		\tfrac\Lambda4
		\|x-x^\star\|^2
	\)
	and
	\(
		\|\nabla\varphi_\gamma(x)\|
	{}\geq{}
		\tfrac\lambda2
		\|x-x^\star\|
	\),
	and in particular
	\(
		\psi'
		\left(
			\varphi_\gamma(x)-\varphi_\gamma(x^\star)
		\right)
		\|\nabla\varphi_\gamma(x)\|
	{}={}
		\tfrac{\rho}{
			2\sqrt{
				\varphi_\gamma(x)-\varphi_\gamma(x^\star)
			}
		}
		\|\nabla\varphi_\gamma(x)\|
	{}\geq{}
		\tfrac{\rho\lambda}{
			2\sqrt{
				\Lambda
			}
		}
	\).
	Letting
	\(
		\rho
	{}={}
		\frac{2\sqrt{\Lambda}}{\lambda}
	\)
	we obtain that \(\psi\) is a KL function for \(\varphi_\gamma\) at \(x^\star\).
\end{proof}
\end{lem}

\end{appendix}

\end{document}